\documentclass[11pt]{article}
\usepackage[letterpaper,margin=1in]{geometry}
\usepackage[T1]{fontenc}
\usepackage{lmodern,amsmath,amssymb,amsthm,mathtools,microtype}
\usepackage{xcolor,tikz,caption,enumitem,needspace,booktabs,longtable}
\usepackage[section]{placeins}
\usetikzlibrary{arrows.meta,positioning,calc,fit,backgrounds}
\usepackage[colorlinks=true,linkcolor=blue!40!black,urlcolor=blue!40!black,citecolor=blue!40!black]{hyperref}
\hypersetup{pdftitle={Asymptotically Optimal Weighted Edge Colorings
via Local Bin Packing},pdfauthor={Rishikesh Gajjala and Arindam Khan}}
\usepackage{aliascnt}
\usepackage[capitalize,nameinlink,noabbrev,nosort]{cleveref}
\usepackage{tcolorbox}
\tcbuselibrary{breakable}

\setlist{itemsep=3pt,topsep=5pt}
\newtheorem{theorem}{Theorem}[section]
\newaliascnt{lemma}{theorem}
\newtheorem{lemma}[lemma]{Lemma}
\aliascntresetthe{lemma}
\newaliascnt{proposition}{theorem}
\newtheorem{proposition}[proposition]{Proposition}
\aliascntresetthe{proposition}
\theoremstyle{definition}\newtheorem{algorithm}{Algorithm}[section]
\theoremstyle{remark}
\crefname{lemma}{lemma}{lemmas}
\Crefname{lemma}{Lemma}{Lemmas}
\crefname{proposition}{proposition}{propositions}
\Crefname{proposition}{Proposition}{Propositions}
\crefname{algorithm}{algorithm}{algorithms}
\Crefname{algorithm}{Algorithm}{Algorithms}
\numberwithin{equation}{section}
\usepackage{thm-restate}
\newcommand{\restatedstatement}[2]{%
  \begingroup
  \let\originalstatement#1%
  \renewcommand{#1}[1][]{%
    \if\relax\detokenize{##1}\relax
      \originalstatement[restated]%
    \else
      \originalstatement[##1; restated]%
    \fi}%
  #2*%
  \endgroup}
\newcommand{\E}{\mathbb E}
\newcommand{\Prb}{\mathbb P}
\newcommand{\bp}{\operatorname{bp}}
\newcommand{\dw}{\Delta_w}
\newcommand{\cw}{\chi'_w}
\newcommand{\one}{\mathbf 1}
\newcommand{\cC}{\mathcal C}

\newcommand{\eps}{\varepsilon}

\definecolor{inkblue}{RGB}{34,87,128}
\definecolor{warmred}{RGB}{172,70,46}
\definecolor{softblue}{RGB}{228,238,246}
\definecolor{softgold}{RGB}{249,239,213}
\definecolor{softgreen}{RGB}{224,239,228}
\tikzset{vertex/.style={circle,draw=black,fill=white,inner sep=2.5pt},
 lab/.style={font=\small},
 panel/.style={draw=inkblue,rounded corners=3pt,fill=softblue,align=center,inner sep=7pt}}
\title{Asymptotically Optimal Weighted Edge Colorings \\via Local Bin Packing}
\author{Rishikesh Gajjala\\\small New York University Abu Dhabi
 \and Arindam Khan\\\small Indian Institute of Science (IISc)}
\date{}
\begin{document}
\pagenumbering{roman}
\hypersetup{pageanchor=false}
\begin{titlepage}
\maketitle
\thispagestyle{empty}
\begin{abstract}

We study \emph{weighted edge coloring}, where each edge of a graph
$G$ has a weight in $(0,1]$, and the goal is to assign colors so
that, at every vertex, the total weight of incident edges of each
color is at most one. The objective is to minimize the number of
colors, denoted by $\chi'_w(G)$. This problem generalizes classical
edge coloring (when all weights are one) and bin packing (when
$G$ is a star). Two natural lower bounds are the maximum weighted
degree $\Delta_w(G)$ and
the maximum local bin-packing number \(\mathrm{bp}(G)\), defined as the maximum, over all vertices, of the minimum number of unit-capacity bins needed to pack the incident edge weights at any vertex.
While
$\lceil\Delta_w(G)\rceil\leq\mathrm{bp}(G)\leq\chi'_w(G)$,
understanding how closely these local lower bounds approximate
the global optimum has been a longstanding open problem.

Our main result establishes an \emph{asymptotic weighted analogue
of Vizing's theorem}: every simple graph admits a weighted edge
coloring using at most $(1+o(1))\mathrm{bp}(G)$ colors. Thus,
local bin packings can be made globally consistent with an
asymptotically vanishing relative loss, even in nonbipartite
graphs. Our main technical contribution is a \emph{local-to-global
rounding framework} that coordinates feasible local packing
distributions into a global partial coloring while leaving only
$o(\mathrm{bp}(G))$ uncolored weight at every vertex.

For bipartite multigraphs, we prove
$\chi'_w(G)\le (2+o(1))\Delta_w(G)$, attaining the optimal
leading coefficient and establishing the Chung--Ross conjecture [SICOMP'91]  asymptotically. For arbitrary
loopless multigraphs, we obtain the exact bound
$\chi'_w(G)\le 3\lceil\Delta_w(G)\rceil$, with optimal leading
coefficient. These results provide weighted analogues of the
classical K\H{o}nig--Shannon bounds.

Going beyond weighted degree, we exploit the compatibility of
edge weights within local bins to prove
\[
    \chi'_w(G)\le
    \begin{cases}
       (1.62+o(1))\,\mathrm{bp}(G),
           & \text{for bipartite multigraphs},\\
       (1.999+o(1))\,\mathrm{bp}(G),
           & \text{for loopless multigraphs}.
    \end{cases}
\]
The bipartite bound combines harmonic bin-packing charges with
compatible selections of edge groups. The general bound additionally
handles dense odd vertex sets through contraction and restoration.
We complement these upper bounds with bipartite multigraphs satisfying $\chi'_w(G)=\frac43\,\mathrm{bp}(G)$, even when the
underlying simple graph is a tree of maximum degree three.

We also construct a simple bipartite counterexample to the
exact weighted-degree conjecture of Feige and Singh [ESA'08].

All five upper bounds have polynomial-time algorithms.
All five existence statements are formalized in Lean~4.

\end{abstract}
\newpage

\vfill
\begin{center}\textbf{Contents}\end{center}
The overview occupies pages~\pageref{sec:introduction}--\pageref{overview:end}.
The appendices contain the detailed upper-bound proofs and their algorithmic details.
Section~\ref{sec:lower-bounds} contains both complete lower-bound proofs.
\begin{center}
\small
\begin{tabular}{@{}p{.08\textwidth}p{.75\textwidth}r@{}}
\toprule
Section & Overview & Page\\
\midrule
\ref{sec:introduction} & \hyperref[sec:introduction]{Introduction} & \pageref{sec:introduction}\\[3pt]
\ref{sec:overview-three} & \hyperref[sec:overview-three]{The exact $3\lceil\Delta_w\rceil$ bound} & \pageref{sec:overview-three}\\[3pt]
\ref{sec:overview-simple} & \hyperref[sec:overview-simple]{Simple graphs and the weighted analogue of Vizing's theorem} & \pageref{sec:overview-simple}--\pageref{ov:simple-end}\\[3pt]
\ref{sec:overview-bipartite} & \hyperref[sec:overview-bipartite]{The tight $2\Delta_w$ coefficient for bipartite multigraphs} & \pageref{sec:overview-bipartite}--\pageref{ov:bipartite-end}\\[3pt]
\ref{sec:lower-bounds} & \hyperref[sec:lower-bounds]{Lower bounds and the Feige--Singh conjecture} & \pageref{sec:lower-bounds}--\pageref{ov:lower-end}\\[3pt]
\ref{sec:overview-constants} & \hyperref[sec:overview-constants]{Improved constants and the K\H{o}nig--Shannon comparison} & \pageref{sec:overview-constants}--\pageref{ov:constants-end}\\[3pt]
--- & \hyperref[sec:ai-overview]{AI Methodology} & \pageref{sec:ai-overview}\\
\bottomrule
\end{tabular}
\end{center}
\smallskip\noindent\textbf{Appendix guide}
\begin{center}
\small
\begin{tabular}{@{}p{.10\textwidth}p{.12\textwidth}p{.65\textwidth}r@{}}
\toprule
Appendix & Main text & Contents & Page\\
\midrule
\ref{sec:elementary} & --- & Preliminaries & \pageref{sec:elementary}\\[3pt]
\ref{sec:general-multi} & \ref{sec:overview-three} & The exact $3\lceil\Delta_w\rceil$ bound & \pageref{sec:general-multi}\\[3pt]
\ref{sec:round-overview} & \ref{sec:overview-simple} & Simple graphs: full rounding proof and implementation & \pageref{sec:round-overview}\\[3pt]
\ref{sec:bip-multi} & \ref{sec:overview-bipartite} & Bipartite multigraphs: the tight $2\Delta_w$ coefficient & \pageref{sec:bip-multi}\\[3pt]
\ref{sec:bipartite-bounds} & \ref{sec:overview-bipartite-bp} & Bipartite bin-packing bounds, through $1.62$ & \pageref{sec:bipartite-bounds}\\[3pt]
\ref{sec:general-bounds} & \ref{ov:page-12} & General bin-packing bounds, through $1.999$ & \pageref{sec:general-bounds}\\[3pt]
\bottomrule
\end{tabular}
\end{center}
\noindent The appendix starts on page~\pageref{sec:appendix-guide}.
\end{titlepage}
\hypersetup{pageanchor=true}
\pagenumbering{arabic}

\section{Introduction}\label{sec:introduction}\label{ov:page-1}

Edge coloring is a classical and fundamental topic in graph theory and algorithms. We study the {\em weighted edge coloring} problem, where we are given an edge-weighted graph $G:=(V,E)$, possibly with parallel edges but without loops,  and a weight function $w: E \rightarrow (0,1]$.
A coloring is a \emph{proper weighted coloring} if, at every vertex, the sum
of the weights of incident edges of each color is at most one. The objective
is to minimize the number of colors, denoted by $\cw(G)$.

This problem is a generalization of classical edge coloring in unweighted graphs (when all edge weights are one) and bin packing (when the graph is a star). Apart from its structural importance, the problem also finds many applications in switching networks, scheduling, frequency assignment, and graph decomposition.
A principal application is multirate routing in three-stage
Clos networks~\cite{Clos53,MT89,CR91}. Here, vertices represent input and
output switches, edges represent unsplittable connection requests, and
colors represent middle-stage switches. Minimizing the number of colors
reduces the interconnection resources needed to support simultaneous
requests in communication and parallel-computing systems. We refer to
Correa and Goemans~\cite{CG07} for a detailed account of this connection.

In (unweighted) edge coloring, a color class forms a matching.  Three classical theorems describe much of the landscape.
K\H{o}nig's theorem gives $\chi'(G)=\Delta(G)$ for bipartite
multigraphs~\cite{Konig16}; Vizing's theorem gives
$\chi'(G)\le\Delta(G)+1$ for simple graphs~\cite{Vizing64}; and Shannon's
theorem gives $\chi'(G)\le\lfloor3\Delta(G)/2\rfloor$ for arbitrary
multigraphs~\cite{Shannon49}. Here $\Delta(G)$ is the maximum vertex degree.

In the weighted setting, two natural local lower bounds are
$\dw(G):=\max_{v\in V}\sum_{e\ni v}w_e$, and $\bp(G):=\max_{v\in V}\bp_v(G)$,
where $\bp_v(G)$ is the minimum number of unit-capacity bins needed to pack
the weights of the edges incident with $v$.
Every proper weighted coloring induces a feasible
bin packing at each vertex, so
$\lceil\dw(G)\rceil\le\bp(G)\le\cw(G)$.
We call the ratio $\cw(G)/\bp(G)$, when $\bp(G)>0$, the
\emph{global-to-local packing gap}: it measures the additional cost of
reconciling locally feasible packings into a single coloring. Even on a
star, the distinction between volume and local packability matters.
A star with $k$ edges of weight $1/2+\eps$ has $\cw(G)=\bp(G)=k$ and
$\dw(G)=k(1/2+\eps)$. In particular, $\cw(G)>2\dw(G)-1$ whenever
$0<\eps<1/(2k)$.

In 1991, Chung and Ross \cite{CR91} conjectured that
every  bipartite
multigraph admits a proper weighted coloring
 using $2 \bp(G) -1$ colors.\footnote{Their original conjecture was formulated
in the context of Clos networks; this graph-theoretic formulation is
explicit in later work~\cite{FS,KS15}.} Feige and Singh~\cite{FS} also proposed a related conjecture that a proper weighted coloring exists using at most $2N-1$ colors, where $N$ is the smallest integer strictly greater
than $\dw(G)$, namely $N=\lfloor\dw(G)\rfloor+1$.
These conjectures ask how efficiently locally feasible
packings can be reconciled into a globally consistent coloring.
The conjectures have motivated a sequence of works in approximation algorithms.

Du, Gao, Hwang, and Kim~\cite{DGHK98} obtained the bound
$\cw(G) \le \frac{41}{16}\bp(G)+O(1)$ for bipartite multigraphs. Correa and Goemans~\cite{CG07} improved this to $2.5480\bp(G)+o(\bp(G))$ and separately established
$2.557\dw(G)+o(\dw(G))$.
Feige and Singh~\cite{FS} subsequently gave a combinatorial algorithm using $\lceil\frac94\dw(G)\rceil$ colors. Khan and Singh~\cite{KS15} gave a guarantee of $\lceil2.2223\bp(G)\rceil$ colors.
Ngo and Vu~\cite{NV03} constructed bipartite instances exhibiting a
$5/4$ gap between $\cw(G)$ and $\bp(G)$.

For arbitrary multigraphs, Feige and Singh~\cite{FS} extended their
algorithm to use at most $\lceil\frac94\dw(G)\rceil+\mu(G)$ colors,
where $\mu(G)$ is the maximum edge multiplicity.
Odd cycles and multiplicities create genuine
obstructions. The triangle example recorded by Feige and Singh~\cite[Section~1.2]{FS}
has $t$ parallel edges on each side, all of weight
$1/2+\varepsilon$. It has $\cw(G)=3t$, $\bp(G)=2t$, and
$\dw(G)=t(1+2\varepsilon)$: its global-to-local packing gap is $3/2$, while its ratio to weighted degree approaches three.
Improved bounds are also known
for certain special graph classes; see Huc~\cite{Huc11} and
Sannyasi~\cite{Sannyasi20}.

Despite this progress, even for simple bipartite graphs, it remained unknown whether the local bin-packing bound approximates the global optimum asymptotically, or whether the factor two in terms of maximum weighted degree is sufficient asymptotically. More broadly, weighted analogues of Vizing's and Shannon's theorems remained elusive.

\subsection{Our contributions}

We establish weighted analogues of Vizing's, K\H{o}nig's,
and Shannon's theorems, resolve the Chung--Ross conjecture
asymptotically, and disprove the exact weighted-degree
conjecture of Feige and Singh.
Our main results are summarized below.

\begin{theorem}\label{thm:main}
The following bounds hold uniformly over the indicated graph classes.
\begin{center}\normalfont
\begin{tabular}{@{}lll@{}}
\toprule
Graph class & Upper bound for $\chi'_w(G)$ & Full proof\\
\midrule
Simple graphs & $(1+o(1))\bp(G)$ & \Cref{sec:round-overview}\\[3pt]
Bipartite multigraphs & $(2+o(1))\Delta_w(G)$ & \Cref{sec:bip-multi}\\[3pt]
Loopless multigraphs & $3\lceil\Delta_w(G)\rceil$ & \Cref{sec:general-multi}\\[3pt]
Bipartite multigraphs & $(1.62+o(1))\bp(G)$ & \Cref{sec:bipartite-bounds}\\[3pt]
Loopless multigraphs & $(1.999+o(1))\bp(G)$ & \Cref{sec:general-bounds}\\
\bottomrule
\end{tabular}
\end{center}
The coefficient one for simple graphs and the total-weight coefficients two and three are optimal.\footnote{An asymptotic bound $(c+o(1))n$ means that, for every $\varepsilon>0$, there is $C_\varepsilon$ with $\chi'_w(G)\le(c+\varepsilon)n+C_\varepsilon$. The constant is independent of the number of vertices in the graph, multiplicity of edges, and the edge weights.}
\end{theorem}

In the following, we explain the results in more detail.

\paragraph{Simple graphs: an asymptotic weighted
analogue of Vizing's theorem.}
Our {\em main} result shows that
$\cw(G)\le\bp(G)+o(\bp(G))$ for every simple graph,
whether or not it is bipartite.
This (asymptotically) extends the  coefficient-one phenomenon
of Vizing's theorem to weighted  graphs.

The local packing problems here are considerably richer
than in the unweighted setting.
When every edge has weight one, a local bin can contain
only one edge.
With arbitrary weights, an optimal local packing may
rely on intricate combinations of different item sizes.
Moreover, these packings cannot be chosen independently:
an edge belongs to the packing problems at both of its
endpoints and must receive the \emph{same} color at both.
{\em Surprisingly}, our result shows that these local packings can nevertheless
be made globally consistent with an asymptotically
vanishing relative loss.

The main technical contribution is a local-to-global
rounding framework.
Our rounding theorem coordinates the local choices
while retaining almost all incident weight at every
vertex.
To coordinate local packing distributions, we use partial
resampling~\cite{HS,MT} and hypergraph matching~\cite{DP}.
The proof has four main steps.
First, each local bin is represented by a short
description consisting of a bounded number of slots
for large edges and reserved space for small edges.
Second, fractional assignments link compatible slots
at the two endpoints of each edge.
Third, we choose one description for each vertex
and color, controlling the resulting local capacities.
Finally, hypergraph matching selects whole large edges,
and a separate rounding procedure assigns the small edges.
A small fresh palette colors the remaining edges.
A key feature is that the loss at every vertex is
$o(\bp(G))$, independently of the number of vertices.
See \Cref{sec:overview-simple} for details.

\paragraph{Tight weighted degree bounds and the Chung--Ross conjecture.}
For bipartite multigraphs we obtain an asymptotically tight guarantee:
$\cw(G)\le(2+o(1))\dw(G)$.
In particular, since $\dw(G)\le\bp(G)$, this establishes the Chung--Ross conjecture asymptotically.

The proof groups parallel edges into feasible bundles, leaving a simple graph of light remainders. K\H{o}nig's theorem properly colors the heavier bundles with approximately $2\dw(G)$ colors. The challenge is that a remaining edge can use only a color missing at \emph{both} endpoints, and the two sets of missing colors may initially have little overlap. Kempe-style exchanges along two-color paths coordinate these missing sets in a weighted average sense, after which the simple-graph local-to-global rounding theorem completes the coloring. Thus the bipartite bound is also a substantial application of our main rounding framework.
See \Cref{sec:overview-bipartite} for details.

For arbitrary loopless multigraphs we obtain the tight bound
$\cw(G)\le3\lceil\dw(G)\rceil$. This proof is considerably simpler and independent of the rounding machinery.
We present a self-contained argument for this guarantee as a gentle starting point in \Cref{sec:overview-three}.

\vspace{-3mm}
\paragraph{Improved bounds from the contents of bins.}
Exploiting local bin packing, we obtain further improvements on upper bounds on $\cw(G)$:
$(1.62+o(1))\bp(G)$ for bipartite multigraphs, and $(1.999+o(1))\bp(G)$
for arbitrary loopless multigraphs.

The bipartite bound builds on our rounding framework. We first use harmonic rounding, ubiquitous in bin packing~\cite{LeeLee85}, to achieve a guarantee of the classical harmonic constant $h_\infty\approx1.691$.
The further improvement to $1.62$ uses three selections of edges
from compatible weight ranges: a heavy edge and a small edge, or two
small edges, can share one color. Coloring these selections separately
reduces the harmonic charge of the remainder.

For general multigraphs, packet coloring must additionally respect \emph{odd-set constraints}. Degree bounds from the previous results already handle  sufficiently large odd sets, so the difficult cases reduce to dense triangles and five-vertex sets. We show that the relevant obstructions can be selected disjointly and contracted. After coloring the remaining edges, we restore their internal edges without conflicting with external colors.
For multigraphs, we also employ
fractional edge-coloring methods~\cite{Edmonds65,SS}.
Finally, this gives a coefficient strictly below two even without bipartiteness; the details are explained in \Cref{sec:overview-constants}.

\vspace{-3mm}
\paragraph{Lower bounds and exact conjectures.}
We give a \emph{simple bipartite} counterexample (see \Cref{sec:counterexample}) to the exact weighted-degree conjecture of Feige and Singh \cite{FS}: it has $\dw(G)=39/20<2$ but $\cw(G)=4$, whereas the conjecture predicts three colors.

Parallel edges make the simple-graph theorem genuinely special. We construct bipartite multigraphs $G_N$ with
$\bp(G_N)=3N$, $\cw(G_N)=4N$,
although the underlying simple graph is a tree of maximum degree three (see \Cref{subsec:4by3lower}). This improves the earlier $5/4$ lower bound of \cite{NV03} to $4/3$.
Our bounds also show that the ratio of the optimal uniform
asymptotic coefficients for general and bipartite multigraphs, relative
to $\bp(G)$, is strictly less than $3/2$:
$\frac{1.999}{4/3}<\frac32$.

\vspace{-3mm}
\paragraph{Algorithms and verification.}
All five upper bounds in \Cref{thm:main} have polynomial-time algorithms.
Their descriptions and running-time analyses appear 
in the appendices. Companion Lean~4 formalizations cover the five
\emph{existence} statements in the table, not the running-time analyses.

\vspace{-2mm}
\subsection{Related results.}
Our algorithms draw on classical bin-packing ideas, including linear
grouping~\cite{fernandez1981bin}, configuration-LP and
rounding~\cite{KK82, HR17}. They also use techniques from online bin packing related to Harmonic grouping ~\cite{LeeLee85, Ramanan89, Seiden02}. For recent developments on online bin packing, see ~\cite{BaloghBDEL18, BaloghBDEL21, HebbarKS24, AyyadevaraD0S25, DJKS26}.
The online variant of weighted edge coloring, where edges
must be colored irrevocably upon arrival, has also been
studied~\cite{CG07,EpsteinWEC23}.

\section{A Weighted Analogue of Shannon's Theorem}
\label{sec:overview-three}\label{ov:page-3}

We begin with an elementary argument showing that every
loopless weighted multigraph admits a coloring with
$\cw(G)\le 3\lceil\dw(G)\rceil$.
Unlike our subsequent results, this
argument requires no probabilistic rounding.
Here, the idea is to properly color a carefully selected subgraph and then insert all remaining edges into the same palette.
It is similar in spirit to the
previous approaches of Feige and Singh~\cite{FS}. See \Cref{sec:general-multi} for omitted details in this section. 

\paragraph{Select and color a bounded-degree skeleton.}
Put $n=\lceil\dw(G)\rceil$ and inspect the edges in nonincreasing order of
weight. Starting with $F=\varnothing$, keep an edge if adding it leaves the
degree of $F$ at most $2n$ at both endpoints. Thus $\Delta(F)\le 2n$, and
Shannon's theorem properly colors $F$ with at most $3n$ colors. Pad the
palette with unused colors to obtain exactly $3n$ colors. A constructive
proof of the required form of Shannon's theorem appears in
\Cref{lem:three-per-two}.

Consider a rejected edge $e=uv$ of weight $a$. When $e$ was inspected, one
endpoint, say $u$, already had $2n$ selected edges. Each was inspected
before $e$ and therefore has weight at least $a$. We call these edges the
\emph{witnesses} of $e$.
The witnesses serve two purposes. First, $a<1/2$: otherwise the $2n$
witnesses together with $e$ would have weight at least $(2n+1)/2>n$ at $u$,
contradicting $\dw(G)\le n$. Second, because the coloring of $F$ is proper,
the witnesses receive $2n$ distinct colors. We never recolor $F$, so each of
these colors continues to carry load at least $a$ at $u$, regardless of
which other edges are inserted.
The following elementary lemma explains why the witnesses suffice.

\begin{restatable}{lemma}{commonroomlemma}
\label{lem:overview-common-room}\label{lem:three-common-room}
Let $n\ge 1$ be an integer and $0<a<1/2$. Suppose the loads $x_c,y_c\ge 0$
of $3n$ colors at two vertices satisfy
$\sum_c x_c\le n$, and $\sum_c y_c\le n$.
If at least $2n$ colors satisfy $x_c\ge a$, then some color $c$ satisfies
$x_c+a\le 1$ and $y_c+a\le 1$.
\end{restatable}

\noindent\emph{Proof sketch.}
Suppose, to the contrary, that no color has room at both endpoints.
Then every color has load greater than $1-a$ at one endpoint. Let $S=\{c:y_c>1-a\}$, and $s=|S|$.
Since $\sum_c y_c\le n$, we have $s<n/(1-a)$; for every $c\notin S$,
we must have $x_c>1-a$.

If $s\le n$, at least $2n$ colors lie outside $S$, giving
$\sum_c x_c>(3n-s)(1-a)\ge 2n(1-a)>n$,
a contradiction. Otherwise $s>n$. At most $n$ colors have $x_c<a$, so at
least $s-n$ colors in $S$ have $x_c\ge a$. Consequently,
$\sum_c x_c
    >(3n-s)(1-a)+(s-n)a >n+\frac{n(1-2a)^2}{1-a}>n$,
where the second inequality uses $s<n/(1-a)$ and $a<1/2$.
This is again a contradiction. The full counting argument is given in
\Cref{sec:common-room-proof}.

\paragraph{Insert the remaining edges.}
Process these edges in any order. Consider $e=uv$ of weight $a$,
with witness endpoint $u$, and let $x_c,y_c$ be the current color loads
at $u,v$. As only edges of $G$ have been colored,
$\sum_c x_c,\sum_c y_c\le\dw(G)\le n$. Moreover, the $2n$ distinct witness
colors still have load at least $a$ at $u$. The lemma thus provides a
color with room for $e$ at both endpoints; assign that color to $e$.
The witnesses remain intact throughout the process, so the same argument
applies to every rejected edge. Thus all edges receive one of the original
$3n$ colors, proving
$\cw(G)\le 3\lceil\dw(G)\rceil$.
The algorithm is deterministic and polynomial-time: sort the edges,
construct and properly color $F$, and scan the palette when inserting each
remaining edge. See \Cref{sec:general-total-proof} for the details. 

\paragraph{The coefficient three is best possible.}
Use the triangle construction recorded by Feige and Singh~\cite[Section~1.2]{FS}:
take $t$ parallel edges on each side, each of weight
$(t+1)/(2t)>1/2$. Any two edges meet and cannot share a color. Hence, $\cw(G)=3t, \dw(G)=t+1$, and  $\bp(G)=2t$.
Therefore $\cw(G)/\dw(G)\to 3$ as $t\to\infty$. The same construction
gives the lower bound $3/2$ relative to local bin packing; see \Cref{app:triangle-lower} for details and 
\Cref{fig:general-triangle} in the appendix for an illustration.

\section{Simple Graphs: Local Packings Almost Suffice}
\label{sec:overview-simple}\label{ov:page-4}

We show that the local bin-packing bound is asymptotically
sufficient for every simple graph, even if the graph is
not bipartite: $\chi'_w(G)\le \bp(G)+o(\bp(G))$.
Since $\bp(G)\le \chi'_w(G)$, this gives an
asymptotically optimal approximation to the minimum
number of colors.

The main difficulty is easy to state. At each vertex,
the incident edges can be packed into
$n:=\bp(G)$ bins. However, these are separate local
packings: an edge $uv$ may belong to one bin at $u$
and a different bin at $v$. A valid coloring must
assign a \emph{single} color to each edge, feasible at
both endpoints. Our goal is to reconcile these local
packings while using only $o(n)$ additional colors.

The key observation is that we need not preserve the
original bins. An edge can be assigned to any color
with sufficient capacity at both endpoints. Moreover,
we may leave a small amount of incident weight
uncolored and handle those edges using fresh colors
at the end. Our approach therefore separates the
problem into four steps: describe local packings compactly,
couple the possible placements of each edge at its two endpoints,
choose local profiles with spare capacity, and finally select
whole edges without violating any capacity.
We will show that each step incurs only a small loss. For the omitted details, see \Cref{sec:round-overview}.

\paragraph{Local-to-global rounding.}
Imagine that, for each color, every vertex proposes
a random subset of incident edges that fits into
one bin. Such proposals are feasible locally, but
the two endpoints of an edge may disagree about
whether that edge should receive the color.
Our main rounding theorem resolves this disagreement:
provided that no edge places too much probability
on any single color, it produces a globally
consistent partial coloring that preserves almost
all of the proposed incident weight at every vertex.
To formalize local feasibility, let
$W_G(v):=\sum_{e\ni v}w_e$ and define
\[
\mathcal K_v:=
\operatorname{conv}\left\{
\mathbf 1_S:
S\subseteq E(v),\
\sum_{e\in S}w_e\le 1
\right\}.
\]
Thus, $\mathcal K_v$ is the convex hull of the
incidence vectors of all feasible bins at $v$.
Equivalently, membership in $\mathcal K_v$ has a
configuration-LP interpretation. Let
$\mathcal F_v$ be the family of feasible subsets
of $E(v)$. A vector $x$ belongs to $\mathcal K_v$
if and only if there are nonnegative values
$\lambda_S$, for $S\in\mathcal F_v$, such that
$\sum_{S\in\mathcal F_v}\lambda_S=1$, and 
$x_e=\sum_{{S\in\mathcal F_v, e\in S}}
\lambda_S$, where $e\in E(v)$.
In other words, $\lambda$ is a probability
distribution over feasible integral bins, and
$x_e$ is the probability that the sampled bin
contains edge $e$. This is stronger than merely
requiring the expected load
$\sum_{e\ni v}w_ex_e$ to be at most one.

Independent choices need not agree across endpoints. Fix at each vertex $v$ a
packing $B_{v,1},\ldots,B_{v,n}$ into $n$ bins,
padding with empty bins if necessary.
Choosing one of these bins uniformly at random
includes each incident edge with probability
$1/n$. This suggests proposing $x_{e,c}=1/n$
for every $e\in E$ and $c\in[n]$.
Every edge then has total proposed probability
$\sum_c x_{e,c}=1$, and every local vector
$(x_{e,c}:e\ni v)$ belongs to $\mathcal K_v$,
since it is the average of $n$ feasible bins.

But suppose that, for each color $c$, the two
endpoints of an edge $e=uv$ independently choose
their local bins. Each endpoint includes $e$
with probability $1/n$, so both include it
with probability only $1/n^2$. Even across
all $n$ colors, the expected number of colors
on which they agree to include $e$ is just $1/n$.
Thus, independent local choices lose almost
all coloring opportunities.

The essential challenge is therefore to
\emph{coordinate the choices at the two endpoints}
before assigning actual colors. Our rounding
theorem accomplishes this while leaving only
an arbitrarily small multiple of $n$ in uncolored
weight at every vertex. 
This per-vertex guarantee is crucial,
since a bound on the total uncolored weight
would not prevent the remaining edges from
concentrating at a single vertex.

\begin{restatable}[Local-to-global rounding]{theorem}{localglobalrounding}
\label{thm:overview-rounding}\label{thm:round}
Fix constants $C_0,C_1,\varepsilon>0$.
Let $G$ be a simple graph, and let $n$ be a
parameter such that $W_G(v)\le n$ for every
vertex. Let $\mathcal C$ contain at most
$C_1n$ colors. Suppose that, for every edge
$e$, color $c$, and vertex $v$, the following
conditions hold:
\[
0\le x_{e,c}\le \frac{C_0}{n},
\qquad
\sum_{c\in\mathcal C}x_{e,c}\le 1,
\qquad
(x_{e,c}:e\ni v)\in\mathcal K_v
\quad(c\in\mathcal C).
\]
Then, for all sufficiently large $n$, there
is a feasible partial coloring using only
pairs $(e,c)$ with $x_{e,c}>0$, whose colored
incident weight at every vertex $v$ is at least
$\sum_{e\ni v}w_e\sum_{c\in\mathcal C}x_{e,c}
-\varepsilon n
-O_{C_0,C_1,\varepsilon}(n^{3/4})$.
Both the threshold and the hidden constant
are independent of the number of vertices
and the minimum positive edge weight.
\end{restatable}

Informally, the theorem converts locally
feasible \emph{fractional proposals} into
globally consistent coloring. Each edge may be proposed for several
colors, but its total proposal is at most one,
and its probability for any individual color
is $O(1/n)$. The theorem preserves almost
all of the proposed incident weight at every
vertex, losing at most
$\varepsilon n+O(n^{3/4})$.
Crucially, it also respects the support of
the proposals: an edge can receive only a
color for which $x_{e,c}>0$. This flexibility
will later allow us to apply the theorem to  the multigraph setting,
where some colors are already occupied and
hence unavailable to certain edges. The full proof is given in
\Cref{sec:round-proof}.

\paragraph{Why the theorem implies the weighted Vizing bound.}
Let $n=\bp(G)$, pad each local packing to $n$
bins, and use $n$ colors. As above, set
$x_{e,c}=1/n$ for every edge and color.
Each local vector is an average of feasible
bins and therefore belongs to $\mathcal K_v$,
while $\sum_c x_{e,c}=1$.

Applying \Cref{thm:overview-rounding} with
$C_0=C_1=1$ leaves at most
$\varepsilon n+O_\varepsilon(n^{3/4})$
incident weight uncolored at any vertex.
The greedy cleanup lemma
(\Cref{lem:general-cleanup}) colors these
remaining edges using at most four times
this quantity plus one fresh color. Hence,
$\cw(G)\le n+4\varepsilon n
+O_\varepsilon(n^{3/4})$.
Letting $n$ grow and then $\varepsilon$
tend to zero proves
$\cw(G)\le(1+o(1))\bp(G)$.
This proves the simple graph bound in \cref{thm:main}. 
See \Cref{sec:simple} for the details. 

\paragraph{Proof overview of \Cref{thm:overview-rounding}.}
The proof has four stages.
First, we replace each local feasible-bin distribution
by a bounded menu of \emph{packing profiles}, while
preserving almost all inclusion probabilities
(\Cref{ov:page-5}).
Second, we {\em couple} the fractional placement options
at the two endpoints of each edge, preserving their
common marginal probability without the $1/n$
loss caused by independent choices
(\Cref{ov:page-6}).
Third, we choose one profile for every vertex--color
pair, ensuring spare capacity in every selected
profile and controlling the total weighted assignment
error at every vertex
(\Cref{ov:page-6-selection}).
Finally, we turn the fractional assignments into
whole-edge decisions: hypergraph matching selects
nonconflicting placements for large edges, while a simple 
linear-algebraic rounding handles small edges
(\Cref{ov:page-7}).
Together, these stages lose at most
$\varepsilon n+O_{C_0,C_1,\varepsilon}(n^{3/4})$
proposed incident weight at every vertex.
The complete proof is given in
\Cref{sec:round-proof}, which proves
\Cref{thm:round}. We now give a high-level overview below. 

\subsection{Bounded types of bins by grouping}\label{ov:page-5}
To coordinate local packings, we first replace each bin by a short description of the \emph{space it offers}, rather than of the particular edges it contains.
Fix a small constant $\beta>0$. Call an edge \emph{large} if its weight is at least $\beta$, and \emph{small} otherwise. A bin contains at most $1/\beta$ large edges. We describe its large-edge capacity by slots of sizes $b_1,\ldots,b_k$ and put aside space $R$ for all small edges, with
$b_1+\cdots+b_k+R\le 1$.
A slot can receive any one large edge whose weight is at most the slot size. The space $R$ can receive several small edges whose combined weight is at most $R$. We call this description a \emph{profile}.

Even with a constant number of large slots per bin, their sizes could take arbitrarily many values. If we allow all of them, later local choices could be too complicated. We therefore replace the possible profiles by a bounded number. This idea is similar in spirit to other grouping techniques widely used in bin packing, such as linear grouping \cite{fernandez1981bin} and harmonic grouping \cite{LeeLee85}.

Simply rounding slot sizes \emph{up} would be dangerous, since the resulting profile might no longer fit in a unit bin. Instead, we round sizes \emph{down}, keeping only a bounded collection of representative sizes.
This can lose some chances to place a large edge. The remedy is to reserve a small extra probability for a profile consisting of one full-size slot; it compensates for the lost slot capacity at every size threshold. 
We similarly round down the space reserved for small edges and assign a small probability to a profile with a unit-capacity reserve for small edges.
By first sacrificing a small fraction of the proposed probabilities, the repairs can be paid for within total probability one.

The compression need not reproduce each original bin. It suffices to preserve fractional placement opportunities for every edge: large slots are limited by item counts, whereas small-edge reserves are limited by total weight. Competing fractional requests will be resolved after the two endpoints have been coordinated.
A profile captures \emph{capacity}, not the identities of edges assigned to it. Thus an edge may use any compatible slot, even if that slot originally represented another edge. The next lemma formalizes the crucial compression step: a bounded collection of valid profiles can preserve a $(1-\delta)$ fraction of \emph{each edge's} inclusion probability, while respecting both slot counts and reserved small-edge weight. The bounded number of profiles is essential later: once every profile has probability bounded below, conditioning on a selected profile does not amplify a single edge's $O(1/n)$ proposal uncontrollably.

\begin{lemma}[Bounded-profile compression; informal version of \Cref{lem:profiles}]
\label{lem:overview-profile}
Fix $\delta\in(0,1)$ and $0<\beta\le\delta/4$.
For every $x\in\mathcal K_v$, there exist
a probability distribution $(\lambda_P)_P$
over at most $M=M(\beta,\delta)$ valid profiles,
and nonnegative assignments $a(e;P,i)$ to
compatible slots (for large edges) or reserved
spaces (for small edges), satisfying the following
properties:
\begin{enumerate}
    \item \emph{Almost all marginals are preserved:}
    for every edge $e\in E(v)$,
    $\sum_{P,i}a(e;P,i)=(1-\delta)x_e$.

    \item \emph{Every profile respects its capacities:}
    for each profile $P$, every large-edge slot $i$
    receives total assigned fraction at most
    $\lambda_P$, and its reserved small-edge space
    receives total assigned weight at most
    $\lambda_P R(P)$.
\end{enumerate}
The number of profiles $M$ depends only on
$\beta$ and $\delta$, independently of the degree
of $v$ and the minimum positive edge weight.
\end{lemma}

\noindent\emph{Proof sketch.}
Round slot sizes and small-edge spaces downward to a bounded grid.
The lost large-slot capacity at any size threshold is compensated by
adding a small probability of a profile with one unit slot; a
unit-reserve profile compensates the lost small-edge capacity.
Scale the original inclusion probabilities by $1-\delta$ to pay for
both repairs. Assign larger edges first to available compatible slots,
and assign small edges proportionally to reserved capacities.
The threshold-count inequalities prove the stated capacities; see \Cref{sec:local} and  \Cref{lem:profiles} for the precise construction.

\subsection{Coupling: coordinating the two endpoints without losing probability}
\label{ov:page-6}

We now coordinate the possible placements of an edge at its
two endpoints. Recall that each profile contains slots for
large edges and reserved space for small edges. We use the
term \emph{place} to refer to either type of capacity in
a particular profile.
Fix an edge $e=uv$ and a color $c$.
After the compression step, suppose that $e$ has retained
inclusion probability $y$ at both endpoints.
At $u$, this probability is distributed among its possible
places, with fractions $a_1,a_2,\ldots$ satisfying
$\sum_i a_i=y$. Similarly, at $v$, the fractions
$b_1,b_2,\ldots$ satisfy $\sum_j b_j=y$.

The challenge is to combine these two local assignments
into \emph{joint placement options}: each option must
specify a place at $u$ and a place at $v$ for the
\emph{same edge and color}. Simply multiplying the
fractions would give total mass $y^2$, losing a factor
of $y$.

Instead, for $y>0$, assign mass
${\frac{a_i b_j}{y}}$
to the option pairing place $i$ at $u$ with place $j$
at $v$. The normalization by $y$ preserves both
endpoints' original assignments. Indeed, for every
place $i$ at $u$,
$\sum_j\frac{a_i b_j}{y}=a_i$,
and symmetrically, the total mass assigned to each
place $j$ at $v$ is $b_j$.
Consequently, the total mass of all paired options
is exactly $y$, rather than $y^2$.
If $y=0$, no pairing is needed.

Thus, we obtain \emph{joint fractional placement options}
without losing any additional inclusion probability.
Crucially, these fractions represent possible assignments,
not pieces of physical edges: no edge has yet been colored.
The next step is to select actual profiles and turn these
fractional options into feasible assignments of whole edges.

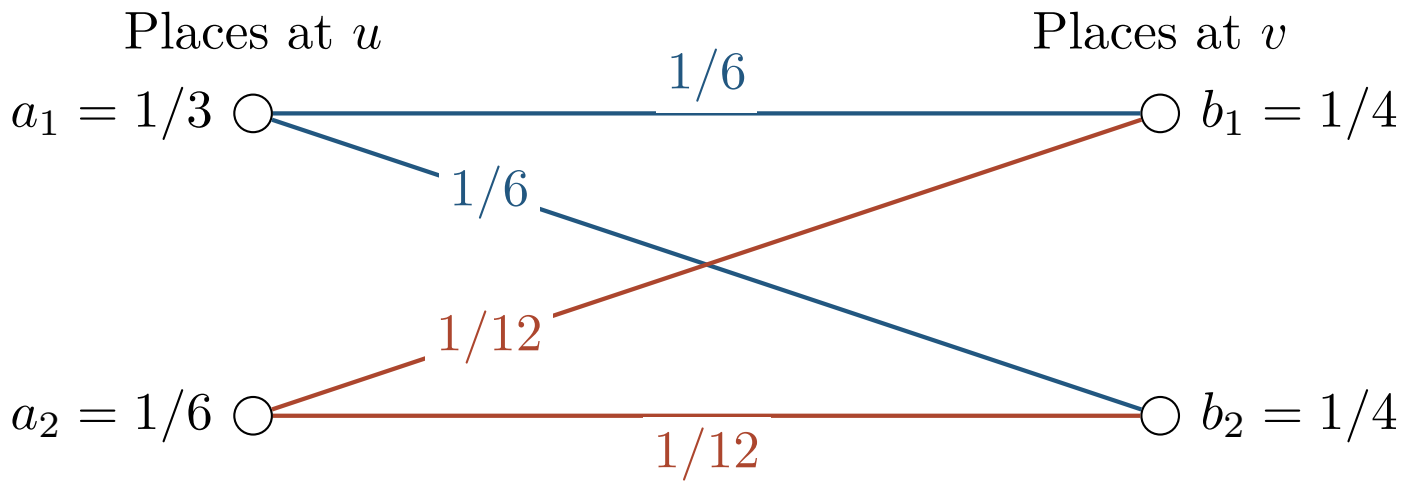
\begin{figure}[ht]
\centering
\begin{tikzpicture}[font=\small,x=1cm,y=1cm]
\node at (0,1.55) {Places at $u$};
\node at (6,1.55) {Places at $v$};
\node[vertex,label=left:{$a_1=1/3$}] (u1) at (0,1) {};
\node[vertex,label=left:{$a_2=1/6$}] (u2) at (0,-1) {};
\node[vertex,label=right:{$b_1=1/4$}] (v1) at (6,1) {};
\node[vertex,label=right:{$b_2=1/4$}] (v2) at (6,-1) {};
\draw[inkblue,thick] (u1)--(v1) (u1)--(v2);
\draw[warmred,thick] (u2)--(v1) (u2)--(v2);
\path (u1)--node[pos=.5,above,fill=white,inner sep=2pt,text=inkblue] {$1/6$}(v1);
\path (u1)--node[pos=.25,fill=white,inner sep=2pt,text=inkblue] {$1/6$}(v2);
\path (u2)--node[pos=.25,fill=white,inner sep=2pt,text=warmred] {$1/12$}(v1);
\path (u2)--node[pos=.5,below,fill=white,inner sep=2pt,text=warmred] {$1/12$}(v2);
\end{tikzpicture}
\caption{Pairing the fractions of one edge and one color, with $y=1/2$. Each link has mass $a_i b_j/y$. The link masses at every place sum to its original fraction; the total mass remains $1/2$.}
\label{fig:overview-linking}
\end{figure}

\subsection{Choosing the local profiles consistently}\label{ov:page-6-selection}
The coupling above makes the two endpoint assignments agree \emph{on average}. It does not yet ensure that a particular selection of local profiles has enough capacity. We must choose one profile at each vertex for each color; these choices affect the possible assignments of all incident edges.

Imagine first choosing these profiles independently. If a paired option uses profiles chosen with probabilities $\lambda_u$ and $\lambda_v$, give it conditional request equal to its paired mass divided by $\lambda_u\lambda_v$. Averaging over the profile choices cancels the denominators and recovers the intended mass. This is why the coupling remains useful after choosing profiles.
But there are two ways an actual choice can fail. \emph{Overcrowding}: many edges request the same selected slot, exceeding its capacity. \emph{Lost opportunities}: an edge has too little total request across all colors, leaving too much incident weight for the final cleanup. We must prevent both failures simultaneously, at every vertex.

We first give every profile  (of the bounded number of types) a small positive probability. Thus no chosen profile is extremely rare, and an individual edge's request in one color remains $O(1/n)$. We also reduce the requests slightly to leave spare capacity. This space is needed later, when fractional options are converted into integral assignments. 

\paragraph{The two invariants we need.}
After choosing one profile at each vertex and color, let $S_e$ be the total
fractional mass of the resulting options for edge $e$, and let $\mu_e$ be
its intended mass before the profiles were fixed (after the deliberate
small losses). We seek two guarantees, simultaneously:
\begin{enumerate}
\item \emph{Capacity slack:} every selected large-edge slot and every
small-edge reserve receives at most $(1-\sigma')$ times its capacity,
for some constant $\sigma'>0$ depending only on the accuracy.
\item \emph{Small weighted discrepancy:} for every vertex $v$,
\begin{equation}\label{eq:overview-discrepancy}
 \sum_{e\ni v}w_e\lvert S_e-\mu_e\rvert=O_\varepsilon(n^{3/4}).
\end{equation}
\end{enumerate}
The first guarantee makes subsequent integral rounding feasible. The
second says that the total \emph{weight} of lost assignment opportunities
is small even when a vertex has arbitrarily many tiny incident edges.
It is deliberately weaker than requiring $S_e\approx\mu_e$ for every edge.

\paragraph{Fixing one group at a time.}
A profile list at a vertex specifies its chosen profile for every
color. We choose these complete lists in groups, keeping all earlier
choices fixed. Partition the vertices into a fixed number
$r=O_\varepsilon(1)$ of groups so that edges within one group have
total weight at most $n/r+O(\sqrt{n\log n})$ at each vertex.
Choose $r$ large enough and set these edges aside for cleanup.
There are now no edges within a group. Consequently, for a vertex in
the current group, its incident assignment mass depends on its own list and on lists in
other groups, but not on the other choices made at the same stage.
This is what allows us to consider candidate lists separately at each vertex.

\paragraph{Measuring the error still to come.}
We need to control the assignment error while fixing the
profile lists one group at a time. After the first $k$
groups have been chosen, let $\mathcal F_k$ denote their
fixed lists. For the remaining groups, imagine that
their lists are still chosen independently from their
original distributions. We use this randomness only
to measure the remaining error; the algorithm will
choose these lists more carefully.
Define
$Q_v^{(k)}
:=
\sum_{e\ni v} w_e
\sqrt{\mathbb E\bigl[(S_e-\mu_e)^2
  \mid \mathcal F_k\bigr]}$.
For each edge, the square-root term bounds its expected error by Cauchy--Schwarz. It measures both
the error caused by the lists already fixed and the
uncertainty from the lists not yet chosen. Once all
lists are fixed, the expectation disappears, and
$Q_v^{(r)}=\sum_{e\ni v}w_e|S_e-\mu_e|$ -- the discrepancy we want to bound.

This measure is useful for two reasons. 
First, we can compute it from the first and second
moments of the contributions from different colors.
Second, if we choose the next group according to its
original distribution, the expected value of
$Q_v^{(k)}$ cannot increase. This follows from the
tower property of conditional expectation and the
concavity of the square root. Thus, we can fix
the groups successively while controlling this
measure.

Initially, each $S_e$ is a sum of $O(n)$ independent
contributions, each of size $O(1/n)$. Hence
$\operatorname{Var}(S_e)=O(1/n)$.
Since $W_G(v)\le n$, we obtain
$Q_v^{(0)}=O(\sqrt n)$.
The following argument shows how to keep this error
sublinear while making all profile choices feasible.

\paragraph{Rejecting unsuitable lists and controlling their neighbors.}
When choosing the profiles of the current group, we must
satisfy two requirements. Each vertex must retain enough
assignment opportunities for its own incident edges.
At the same time, the choices made by the group must not
overload capacities or create too much error at neighboring
vertices.

We handle the first requirement by rejecting unsuitable
profile lists. For a vertex $u$ in the current group,
we accept a candidate list only if its resulting
potential $Q_u$ is at most $n^{3/4}$.
Because no two vertices in the current group are adjacent,
this test depends only on the candidate at $u$, the
previously fixed lists, and the distributions of the
lists still to be chosen. In particular, it does not
depend on the choices of other vertices in the current
group.
Our induction maintains
$Q_u=O_\varepsilon(\sqrt{n\log n})$
before $u$ is chosen. By Markov's inequality, a random
candidate is therefore rejected with probability at most
$O_\varepsilon(n^{-1/4}\sqrt{\log n})=o(1)$.
Thus, we can repeatedly sample until the test succeeds,
while discarding only a negligible fraction. Since the acceptance tests are
separate for different vertices in the current group,
their accepted lists remain independent.

We next choose these accepted lists so that they also
satisfy the requirements at neighboring vertices.
For this, we use Harris--Srinivasan partial resampling
\cite{HS} (stated as \Cref{prop:HS}).
The procedure resamples selected vertex lists whenever
a constraint is violated. It simultaneously controls
the load assigned to every slot or reserve and the
increase in the error potential at vertices outside
the current group. Importantly, these guarantees do
not require a union bound over all vertices.
We impose the load constraints for every possible
profile at a receiving vertex. They therefore remain
valid even if that vertex's profile has not yet
been chosen.

Here the assumption that $G$ is simple is crucial.
Changing the list at $u$ affects any fixed neighbor
$v$ through at most one edge, while the total weight
of edges incident to $u$ is at most $n$.
Together with the bounded number of profile types,
these properties limit how much one choice can affect
the constraints. This is what makes partial resampling
applicable with error bounds depending only on $n$
and the accuracy, rather than the size of the graph.

\paragraph{Keeping both guarantees through the last group.}
The key is that neither rejection sampling nor partial
resampling introduces much additional error.
Rejecting a probability $p=o(1)$ increases the
expectation of any nonnegative quantity by at most
a factor $1/(1-p)=1+o(1)$.
Partial resampling adds only
$O_\varepsilon(\sqrt{n\log n})$ to the error
bound at each stage.
Before a vertex's list is chosen, we need the smaller error bound
so that rejection is rare. After its list is chosen, it suffices
to preserve the larger $O_\varepsilon(n^{3/4})$ bound needed for cleanup.
Since the number of groups is fixed, these bounds
maintain $Q_v=O_\varepsilon(\sqrt{n\log n})$
until the list at $v$ is chosen. Its acceptance test
then ensures $Q_v\le n^{3/4}$.
Although later choices may change the assignments
of edges, the resampling constraints
continue to control its error. Consequently,
its final discrepancy is $O_\varepsilon(n^{3/4})$.

The capacity constraints are preserved for a similar
reason. The load error introduced at each stage is
$o(1)$ after scaling back to the original units.
Every positive slot or reserve has capacity bounded
below by a fixed constant. Therefore, the total
error over the constant number of stages can be
absorbed by the spare capacity introduced earlier.

Thus, after all profiles have been chosen, every
selected slot and reserve retains positive capacity
slack, and every vertex $v$ satisfies
$\sum_{e\ni v} w_e|S_e-\mu_e|
= Q_v^{(r)}
= O_\varepsilon(n^{3/4})$.
The full induction is given in
\Cref{sec:round-proof}, while the computation of
the required moments is explained in
\Cref{sec:rounding-implementation}.

\subsection{Turning the coordinated assignments into a coloring}\label{ov:page-7}
After the profile choices, each edge has fractional options that name \emph{the same color at both endpoints} and compatible places in the chosen profiles. If an edge has total proposed mass greater than one, we scale its options down to total mass one; the preceding error estimate pays for this reduction. We must now choose whole options without assigning two competing edges to the same place.

The two types of places impose different constraints. A large-edge slot
can contain \emph{at most one edge}, so selecting large-edge options
requires a matching. A small-edge reserve instead constrains only the
\emph{sum of weights}, and therefore permits a simpler load-preserving
rounding. Both methods choose an entire outcome for each edge.

\paragraph{Large edges: choose nonoverlapping options.}
A large-edge option consists of three objects: the edge
$e=uv$ itself, one slot for color $c$ at $u$, and one
slot for the same color at $v$. Selecting pairwise
disjoint triples ensures that each edge receives at
most one color and each slot is used at most once.
Thus, the selected options give a feasible coloring.
Some edges have total fractional assignment mass
less than one. We therefore allow another outcome:
leaving the edge uncolored. We represent this by
\emph{extra slots} at its endpoints. These are only
bookkeeping devices, not additional colors or usable
packing space. We create extra slots in each weight
range in proportion to the missing assignment mass,
with a small rounding allowance. This controls the
total weight left uncolored at each vertex.

We form a $3$-uniform hypergraph with one \emph{job}
vertex for each large edge and one \emph{resource}
vertex for each ordinary or extra slot. A hyperedge
joins a job to the two slots of one possible outcome.
A matching covering every job therefore chooses
exactly one outcome per large edge, without using
any slot twice.
We might ask why such a matching exists.
Note that every job has many possible outcomes, while individual
resources and pairs of vertices occur in relatively
few options. The simplicity of $G$ is important:
two ordinary slots at opposite endpoints can be
linked by at most one original edge. Moreover, each
edge--color proposal has fractional mass $O(1/n)$,
and the extra slots are distributed to keep their
individual overlaps small.

To apply the Delcourt--Postle matching theorem (\Cref{prop:DP} in the appendix), we turn fractional option masses into
integer hyperedge multiplicities using the scale
$N=\lceil n\rceil^6$. Let $d(x)$ be the degree of
vertex $x$, and $\Delta_2$ be the maximum number
of hyperedges containing any pair of distinct
vertices (the maximum pair-codegree), counting
multiplicities. For fixed $\eta>0$, the appendix
establishes
$d(z) = N-O(n^2)$ for each job $z$; 
$d(b) \le (1-\eta)N$ for each resource $b$; and 
$\Delta_2 = O(Nn^{-2/3})$.
The first bound follows from the unit total mass
of each job's options, including the uncolored
outcome, with only $O(n^2)$ copies lost to rounding.
Capacity slack gives the second bound. The third follows
from the small option masses, the simplicity of $G$,
and the distribution of extra slots.
Fix $0<\zeta<1/9$, let $a>0$ be the exponent supplied by
\Cref{prop:DP}, and set $D=\lfloor(1-\eta/2)N\rfloor$.
For sufficiently large $n$, the job degrees are at least
$(1+D^{-a})D$, since their fixed relative surplus above $D$
exceeds $D^{-a}$. The resource degrees are at most $D$, and
$\Delta_2=O(Nn^{-2/3})=o(D^{1-\zeta})$, so
$\Delta_2\le D^{1-\zeta}$ for sufficiently large $n$.
The Delcourt--Postle theorem~\cite{DP} therefore gives a matching
covering every job. The extra slots also bound
the total uncolored weight at each vertex by
the original missing fractional mass, up to
a small error, using a constant number of
narrow weight intervals. 
Details and
the algorithmic implementation appear in
\Cref{sec:round-proof,sec:matching-algorithm}.

\paragraph{Small edges: round loads, not individual slots.}
A small-edge option uses one reserved space at each
endpoint, contributing weight less than $\beta$ to each.
Unlike large edges, small edges do not require a matching:
several edges may share a reserve, provided their total
weight fits. We therefore need only control the total
load in each reserve and ensure that every edge is
either assigned a color or left uncolored.
Every positive reserve has size at least
$r_0:=\delta/4$ (\Cref{lem:profiles}).
After normalizing its capacity to one, each option
contributes at most $\tau=\beta/r_0$ to each of at most
two constraints. The two-constraint rounding lemma
(\Cref{lem:two}) makes one choice for each small edge,
increasing every normalized load by at most $2\tau$.
By choosing $\beta$ sufficiently small relative to
the available slack, all reserves remain feasible.
We also track the total weight left uncolored at
each vertex, which increases by only a small
additive amount during rounding.

\paragraph{Accounting for the losses.}
We now combine the losses from the different steps.
First, dividing the vertices into groups, simplifying
the local bins, and leaving some spare capacity
lose only an arbitrarily small fraction of the
proposed incident weight.
Second, after choosing the profiles, an edge may
have total fractional assignment greater than one.
We reduce its assignment to one in this case.
The discrepancy bound \eqref{eq:overview-discrepancy} shows that this
reduction loses at most $O_\varepsilon(n^{3/4})$
total weight at each vertex.
Finally, converting fractional assignments into
whole edges may leave some edges uncolored.
For large edges, the extra slots ensure that the
weight left uncolored is bounded by the weight
not covered by the fractional assignments,
up to a factor arbitrarily close to one and an
$O_\varepsilon(n^{2/3})$ additive error.
For small edges, rounding adds only a bounded
error at each vertex.
By choosing the fixed parameters sufficiently
small, the total loss at every vertex is at most
$\varepsilon n+O_\varepsilon(n^{3/4})$.
Thus,
the total weight of colored edges at each vertex $v$ is at least 
$\sum_{e\ni v}w_e\sum_{c\in\mathcal C}x_{e,c}
-\varepsilon n-O_\varepsilon(n^{3/4})$,
as claimed in \Cref{thm:round}.
The complete proof and algorithm appear in
\Cref{sec:round-proof} and \Cref{subsec:simpalg}.
 
\label{ov:simple-end}

\section{Tight Weighted Guarantees for Bipartite Multigraphs}
\label{sec:overview-bipartite}\label{ov:page-8}

We prove the asymptotically tight bound
$\chi'_w(G)\le(2+o(1))\Delta_w(G)$ for bipartite multigraphs.
The idea is to group parallel edges into feasible bundles, color the
heavier bundles, and use our local-to-global rounding theorem to color
almost all remaining edges. The difficulty is that a remaining edge
can use only a color missing at \emph{both} endpoints.

\paragraph{Grouping parallel edges.}
Put $n=\lceil\Delta_w(G)\rceil$. On each pair of vertices, start with
one group per edge and repeatedly merge two groups of weight at most
$1/2$. Each merged group has weight at most one. When no further merge
is possible, at most one group of weight at most $1/2$ remains on each
pair. Let $H$ contain the heavier groups and $R$ the remaining groups,
regarding each group as one weighted edge. Thus $H$ is a bipartite
multigraph, while $R$ is simple and every edge of $R$ weighs at most
$1/2$. Here, giving a group one color means assigning that color to all its
original edges.

\paragraph{Coloring the heavy groups.}
Since every edge of $H$ weighs more than $1/2$, we have
$d_H(v)\le 2n$. By K\H{o}nig's theorem, $H$ admits a proper coloring
with at most $2n$ colors. Fix $0<\varepsilon<1$ and pad the palette to
the least even integer $K\ge(2+\varepsilon)n$; for sufficiently large
$n$, we have $K\le4n$. Let $F_v$ be the set of colors missing at $v$,
and write $f_v=|F_v|=K-d_H(v)$. Since the heavy groups contribute
more than $d_H(v)/2$ weight at $v$, we have
 $\frac12 d_H(v)+W_R(v)\le n$,
 and
 $f_v\ge\varepsilon n+2W_R(v)$.
Each vertex therefore has ample missing colors relative to its
remaining incident weight. But the sets $F_u$ and $F_v$ may overlap
very little, so this local count alone does not provide colors for
an edge $uv\in R$.

\paragraph{Making missing colors agree.}
A \emph{Kempe switch} exchanges two colors along an
entire two-color path of $H$. It preserves proper
coloring and changes missing colors only at the path
endpoints; in particular, every $f_v$ remains
unchanged. See Misra and Gries~\cite{MisraGries92} for this
operation in a constructive proof of Vizing's theorem.

We use such switches not to guarantee large overlap
for each edge separately, but to construct small
fractional color proposals whose total uncovered
weight is small at every vertex. In each of
$O(\log n)$ stages, we pair the colors and decide
which two-color paths to switch. These switches
help distribute the missing colors more evenly,
so that most of the remaining incident weight
comes from edges with sufficiently many common
missing colors. Harris--Srinivasan resampling
then allows us to choose the switches while
ensuring these guarantees simultaneously at
every vertex.

\begin{lemma}[Missing-color coordination; special case of \Cref{lem:weighted-moment-mix}]
\label{lem:overview-missing}
In the setting above, for all sufficiently large $n$, Kempe switches
in $H$ produce missing-color sets $F_v$ and fractions $x_{e,c}$,
for $e\in E(R)$ and $c\in[K]$, such that:
\begin{enumerate}
 \item \emph{Common missing colors:} $x_{e,c}>0$ implies
 $c\in F_u\cap F_v$ for $e=uv\in R$.
 \item \emph{Small proposals:} $0\le x_{e,c}=O_\varepsilon(1/n)$
 and $\sum_c x_{e,c}\le1$ for every $e\in R$.
 \item \emph{Small load and loss:} for every vertex $v$ and color $c$,
 $\sum_{e\in E_R(v)}w_ex_{e,c}
  \le\frac{W_R(v)}{f_v}+o(1)$,
 and for every vertex $v$,
 $\sum_{e\in E_R(v)}w_e\Bigl(1-\sum_cx_{e,c}\Bigr)=o(n)$.
\end{enumerate}
The error bounds are independent of graph order and the minimum positive edge weight.
\end{lemma}

The proof of \Cref{lem:overview-missing} follows by applying
\Cref{lem:weighted-moment-mix} to the residual graph with charges
$q_e=2w_e$ and size parameter $2n$.
See \Cref{sec:missing-colors} for the complete argument. 

Now we complete the proof of
$\chi'_w(G)\le(2+o(1))\Delta_w(G)$. The complete proof is given in Appendix~\ref{sec:bip-multi}.
Section~\ref{sec:parallel-grouping} describes how parallel
edges are grouped into heavy bundles and a simple residual
graph. Section~\ref{sec:missing-colors} proves the
missing-color coordination lemma. 
Section~\ref{sec:mixing-implementation} explains how
the switches can be found in expected polynomial time.
The subsequent subsection combines the missing-color
guarantees with Theorem~\ref{thm:round} to color the residual edges.
Finally, Section~\ref{sec:bipartite-total-proof}
establishes the final bound
and the algorithm. Below we give a short proof sketch. 

\begin{proof}[Proof sketch]
Apply \Cref{lem:overview-missing} to $R$. Since
$f_v\ge\varepsilon n+2W_R(v)$ and $f_v\le K\le4n$,
we have
 $\frac{W_R(v)}{f_v}
 \le\frac12-\frac{\varepsilon n}{2f_v}
 \le\frac12-\frac\varepsilon8$.
Thus, for all sufficiently large $n$,
$\sum_{e\in E_R(v)}w_ex_{e,c}\le 1/2-\varepsilon/16$.
Every edge of $R$ weighs at most $1/2$, so  we can  ensure (via \Cref{lem:half-load}) that $(x_{e,c}:e\ni v)\in\mathcal K_v$: each local vector is
an average of feasible integral bins.

We can therefore apply \Cref{thm:overview-rounding} to the simple
graph $R$, with size parameter $n$, $C_1=4$, and
$C_0=O_\varepsilon(1)$. With any fixed rounding accuracy $\delta>0$,
the uncolored incident weight at $v$ is at most
$\sum_{e\in E_R(v)}w_e\Bigl(1-\sum_cx_{e,c}\Bigr)
      +\delta n+O_{\varepsilon,\delta}(n^{3/4})
 =\delta n+o(n)$.
By the support condition in \Cref{lem:overview-missing}, every color
used on $R$ is missing at both endpoints in $H$. Hence the partial
colorings are compatible. The greedy cleanup lemma (\Cref{lem:general-cleanup}) colors the
remaining edges with at most $4\delta n+o(n)$ fresh colors, and
expanding groups preserves feasibility. Therefore
$\chi'_w(G)\le(2+\varepsilon+4\delta)n+o(n)$.
Letting $\varepsilon,\delta>0$ be arbitrarily small and using
$n=\Delta_w(G)+O(1)$ proves the result.
The full completion argument is \Cref{lem:finish-packets};
see also \Cref{thm:bipartite-total-weight}.
\end{proof}
\label{ov:bipartite-end}

\section{Lower Bounds and the Feige--Singh Conjecture}
\label{sec:lower-bounds}
In this section, we~give two constructions showing how agreement between endpoints can require more colors than local packing alone.
The first is~a~simple bipartite graph that disproves the exact total-weight conjecture.
The second has $\chi'_w(G)/\bp(G)=4/3$ for arbitrarily large $\bp(G)$, even when the underlying simple graph is~a~tree of~maximum degree three.

Here, we want to stress that the total weight is a different benchmark.
On a simple star, $t$ edges of weight $1/2+1/(2t)$ require $t$ colors, although their total weight is $(t+1)/2$. Thus a coefficient below two in terms of $\Delta_w$ is impossible even for simple bipartite graphs. This does not contradict the simple-graph theorem: the same star has $\bp=t$.

\subsection{A~Counterexample to~the Feige--Singh Exact Total-Weight Conjecture}\label{sec:counterexample}
The total-weight conjecture asks for at~most $2N-1$ colors, where $N=\lfloor\Delta_w(G)\rfloor+1$~\cite[Conjecture~2]{FS}.
We~construct a~simple bipartite graph with $\Delta_w(G)<2$ that needs four colors.

\begin{proposition}\label{prop:counterexample}
	There is~a~simple bipartite weighted graph with
	\[
	\Delta_w(G)=1.95<2,\qquad \bp(G)=3,\qquad \chi'_w(G)=4.
	\]
\end{proposition}
\begin{proof}
	Take the two vertex sets
$\{P,a_1,a_2,a_3,v\}$, and $\{b_1,b_2,b_3,u\}$.
	The edges and weights are
	\[
	\begin{array}{c|c|c}
		\text{edge}&\text{indices}&\text{weight}\\ \hline
		p_i=Pb_i&i\in\{1,2,3\}&0.65\\
		q_{ij}=a_jb_i&i,j\in\{1,2,3\},\ i\ne j&0.6\\
		x_j=a_ju&j\in\{1,2,3\}&0.45\\
		e=vu&&0.6
	\end{array}
	\]
	See \Cref{fig:counterexample}.
	The weighted degrees at~$P,u,b_i,a_j,v$, respectively, are $1.95, 1.95, 1.85, 1.65, 0.6$.
	Thus $\Delta_w(G)=1.95$.
	Every vertex has a~packing into three bins: all but $u$ have degree at~most three, and at~$u$ two of~the three $0.45$ edges share a~bin while the other two edges get one bin each.
	The three edges at~$P$ have pairwise sums greater than one, so~they require three bins.
	Hence $\bp(G)=3$.

\begin{figure}[htbp]
		\centering
\begin{tikzpicture}[x=1cm,y=1cm,font=\small]
\node[vertex,fill=softblue,label=left:$P$] (P) at (0,0) {};
\foreach \ii/\yy in {1/1.4,2/0,3/-1.4}{
 \node[vertex,fill=softgold,label=above:$b_{\ii}$] (b\ii) at (2.6,\yy) {};
 \node[vertex,fill=softblue,label=above:$a_{\ii}$] (a\ii) at (6.2,\yy) {};
}
\node[vertex,fill=softgold,label=right:$u$] (u) at (8.8,0) {};
\node[vertex,fill=softblue,label=below:$v$] (v) at (8.8,-2.4) {};
\foreach \ii in {1,2,3}{
 \draw[inkblue,thick] (P)--(b\ii);
 \draw[warmred,thick,dashed] (a\ii)--(u);
 \foreach \jj in {1,2,3}{
  \ifnum\ii=\jj\else
   \draw[black!65,thick] (b\ii)--(a\jj);
  \fi
 }
}
\draw[black!65,thick] (u)--node[right] {$0.6$}(v);
\node[inkblue] at (1.3,2.1) {$p_i:\ 0.65$};
\node at (4.4,2.1) {$q_{ij}:\ 0.6$};
\node[warmred] at (7.5,2.1) {$x_j:\ 0.45$};
\node[align=center] at (4.4,-2.2) {Every middle pair is joined\\except $a_i b_i$.};
\end{tikzpicture}
		\caption{A simple bipartite graph with maximum weighted degree below
		two that needs four colors. Blue-filled vertices form one side and
		yellow-filled vertices the other. Edge styles identify weights, not
		an edge coloring. Crossings of drawn edges are not vertices.}
		\label{fig:counterexample}
	\end{figure}
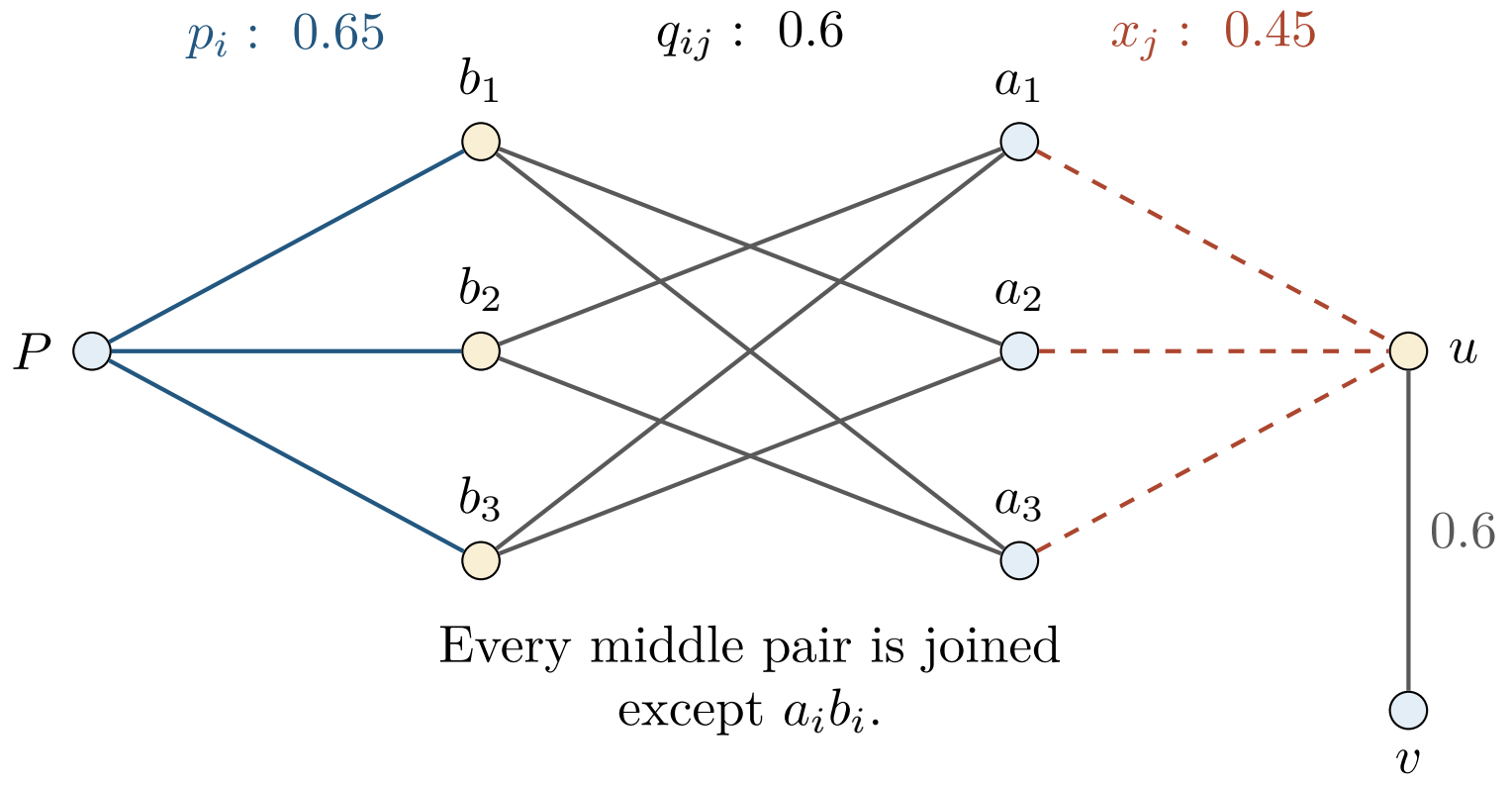

	Suppose three colors suffice.
	At~each of~$P,a_1,a_2,a_3$, every pair of~incident edges has total weight greater than one.
	Thus all three colors occur once at~each of~these four vertices.
	Counting the twelve edges other than $e$ at~their endpoints in~this set, each color occurs exactly four times.

	At~each of~$b_1,b_2,b_3$, the three incident weights also have pairwise sums greater than one, so~all three colors occur once there.
	These vertices account for three occurrences of~each color.
	The remaining edges $x_1,x_2,x_3$ must therefore account for one occurrence of~each color.
	At~$u$, the edge $e$ cannot share a~color with any $x_j$, because $0.6+0.45>1$.
	All three colors are blocked, a~contradiction.

	Finally, the ordinary maximum degree is~four.
	K\H{o}nig's theorem gives a~proper edge coloring with four colors, which is~feasible since all weights are at~most one.
	Thus $\chi'_w(G)=4$.
\end{proof}

For this graph, $N=2$, so~the conjectured bound is~$2N-1=3$.
The Chung--Ross conjecture uses $\bp(G)=3$ and permits 5 colors.
Thus, the construction disproves the exact total-weight conjecture while leaving the Chung--Ross conjecture and the asymptotic bound of~\Cref{thm:bipartite-total-weight} unaffected.

\subsection{A Four-Thirds Lower Bound}
\label{subsec:4by3lower}
For simple trees, local bin packings can be aligned along the tree,
giving $\cw(G)=\bp(G)$; see the proof of Theorem~10
in~\cite{Sannyasi20}. For multigraphs whose underlying simple graph is
a tree, Sannyasi proved the upper bound
$1.693\,\bp(G)+12$~\cite[Theorem~11]{Sannyasi20}.
The following construction shows that parallel edges force a factor
of $4/3$, even when the underlying tree has maximum degree three.
Our bipartite upper bound gives an asymptotic coefficient of $1.62$
for this class.
\begin{proposition}
	\label{prop:bipartite-four-thirds}
	For every integer $N\ge1$, there is~a~bipartite weighted multigraph $G_N$ whose underlying simple graph is~a~tree of~maximum degree three and whose weights are rational and greater than $1/3$, such that
	$\bp(G_N)=3N,\qquad \cw(G_N)=4N$.
\end{proposition}
\begin{proof}
	Set $m=2N$.
	Take a~path $v_0v_1\cdots v_m$ and add a~private leaf $u_i$ adjacent to~each $v_i$.
	For $0\le i<m$, give the path edge $v_iv_{i+1}$ weight
	$w_i=\frac13+\frac{m-i}{6(m+1)}$;
	then $1/2>w_0>\cdots>w_{m-1}>1/3$.
	Give $v_iu_i$ weight $1-w_i$ for $i<m$ and give $v_mu_m$ weight $1$.
	Replace every edge just described by~$2N$ parallel copies of~the same weight.

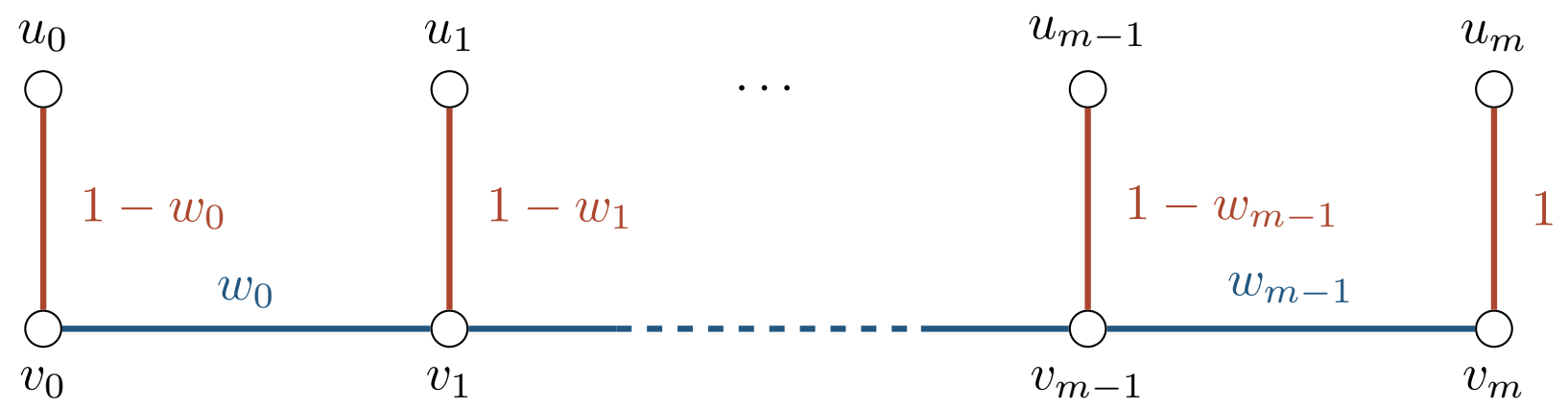
\begin{figure}[htbp]
	\centering
\begin{tikzpicture}[x=1cm,y=1cm,font=\small]
\foreach \name/\xx/\lab in {a/0/0,b/2.8/1,c/7.2/{m-1},d/10/m}{
 \node[vertex,label=below:$v_{\lab}$] (\name) at (\xx,0) {};
 \node[vertex,label=above:$u_{\lab}$] (u\name) at (\xx,1.65) {};
}
\draw[very thick,inkblue] (a)--node[above] {$w_0$}(b);
\draw[very thick,inkblue] (b)--(3.95,0);
\draw[very thick,inkblue,dashed] (3.95,0)--(6.05,0);
\draw[very thick,inkblue] (6.05,0)--(c);
\draw[very thick,inkblue] (c)--node[above] {$w_{m-1}$}(d);
\draw[very thick,warmred] (a)--node[right=3pt] {$1-w_0$}(ua);
\draw[very thick,warmred] (b)--node[right=3pt] {$1-w_1$}(ub);
\draw[very thick,warmred] (c)--node[right=3pt] {$1-w_{m-1}$}(uc);
\draw[very thick,warmred] (d)--node[right=3pt] {$1$}(ud);
\node at (5,1.65) {$\cdots$};
\node[align=center] at (5,-1) {Every segment represents $2N$ parallel edges.\\
Path weights decrease from left to right; $m=2N$.};
\end{tikzpicture}
	\caption{The exact $4/3$ family. A leaf edge fits with one path edge
	to its right, but cannot fit with a path edge to its left.}
	\label{fig:bp-four-thirds}
\end{figure}

	At~an~internal vertex $v_i$, pair each of~its $2N$ leaf edges with one of~the $2N$ path edges to~its right; these pairs weigh exactly one.
	Pair the $2N$ path edges to~its left two per bin, using $N$ more bins.
	The endpoint $v_0$ needs only $2N$ bins.
	At~$v_m$, the $2N$ unit leaf edges use their own bins and the path edges use another $N$ bins.
	This is~optimal there, because a~unit edge cannot share a~bin and three path edges weigh more than one.
	Thus $\bp(G_N)=3N$.

	We~show that one color contains at~most $m+1$ edges.
	It~contains at~most two copies of~any path edge, since every $w_i>1/3$.
	Assign a~single colored copy in~the bundle $v_iv_{i+1}$ to~$v_{i+1}$; if~two copies are colored, assign one to~each endpoint.
	No~vertex can receive two assignments from path edges: that would require a~copy from the bundle on~its left and two from the bundle on~its right, whose three weights exceed one.
	Nor can a~vertex receiving such an~assignment have a~colored leaf edge.
	A~path edge from the left has weight $w_{i-1}>w_i$ and exceeds the space $w_i$ left by~its leaf edge; an~assignment from the right means two $w_i$ edges are already present.
	At~$v_m$ the leaf weight is~one.
	Finally assign each colored leaf edge to~its path vertex.
	Each leaf bundle contributes at~most one edge, as its weights exceed $1/2$.
	All assignments are distinct among $v_0,\ldots,v_m$, proving the claim.

	There are $2N(2m+1)$ edges in~total, so
	$\cw(G_N)\ge
	\left\lceil\frac{2N(2m+1)}{m+1}\right\rceil
	=\left\lceil4N-\frac{2N}{2N+1}\right\rceil=4N$.
	For the reverse inequality, use $2N$ colors on~the leaf bundles, putting one copy of~every leaf bundle in~each color.
	Use another $2N$ colors on~the path bundles in~the same way.
	At~each internal vertex the two path edges in~a~color weigh less than one.
	Hence both palettes are feasible and $\cw(G_N)=4N$.
\end{proof}

\label{ov:lower-end}

\section{Improved Constants from the Contents of Bins}\label{sec:overview-constants}\label{ov:page-11}
In this section, we show that local bin packing bounds allow stronger guarantees  than total weight.

\subsection{Bipartite multigraphs: harmonic charge and compatible selections}\label{sec:overview-bipartite-bp}
First, we start with bipartite multigraphs and for the complete proofs see \Cref{sec:bp-bounds}.

We use harmonic rounding~\cite{LeeLee85} from bin packing: for
$1/(k+1)<w\le1/k$, we assign charge $a(w)=\frac1k$.
Any $k$ edges from this range on the same pair of endpoints form a feasible packet of charge one. Full packets can be properly colored. The remaining objects have bounded multiplicity after tiny parallel edges have first been grouped. The missing-color and rounding arguments then insert them, with their remaining charge paying for the space used. This yields
$\chi'_w(G)\le\max_v\sum_{e\ni v}a(w_e)+o(\bp(G))$.
The largest harmonic charge of a unit bin is
$h_\infty=1+\frac12+\frac16+\frac1{42}+\cdots \approx 1.691$,
where successive denominators satisfy $r_{j+1}=r_j(r_j+1)$. This is the classical harmonic bin-packing constant of Lee and Lee~\cite{LeeLee85}; see also Epstein~\cite{Epstein}. Thus harmonic charge alone gives $(h_\infty+o(1))\bp(G)$ colors. \Cref{lem:harmonic-load} proves the coloring bound, including the treatment of arbitrarily small weights.

\begin{lemma}[Harmonic-charge coloring; informal version of \Cref{lem:harmonic-load}]
\label{lem:harmonic-charge-coloring}
For a bipartite weighted multigraph $J$, let
$A(J) :=
\max_{v \in V(J)}
\sum_{e \ni v}
\frac{1}{\lfloor 1/w_e \rfloor}$.
If $W_J(v) \le n$ at every vertex, then
$\chi'_w(J) \le A(J) + o(n)$.
\end{lemma}

\begin{proof}[Proof sketch]
On each parallel pair, pack $k$ edges of harmonic class
$(1/(k+1),1/k]$ into a single feasible packet, then
properly color full packets using K\H{o}nig's theorem.
Only boundedly many non-full packets remain per pair
after first grouping tiny edges. The general missing-color lemma (\Cref{lem:weighted-moment-mix})
and the bounded-multiplicity version of the rounding
theorem (see \Cref{sec:bounded-rounding}) insert almost all residual weight into the
palette; a fresh $o(n)$ palette finishes the rest.
\end{proof}

See \Cref{thm:bp-harmonic-first} for the complete proof of $(h_\infty+o(1))\bp(G)$ colors.

\paragraph{A small selection illustrates the improvement.}
Take the two ranges $H=(1/2,5/8]$ and $S=(1/3,3/8]$.
These ranges are compatible: an $H$ edge fits with an $S$ edge, and two
$S$ edges fit together. Select roughly one sixth of each range at every
vertex. Splitting vertices and applying K\H{o}nig's theorem makes this
selection on whole edges, with an error of at most one per range at each vertex.
The selected graph has at most $n/6+O(1)$ heavy edges and $n/3+O(1)$ total edges at each vertex, where $n=\bp(G)$. Split a vertex into two copies, put the heavy incidences at one copy, and distribute the small incidences to balance their degrees. A proper coloring of the copied graph gives a feasible weighted coloring with $n/6+O(1)$ colors.
To account for the removed edges, define the discounted charge
$q(w)=\frac56a(w)$ for $w\in H\cup S$ and $q(w)=a(w)$ otherwise.
Every original unit bin has total $q$-charge at most $3/2$.
The balanced selection therefore leaves harmonic charge at most
$3n/2+O(1)$ at every vertex.
This proves the short $5/3$ bound; \Cref{thm:bp-five-thirds-short} gives the calculation.

For the $1.62$ bound, we make three selections from compatible
pairs of weight ranges, assigning each edge to at most one selection.
The selected groups together use $(7/30)n+O(1)$ colors, while the
remaining harmonic charge is at most $(104/75)n+O(1)$ at every vertex.
Applying the harmonic-charge coloring bound to the remainder gives
$7/30+104/75=81/50=1.62$ as the leading coefficient.
The three pairs of ranges and their selection fractions are listed in
\Cref{tab:bip-three-selections}.
The detailed proofs are given in  \Cref{sec:bipartite-improved}.

\subsection{General multigraphs: handling dense odd sets}
\label{ov:page-12}

We now turn to general multigraphs. As in the bipartite case, our starting point is to group parallel edges into feasible packets, each of total weight at most one. We then treat each packet as a single edge and seek a proper edge coloring of the resulting multigraph. Such a coloring immediately gives a feasible weighted coloring of all edges contained in the packets. The remaining edges will be handled using our rounding framework.
For omitted details, see \Cref{sec:general-bounds}.

\paragraph{Why do odd sets matter?}
Unlike bipartite graphs, a general multigraph may require substantially more colors than its maximum degree. The reason is that each color class is a matching. Thus, for any odd set $S$ of vertices, a single color can contain at most $(|S|-1)/2$ edges with both endpoints in $S$. Consequently, the number of colors must be at least
$\frac{2|E_H(S)|}{|S|-1}$,
where $H$ denotes the multigraph of packets. Edmonds' matching theorem~\cite{Edmonds65} shows that these odd-set bounds, together with the maximum-degree bound, completely characterize the fractional chromatic index. Moreover, Sanders and Steurer's theorem~\cite{SS} converts this fractional bound into a proper edge coloring with only a lower-order additive loss.

Fortunately, most odd sets are easy to handle: the degree bound handles all odd sets of size at least seven. Our harmonic-charge analysis bounds the number of packets incident to any vertex by approximately $1.7\bp(G)$. For an odd set $S$ of size $s$, counting the degrees of its vertices gives
$\frac{2|E_H(S)|}{s-1}
    \le \frac{s}{s-1}\Delta(H)$.
When $s\ge 7$, the right-hand side is at most
$\frac76(1.7+o(1))\bp(G)
    < 1.999\bp(G)$
for sufficiently large $\bp(G)$. Thus, \emph{only triangles and five-vertex sets can obstruct our desired bound}.

\paragraph{Measuring the cost of a parallel bundle.}
We need a measure of how many packets a bundle requires before deciding
which sets to contract. Grouping and rounding leave a bounded number of
types on each pair, where a type specifies an endpoint pair and a weight.
A bundle pattern is a multiset of these types with total weight at most
one. Allow fractional amounts of patterns, requiring their combined type
counts to reproduce the bundle. The minimum total amount is its
\emph{fractional bin-packing cost}. Rounding down the amounts produces
whole packets; the few remaining objects on each pair are handled later.

For any fixed $\delta>0$, preprocessing gives a common local bin count
$N\le(1+\delta)\bp(G)+O_\delta(1)$.
For the density calculations, divide every bundle cost by $N$
and denote the result by $t_{uv}$. The number of whole packets on pair
$uv$ is at most $Nt_{uv}$. The local packings give
$t_{uv}\le1$, incident cost at most $M=1.7$ at every vertex, and
internal cost at most two in every triangle
(\Cref{lem:general-two-bags,lem:subtwo-charge}). These bounds control
both the sets we contract and the edges that leave them.

\paragraph{Contracting the dense sets.}
Put $\kappa=1.999$. A triangle is \emph{dense} if its internal cost
exceeds $\kappa$; a five-vertex set is dense if its internal cost exceeds
$2\kappa$. These are exactly their respective odd-set thresholds for a
palette of $\kappa N$ colors. We select all dense five-vertex sets first,
then the dense triangles disjoint from them.

Why can these contractions be made together? If two dense five-vertex
sets overlapped, their combined internal cost would exceed what the
degree and triangle bounds allow. A dense five-vertex set also has so
much cost inside it that little cost can leave it. This rules out a
dense triangle crossing its boundary. The remaining dense triangles
are disjoint for the same degree-counting reason.
\Cref{lem:subtwo-blocks} checks these claims and shows that contracting
the selected sets leaves every odd-set bound satisfied. After discarding
internal packets, Edmonds' and Sanders--Steurer's theorems therefore give
a proper coloring of the quotient packets with $(\kappa+o(1))N$ colors.

\paragraph{Reusing external colors inside a block.}
We next insert the remaining {\em external} objects and restore the internal
edges of each contracted block. The key is to remember which original
vertex an external color uses. In a proper packet coloring, a color
appears on at most one packet incident to a block, so it occupies just
one original vertex there. When inserting the remaining external
objects, the labelled rounding lemma (\Cref{lem:bp-labelled-round})
preserves the single-endpoint property and bounds the number of colors
with each label.

For example, if a color is used externally only at vertex $v$ of a
five-vertex block, it remains available for internal edges on the other
four vertices. An \emph{internal pattern} is a multiset of internal
edge types whose incident weight is at most one at every vertex.
We may assign such a pattern to a color labeled $v$ provided it avoids
$v$. Colors unused externally have no such restriction.
Thus restoration asks for a packing into one unrestricted class of
colors and several classes, each required to avoid a specified vertex.
The local charge estimate supplies a budget for each class: at an
original vertex, the internal bundle cost plus the external charge is
at most $1.7N$. Labelled rounding bounds the number of colors with that
vertex's label by its external charge, up to small rounding errors.

Triangles admit a direct construction using their local star packings.
For five-vertex blocks, we prove that the restricted classes still have
enough capacity by testing prices. If the internal demand cannot fit
fractionally into the available color classes, linear programming
duality supplies prices for which the demand costs more than those
classes can carry. We show that no such prices exist.
Give nonnegative prices to the
internal types and normalize the maximum feasible internal-pattern price
to one. The five-vertex price inequality
(\Cref{lem:subtwo-five-local}) bounds the total demand price by $1.95N$.
The palette lemma (\Cref{lem:subtwo-five-palette}) also accounts for
the loss from requiring some patterns to avoid a vertex: a color labeled
$v$ contributes the largest price of a pattern avoiding $v$, rather
than the unrestricted value one. For classes within the endpoint budgets,
it shows that the combined available price of the $\kappa N$ colors
covers the demand for every price vector.
Linear programming duality then gives a fractional packing respecting
all endpoint labels.

Finally, set aside the boundedly many tiny remainder objects in each
block. Each block has at most five vertices and, for fixed accuracy,
boundedly many remaining types, all of weight bounded away from zero.
Rounding down its fractional pattern amounts therefore leaves only
boundedly many objects to color separately.
These internal leftovers and the patterns moved out of excess label
classes use fresh colors. Since the blocks are disjoint, they can
reuse the same fresh palette. The remaining external objects and
preprocessing set-asides are colored using the global greedy cleanup
bound. Taking the accuracy parameters sufficiently small gives the total
$(1.999+o(1))\bp(G)$ bound.

The intermediate $\tfrac54 h_\infty$ bound in \Cref{sec:general-bp}
develops the packet construction; 
\Cref{sec:label-rounding} gives the rounding theorem needed for the multigraph setting; 
 \Cref{sec:general-subtwo} gives the final 1.999 bound (in particular, \Cref{sec:subtwo-proof} gives the proof); and \Cref{sec:label-algorithm}
gives the polynomial-time implementation.

\paragraph{The K\H{o}nig--Shannon comparison.}
Let $c_{\rm bip},c_{\rm gen}$ be the optimal uniform asymptotic coefficients relative to $\bp(G)$.
Then the bounds $c_{\rm gen}\le1.999$, $c_{\rm bip}\ge4/3$ give
\begin{equation}\label{eq:bp-konig-shannon-gap}
 1\le\frac{c_{\rm gen}}{c_{\rm bip}}\le\frac{1.999}{4/3}<\frac32.
\end{equation}
The individual coefficients remain in the intervals $[4/3,1.62]$ and $[3/2,1.999]$, respectively. For the total-weight parameter, the optimal coefficients are exactly two and three, retaining the ratio $3/2$.
It remains an intriguing open question whether the optimal asymptotic coefficients for bipartite and general multigraphs coincide, i.e., whether ${c_{\rm gen}}={c_{\rm bip}}$.

\label{ov:constants-end}

\paragraph{Acknowledgment.}
The second author would like to thank Mohit Singh for introducing him to the problem and for the inspiring initial discussions. He also thanks Debarsho Sannyasi, Pradyot Mohanty, and Shreyans Nagori for their contributions to the problem during their internships at the Indian Institute of Science.

The second author also thanks Google India Research Award, Ittiam Systems CSR grant, and the Walmart Center for Tech Excellence at IISc (CSR Grant WMGT-23-0001) for supporting the research.

\label{overview:end}
\newpage
\section*{AI Methodology}
\phantomsection\label{sec:ai-methodology}\label{sec:ai-overview}
\addcontentsline{toc}{section}{AI Methodology}
\begin{flushleft}
\renewcommand{\arraystretch}{1.3}
\begin{tabular}{@{}l@{\qquad}c@{}}
\textbf{Activity} & \textbf{AI used} \\
Asking the research question & No \\
Coming up with the approach & Yes \\
Proof development & Yes \\
Writing and exposition & Yes \\
Checking for bugs & Yes \\
Other supporting tasks & Yes \\
\end{tabular}
\end{flushleft}

\paragraph{Details.}
The research questions were chosen entirely by the authors.

The lower-bound constructions and the counterexample to the exact total-weight conjecture stated by Feige and Singh arose through iterative discussions between the authors and GPT-5.5.

The authors subsequently proposed ideas for a simpler randomized-rounding argument and developed these using the nibble method, obtaining a short proof of $\chi'_w(G)\le(2+o(1))\bp(G)$ for simple graphs. Further refinement with GPT Astra, with limited additional mathematical input from the authors during this stage, led to the sharper bound $\chi'_w(G)\le(1+o(1))\bp(G)$.

The authors suggested using harmonic bounds from bin packing to improve the guarantees for multigraphs. Through further discussions with the models, this idea was developed into harmonic-charge estimates connecting local bin contents to the global coloring cost. The main remaining obstacle was obtaining a coefficient below two relative to $\bp(G)$ for general multigraphs. The argument overcoming this obstacle was developed over several days of iterative proof search with GPT-6 in Codex. It combines contraction of dense vertex sets with control of endpoint labels, allowing internal edges to reuse colors assigned to external edges.

The initial drafts presenting all these results exceeded 200 pages. The authors subsequently simplified several arguments substantially and reorganized the presentation. Their contributions included proposing mathematical ideas, guiding the proof search, checking the resulting arguments, and developing shorter and clearer proofs.

ChatGPT also assisted with exposition, checking for bugs, exact-arithmetic checks, and minor corrections. The authors take responsibility for the final statements, proofs, and presentation.

\phantomsection\label{sec:verification-scope}
The companion repository contains AI-assisted Lean~4 and Mathlib
formalizations of the five upper-bound statements in \Cref{thm:main}:
\begin{center}
\url{https://github.com/rishigajjala/weighted-edge-coloring-challenges}
\end{center}
The running-time analyses are not formalized.

\clearpage
\bibliographystyle{alphaurl}
\bibliography{refs}

\clearpage
\appendix
\section*{Guide to the Appendix}
\phantomsection\label{sec:appendix-guide}
\addcontentsline{toc}{section}{Guide to the Appendix}
The table maps each main-text section to its full upper-bound proofs and
implementation details. Shared tools are collected in Appendix~\ref{sec:elementary};
both complete lower-bound proofs are in Section~\ref{sec:lower-bounds} of the main text.
Statements restated from the main text retain their original numbers.

\begin{center}
\small
\begin{tabular}{@{}p{.11\textwidth}p{.12\textwidth}p{.64\textwidth}r@{}}
\toprule
Appendix & Main text & Contents & Page\\
\midrule
\ref{sec:elementary} & --- & Preliminaries & \pageref{sec:elementary}\\[5pt]
\ref{sec:general-multi} & \ref{sec:overview-three} & The exact $3\lceil\Delta_w\rceil$ bound & \pageref{sec:general-multi}\\[5pt]
\ref{sec:round-overview} & \ref{sec:overview-simple} & The simple-graph bound and the full rounding proof & \pageref{sec:round-overview}\\[5pt]
\ref{sec:bip-multi} & \ref{sec:overview-bipartite} & The tight $2\Delta_w$ coefficient and missing-color coordination & \pageref{sec:bip-multi}\\[5pt]
\ref{sec:bipartite-bounds} & \ref{sec:overview-bipartite-bp} & Bipartite bin-packing bounds, through $1.62$ & \pageref{sec:bipartite-bounds}\\[5pt]
\ref{sec:general-bounds} & \ref{ov:page-12} & General bin-packing bounds, through $1.999$ & \pageref{sec:general-bounds}\\
\bottomrule
\end{tabular}
\end{center}
The \hyperref[sec:ai-methodology]{AI Methodology} statement and references follow
the main text.

\clearpage
\section{Preliminaries}\label{sec:elementary}\label{full-proofs:start}
We describe standard packing and rounding tools, together with the variants
needed in our proofs. We give short proofs when the precise bounds or
algorithmic details are useful.
As in the main text, all graphs are finite and loopless, parallel edges
are counted separately, and every edge weight lies in $(0,1]$.
For a vertex $v$, $W_G(v)=\sum_{e\ni v}w_e$; $\Delta_w(G)$ is the
maximum of these sums, and $\bp(G)$ is the largest local bin-packing optimum.
A feasible color has incident weight at most one at every vertex.

\subsection{Input Model and Algorithmic Guarantees}
\label{sec:algorithm-scope}

We~now specify the computational guarantees of~the constructions.
The input lists every edge separately and gives its weight as~a~rational number in~binary.
Every coloring algorithm returns the empty coloring when its input graph has no edges.
All asymptotic upper bounds, including the general $1.999$ bound, have randomized implementations with expected polynomial running time for each fixed accuracy.
The exact $3\lceil\Delta_w(G)\rceil$ construction is~deterministic and polynomial, and has no~accuracy parameter.
For the $1.999$ construction, \Cref{sec:label-algorithm} implements the endpoint-label step and the restoration of contracted blocks.

The algorithms need only approximate local bin packings.
For a~fixed $\xi>0$, the number of~bins used at~each vertex is~at~most $(1+\xi)$ times its optimum plus a~constant depending on~$\xi$.
Choosing $\xi$ sufficiently small absorbs this loss into the desired accuracy.
The subsequent rounding uses computable second moments, partial resampling, and a~matching among three-element options.
The implementations appear with their proofs: \Cref{sec:rounding-implementation,sec:mixing-implementation,sec:harmonic-inputs,sec:label-algorithm}.

All thresholds and constants depend on~the chosen accuracy before the graph is~given.
Below the relevant threshold, greedy cleanup uses at~most an~accuracy-dependent constant number of~colors.
To~recover the stated $(c+o(1))n$ form from guarantees $(c+1/k)n+C_k$, choose increasing thresholds $n_k$ with $C_k/n\le1/k$ for $n\ge n_k$.
Using the $k$th construction between successive thresholds makes the excess over $cn$ at~most $2n/k$, and hence $o(n)$, uniformly over the graph class.
The running-time assertion remains one for each fixed accuracy; it~does not assert polynomial dependence on~the reciprocal of~the accuracy.

\subsection{Completing Partial Colorings}

We begin with an elementary greedy counting bound.
\begin{lemma}\label{lem:general-cleanup}
Every finite loopless weighted multigraph $J$ with maximum weighted
degree at most $R$ satisfies $\chi'_w(J)\le\lfloor4R\rfloor+1$.
\end{lemma}
\begin{proof}
Process the edges in nonincreasing order of weight. For an edge of
weight $b\le1/2$, a color forbidden at an endpoint has load greater
than $1-b\ge1/2$. For $b>1/2$, a forbidden color is nonempty, and
its already colored edges each have weight at least $b>1/2$.
Thus fewer than $2R$ colors are forbidden at each endpoint. A palette
of $\lfloor4R\rfloor+1$ colors always has a color available at both.
\end{proof}

K\H{o}nig's theorem~\cite{Konig16} gives a~proper edge coloring of~a~bipartite multigraph with as~many colors as~its maximum degree.
The first bound below is the bounded-weight guarantee of Melen and
Turner's cap algorithm~\cite[Section~4]{MT89}; see also~\cite[Section~1.1]{FS}.
We give its vertex-splitting proof, then separate the heavy and light
edges to obtain the second bound.

\Needspace{13\baselineskip}
\begin{lemma}\label{lem:cleanup}
If $G$ is bipartite, every edge has weight at most $b<1$, and
$W_G(v)\le R$ for every vertex, then
\[
 \chi'_w(G)\le\left\lceil\frac{R}{1-b}\right\rceil.
\]
Without the bound $w_e\le b$, at most $2\lceil2R\rceil$ colors suffice.
\end{lemma}
\begin{proof}
There are no edges if $R=0$. Otherwise put $q=\lceil R/(1-b)\rceil$.
At each vertex list its edges in decreasing order of weight and cut the
list into consecutive blocks of $q$ edges; the last block may be shorter.
Replace the vertex by one new vertex for each block. An original edge now
joins the two new vertices representing its blocks at its two ends.
The new graph is bipartite and has maximum degree at most $q$. Give it
a proper $q$-edge-coloring and transfer the colors back.

At an original vertex with sorted weights $a_1\ge a_2\ge\cdots$,
a color appears at most once in each block. Its load is at most
\[
 a_1+a_{q+1}+a_{2q+1}+\cdots
 \le b+\frac1q\sum_j a_j
 \le b+\frac Rq\le1.
\]
For the middle inequality, charge each term after $a_1$ to the $q$
no-smaller weights in the preceding block. Thus the coloring is feasible.

For the last assertion, properly color the edges heavier than $1/2$
with at most $\lceil2R\rceil$ colors. Apply the first assertion, with
$b=1/2$, to the other edges, using another $\lceil2R\rceil$ colors.
\end{proof}

\subsection{Choosing Whole Options}
The next lemma rounds one option for each job while controlling every capacity.
In~our application, jobs are edges, and a~color option uses capacity at~the two endpoints.
The proof uses the iterative-relaxation method for sparse side constraints
as in~\cite{KLS08}: keep the choice equations and drop a capacity equation
once its remaining error is small. We give the counting argument for the
weighted, one-sided error $2\tau$ needed here.

\begin{lemma}\label{lem:two}
Suppose each job has a probability distribution over finitely many options.
Each option uses a nonnegative amount from at most two capacity constraints,
and each amount is at most $\tau$. If $b_R$ is the expected amount used
in constraint $R$, then there is a choice of one option per job using
at most $b_R+2\tau$ in every constraint.
\end{lemma}
\begin{proof}
Let $x_j$ be the current value of option $j$. Values zero and one will
remain fixed. Initially the values are the given probabilities.
Keep the equation that the options for each unfinished job sum to one.
For a capacity constraint, also keep its equation fixing the current load
until the first time
\[
 \sum_{j\text{ using that constraint},\ 0<x_j<1}(1-x_j)\le2.
\]
At that time drop the equation permanently. The sum measures how much the
remaining option values could increase if all of them were rounded to one.

Suppose $N>0$ option values are still strictly between zero and one, and
$r$ jobs are unfinished. Each such job has at least two of these values,
so $r\le N/2$. Their values sum to $r$, hence $\sum_j(1-x_j)=N-r$.
If $s$ capacity equations are still kept, counting each option in at most
two constraints gives $s<N-r$. Thus fewer than $N$ linear equations are
being kept for $N$ variables. Linear algebra gives a nonzero direction
in which all these equations stay unchanged. Move in that direction until
one variable first reaches zero or one. Repeat; each step fixes another
variable, so the process ends.

When a capacity equation is dropped, its load is still $b_R$. Its possible
increase afterward is at most $\tau\sum_j(1-x_j)\le2\tau$. A constraint
kept to the end has load exactly $b_R$. Every job has exactly one option
with final value one, which proves the result.
\end{proof}

\Needspace{20\baselineskip}
\subsection{Two Selection Theorems}\label{sec:inputs}
The rounding argument must control many local quantities at~once, with error bounds independent of~the number of~vertices.
The theorem of~Harris and Srinivasan controls the sums arising from independent choices.
The theorem of~Delcourt and Postle then selects disjoint options for the large edges.

\begin{proposition}[Controlling many sums, Harris--Srinivasan~{\cite[Theorem~4.7]{HS}}]\label{prop:HS}
Let $Z_i$ be independent random choices, each with a finite set of possible
values. For each index $t$ let
\[
 F_t=\sum_i f_{t,i}(Z_i),\qquad 0\le f_{t,i}(z)\le1.
\]
Define
\[
 R=\max_t\E F_t,\qquad
 D=\max_{i,z}\sum_t f_{t,i}(z),\qquad C=\max\{2,R,D\}.
\]
There is one choice of all the $Z_i$ for which every sum satisfies
\[
 F_t\le\E F_t+
 O\bigl(\sqrt{R\log(C+2)}+\log(C+2)\bigr).
\]
The constant hidden in $O(\cdot)$ is absolute.
\end{proposition}

Here $R$ bounds the mean of~one sum, whereas $D$ bounds the combined contribution of~one choice to~all sums.
Neither parameter counts the number of~sums.
This distinction allows the later rounding errors to~be independent of~the order of~the graph.

For completeness, this form follows from Theorem~4.7 of~\cite{HS} by
using upper bounds $c_t=\max\{1,\E F_t\}$, the parameter
$\widehat D=C+1$, and its error parameter $1/\widehat D$.
If $c_t\le\log\widehat D$, that theorem gives $O(\log(C+2))$.
Otherwise it gives
$c_t+O(\sqrt{c_t\log(C+2)}+1)$. Since $c_t\le R+1$, both cases
give the displayed bound.

\begin{proposition}[Choosing disjoint sets, Delcourt--Postle]\label{prop:DP}
Let $A$ and $B$ be disjoint sets of elements. Consider a finite list of
sets, each containing exactly one element of $A$ and at most two elements
of $B$. The same set may occur several times in the list. Write $d(x)$
for the number of listed sets containing $x$, and $\Delta_2$ for the
largest number of listed sets containing any given pair of distinct
elements. All counts include repetitions.

Fix $\zeta>0$. There is a constant $a>0$ such that, for all sufficiently
large $D$, the conditions
\[
 d(x)\ge(1+D^{-a})D\quad(x\in A),\qquad
 d(y)\le D\quad(y\in B),\qquad
 \Delta_2\le D^{1-\zeta}
\]
guarantee a collection of pairwise disjoint listed sets covering every
element of $A$.
\end{proposition}

An~element of~$A$ will represent an~edge to~be colored, and an~element of~$B$ will represent an~available slot at~one of~its endpoints.
A~listed set specifies one option for that edge.
A~matching therefore chooses each edge once and uses each slot at~most once.
The degree inequalities provide more options for each edge than compete for a~single slot, while the pair-count bound limits repeated interactions.
The proposition is~\cite[Theorem~2.6]{DP} for sets of~size at~most three and an~empty family of~forbidden configurations.

\subsection{Rounding a Bin-Packing Linear Program}

For an item multiset $I$, write $\operatorname{bp}(I)$ for its minimum number of unit bins and $W(I)=\sum_{i\in I}w_i$ for its total weight.
Let $\operatorname{LP}(I)$ be the minimum mass of integer bin patterns covering its type counts. Types may be
repeated within a pattern. Each actual item can instead be treated as
a separate type; this convention is useful in the following rank
rounding argument. Covering and equality versions have the same
optimum, because surplus items can be deleted.
The next lemma combines the linear grouping of
Fernandez de la Vega and Lueker~\cite{fernandez1981bin} with rounding a
basic solution of the configuration LP; see~\cite[Lemma~2, Corollary~1,
and Section~4]{KK82}.
\begin{lemma}\label{lem:bin-rounding}
Fix $0<\tau\le1$ and $\nu>0$. Every finite item multiset $I$ with
weights in $[\tau,1]$ satisfies
\[
 \operatorname{bp}(I)\le(1+\nu)\operatorname{LP}(I)+O_{\nu,\tau}(1).
\]
For fixed $\nu,\tau$ and rational weights, such a packing can be found
in polynomial time.
\end{lemma}
\begin{proof}
The empty instance is immediate. Otherwise let
$L=\operatorname{LP}(I)>0$. The total item weight is at most $L$, so
$|I|\le L/\tau$. Sort the items from largest to smallest and divide
them into consecutive groups of $g=\lceil\nu L\rceil$ items, with a
possibly smaller final group. Set aside the first group for singleton
bins. In each later group, round every weight up to that group's largest
weight.

Assign the rounded items of a group injectively to the original items
of the preceding group. No rounded item is larger than its assigned
original item, and these assignments are disjoint between groups.
In a fractional packing of $I$, replace every occurrence of an assigned
original type by its rounded type, and delete unassigned types. This
shows that the rounded remainder has LP value at most $L$.

There are at most $t=\lceil1/(\nu\tau)\rceil$ rounded types.
An optimal basic solution of their pattern LP has at most $t$ positive
variables, since there are only $t$ demand constraints. Rounding these
variables upward gives an integral pattern cover of cost at most $L+t$.
Delete surplus occurrences and restore the original, no-larger weights.
Including the first group gives
\[
 \operatorname{bp}(I)\le L+t+g\le(1+\nu)L+O_{\nu,\tau}(1).
\]
For the algorithmic statement, a pattern has at most
$\lfloor1/\tau\rfloor$ items. Enumerating all patterns and solving the
resulting LP is polynomial for fixed $\tau$. The grouping and the
second, bounded-type LP are polynomial as well.
\end{proof}

\subsection{Local Feasible Distributions}
\label{sec:feasible-choices}
The rounding theorem requires a~distribution on~feasible sets at~each vertex and color.
The next two lemmas construct such distributions from a~load bound or~from bounds on~the number of~items in~each size class.
The first uses a weighted variant of the probability-preserving pairwise
rounding step in~\cite{GKPS06}; the proof below establishes the required
weight bound.

\begin{lemma}\label{lem:half-load}
Suppose $0<w_i\le1/2$ for a finite set of items $I$. If
\[
 0\le x_i\le1\quad(i\in I),\qquad \sum_iw_ix_i\le1/2,
\]
then there is a random subset $S\subseteq I$ such that
$\sum_{i\in S}w_i\le1$ always and $\Prb(i\in S)=x_i$ for every $i$.
Equivalently, $x$ is a weighted average of the zero--one vectors of sets
of total weight at most one.
\end{lemma}
\begin{proof}
We turn the entries of $x$ into zeros and ones by random choices that
preserve the expected value of each entry. Call an entry \emph{fractional}
if it is strictly between zero and one.

Suppose two entries $x_i,x_j$ are fractional. Replace them temporarily by
\[
 x_i(t)=x_i+t w_j,\qquad x_j(t)=x_j-t w_i.
\]
Their combined weight stays fixed. The values of $t$ for which both
entries lie in $[0,1]$ form an interval $[-a,b]$, with $a,b>0$. At each
end of this interval, at least one of the two entries is zero or one.
Choose $t=b$ with probability $a/(a+b)$ and $t=-a$ with probability
$b/(a+b)$. Then $\E t=0$, so the expected value of each entry after this
choice is its value before the choice. All other entries stay fixed.

Repeat while at least two fractional entries remain. Each step reduces their
number, preserves total weight, and preserves all entry values in
expectation. Thus the process stops after at most $|I|$ steps, with total
weight still at most $1/2$ and with at most one fractional entry. If that
entry has value $q$, round it to one with probability $q$ and to zero
otherwise. This last choice also preserves its expected value. It adds
at most $\max_i w_i\le1/2$ to the total weight, so the final zero--one
vector always describes a set of weight at most one. Repeatedly taking
conditional expectations shows that the probability its $i$th entry is
one is the original value $x_i$, as required.
\end{proof}

The rounding within each class in the next lemma is the cardinality case
of dependent rounding~\cite[properties (P1)--(P2)]{GKPS06}, applied to a star.
We combine it with a random choice of class.
\begin{lemma}
\label{lem:bp-classchoice}
Partition a finite set of items into classes $I_2,\ldots,I_T$, with
every item in $I_k$ having weight at most $1/k$. Let $x_i\ge0$, and
put
\[
 p_k=\max\left\{\frac1k\sum_{i\in I_k}x_i,
                       \max_{i\in I_k}x_i\right\},
\]
where the maximum over an empty class is zero. If $\sum_{k=2}^Tp_k\le1$,
there is a random feasible set with inclusion probability $x_i$ for
every item $i$.
\end{lemma}
\begin{proof}
Choose class $k$ with probability $p_k$, and choose the empty set with
the remaining probability. When $p_k>0$, the conditional numbers
$z_i=x_i/p_k$ satisfy $0\le z_i\le1$ and $\sum_{i\in I_k}z_i\le k$.
Apply the pairwise procedure from the proof of \Cref{lem:half-load},
now moving two fractional entries by equal amounts in opposite directions.
This preserves their sum and each entry's expectation. After the last
fractional entry is rounded, the cardinality is at most
$\lceil\sum z_i\rceil\le k$, and the inclusion probabilities are the
original $z_i$.

All chosen items have weight at most $1/k$, so the chosen set is
feasible. Multiplication by the probability $p_k$ gives unconditional
inclusion probability $x_i$.
\end{proof}

\subsection{Constructing the Local Inputs}
\label{sec:local-input-algorithms}
The coloring algorithms begin with local packings and distributions on~feasible sets.
We~construct these inputs explicitly, using rational arithmetic throughout.

\paragraph*{Approximate local packings}
We use the classical asymptotic approximation scheme based on linear
grouping and greedy insertion of small items~\cite{fernandez1981bin,KK82}.
For every fixed $\xi>0$, a finite item multiset $I$ can be packed in
polynomial time into
\[
 (1+\xi)\operatorname{bp}(I)+O_\xi(1)
\]
bins. Choose fixed $\nu,\tau>0$ so that
$1+\nu\le1+\xi$ and $(1-\tau)^{-1}\le1+\xi$.
Apply \Cref{lem:bin-rounding} to the items of weight at least $\tau$.
Their LP value is at most $\operatorname{bp}(I)$, so this gives at most
$(1+\nu)\operatorname{bp}(I)+O_{\nu,\tau}(1)$ bins.
Insert the lighter items greedily. Whenever a new bin is needed,
every previous bin has load greater than $1-\tau$. Hence the final
number is at most
\[
 \max\left\{(1+\nu)\operatorname{bp}(I)+O_{\nu,\tau}(1),
                   \frac{W(I)}{1-\tau}+1\right\}.
\]
Since $W(I)\le\operatorname{bp}(I)$, this proves the claim.
A uniform choice among the resulting bins is an explicit distribution
on actual sets, giving each item probability exactly the reciprocal
of the number of bins. Empty bins may be added to equalize this number
between vertices.

The next lemma uses an algorithmic Carath\'eodory decomposition
\cite[Theorem~6.5.11]{GLS88}, followed by rounding the at most one
fractional coordinate of each vertex of the single-capacity polytope.
We give the details to record the support sizes and running time.
\begin{lemma}
\label{lem:algorithmic-local-load}
Suppose $0<w_i\le b\le r$, $0\le x_i\le1$, and
$\sum_{i=1}^m w_ix_i\le r-b$.
Then $x$ is a convex combination of at most $2m+2$ zero--one vectors, each of weight at most $r$.
Similarly, if $x\in[0,1]^m$ and $\sum_i x_i\le k$ for a nonnegative integer $k$, then $x$ is a convex combination of at most $m+1$ incidence vectors of sets with at most $k$ items.
For rational input data, both decompositions can be found in polynomial time.
\end{lemma}
\begin{proof}
Consider the bounded polytope
\[
 P=\{y\in[0,1]^m:\textstyle\sum_i w_iy_i\le r-b\}.
\]
We first express $x$ as a convex combination of at most $m+1$ vertices of $P$. At the current point, keep all
inequalities that hold with equality; they describe its smallest face.
Find a vertex $z$ of this face by successive linear optimization of
the coordinates. If the current point is not already $z$, extend the
ray from $z$ through that point until another inequality becomes tight.
The current point lies between $z$ and this new boundary point.
Record that convex combination and continue at the boundary point.
The dimension of the current face drops at every step, so there are
at most $m$ steps. The polytope has $2m+1$ inequalities.
For rational input data, the linear programs and ray intersections can be computed in polynomial time with polynomial encoding length.

A vertex of $P$ has at most one coordinate strictly between zero and
one. Otherwise two such coordinates could be moved slightly in opposite
weighted directions while staying in $P$, contradicting that it is a
vertex. For each vertex in the decomposition, replace its possible
fractional coordinate $y_i$ by one with probability $y_i$ and by zero
otherwise. This preserves its mean. Rounding up increases its weight
by at most $b$, so both resulting sets have weight at most $r$.
Each vertex gives at most two sets, proving the first claim.

For the second claim use
$\{y\in[0,1]^m:\sum_i y_i\le k\}$ instead.
The same two-coordinate argument leaves at most one fractional
coordinate at a vertex. If the sum constraint is not tight, that
coordinate can be moved alone; if it is tight, its value is an integer
because $k$ and the other coordinates are integers. Thus every vertex
is already the required zero--one vector.
\end{proof}

\clearpage
\section{The Exact Total-Weight Bound for General Multigraphs}
\label{sec:general-multi}
This appendix supplies the proof and deterministic algorithm for
Section~\ref{sec:overview-three}.

\subsection{Properly Coloring a Bounded-Degree Multigraph}
The following elementary bound supplies the proper coloring used in the
exact general total-weight bound. It is a slightly weaker form of
Shannon's bound~\cite{Shannon49}, with a short constructive proof.

\begin{lemma}\label{lem:three-per-two}
Every multigraph without loops and of maximum degree $d$ has a proper
edge coloring with at most $3\lceil d/2\rceil$ colors.
\end{lemma}
\begin{proof}
There is nothing to do if $d=0$. Put $r=\lceil d/2\rceil$.
Take two disjoint copies of the graph. For every original vertex $v$,
join its two copies by $2r-d(v)$ new parallel edges. Every vertex of
the enlarged graph now has degree $2r$.

Orient an Euler tour in each connected component. Each vertex has
exactly $r$ incoming and $r$ outgoing edges. Make two new vertices
$v_{\rm out},v_{\rm in}$ for each vertex $v$, and replace each directed
edge $u\to v$ by the edge $u_{\rm out}v_{\rm in}$. This graph is
bipartite and has maximum degree $r$. K\H{o}nig's theorem colors its
edges properly with $r$ colors.

Transfer this coloring to the enlarged graph. At a vertex, each color
now appears on at most one incoming edge and at most one outgoing edge.
The edges of that color therefore form paths and cycles, which can be
properly colored with three colors. Use a separate set of three colors
for each of the $r$ colors, and then restrict to the original graph.
\end{proof}

\Needspace{10\baselineskip}

\subsection{Finding Room for a Remaining Edge}\label{sec:common-room-proof}
\restatedstatement{\lemma}{\commonroomlemma}
\begin{proof}
	Let $W=\{c:x_c\ge a\}$, so $|W|\ge2n$.
	Assume that no~color has room at~both vertices.
	Then, for each color $c$, either $x_c>1-a$ or~$y_c>1-a$.

	Let
	\[
	S=\{c:y_c>1-a\},\qquad s=|S|.
	\]
	If~$S$ is~nonempty, then
	\[
	s(1-a)<\sum_{c\in S}y_c\le n.
	\]
	The same conclusion $s<n/(1-a)$ is~immediate when $S$ is~empty.
	Every color outside $S$ has $x_c>1-a$.

	First suppose that $s\le n$.
	There are at~least $2n$ colors outside $S$, so
	\[
	\sum_c x_c>(3n-s)(1-a)\ge2n(1-a)>n,
	\]
	a~contradiction.

	It~remains to~consider $s>n$.
	At~most $n$ of~all $3n$ colors lie outside $W$, since $|W|\ge2n$.
	Hence $|S\cap W|\ge s-n$.
	The sets $[3n]\setminus S$ and $S\cap W$ are disjoint, and therefore
	\begin{align*}
		\sum_c x_c
		&>(3n-s)(1-a)+(s-n)a.
	\end{align*}
	The expression on~the right has coefficient $2a-1<0$ in~$s$.
	Using $s<n/(1-a)$, it~is~therefore larger than
	\[
	\left(3n-\frac{n}{1-a}\right)(1-a)
	+\left(\frac{n}{1-a}-n\right)a
	=n+\frac{n(1-2a)^2}{1-a}>n.
	\]
	This again contradicts $\sum_c x_c\le n$.
\end{proof}

\subsection{The Exact Bound}\label{sec:general-total-proof}
\begin{theorem}\label{thm:general-total-weight}
Every weighted multigraph without loops satisfies
\[
\chi'_w(G)\le3\lceil\Delta_w(G)\rceil.
\]
Such a coloring can be found by a deterministic polynomial-time algorithm.
The leading coefficient three is best possible.
\end{theorem}
\begin{proof}
The empty graph needs no colors. Otherwise let
$n=\lceil\Delta_w(G)\rceil\ge1$.
	Order all edges from largest weight to~smallest.
	Starting with $F=\varnothing$, inspect the edges in~this order.
	Put the current edge into $F$ if~its addition leaves the degree of~$F$ at~most $2n$ at~both endpoints.
	Thus
	\[
	\Delta(F)\le2n.
	\]
	If~an~edge $e\notin F$ was rejected, then at~one of~its endpoints, say $u$, there were already $2n$ edges of~$F$.
	All of~them were inspected before $e$, so~each has weight at~least $w_e$.
	We~call them the witnesses for $e$ at~$u$.
	The witnesses remain in~$F$.

	By~Shannon's theorem, in~the form of~\Cref{lem:three-per-two}, the multigraph $F$ has a~proper edge coloring with at~most $3n$ colors.
	Pad it~with unused colors to~obtain a~palette of~exactly $3n$ colors.
	This is~already a~feasible weighted coloring of~$F$, since at~most one edge of~each color meets any vertex.

	We~insert the edges outside $F$ one at~a~time without changing the colors on~$F$.
	Consider an~edge $e=uv$ of~weight $a$.
	First, $a<1/2$.
	Indeed, at~a~witness endpoint the $2n$ witnesses and $e$ are distinct incident edges.
	If~$a\ge1/2$, their total weight would be~at~least $2n a+a>n$, contrary to~$W_G(u)\le n$.

	Let $x_c$ and $y_c$ be~the current loads of~color $c$ at~$u$ and $v$.
	Only a~subset of~the incident edges has been colored, so
	\[
	\sum_c x_c\le W_G(u)\le n,
	\qquad
	\sum_c y_c\le W_G(v)\le n.
	\]
	The $2n$ witness edges at~$u$ have distinct colors because the coloring of~$F$ is~proper.
	In~each of~those colors the present load at~$u$ is~at~least the weight of~its witness, and hence at~least $a$.
	\Cref{lem:three-common-room} now gives a~color with room for $e$ at~both endpoints.
	Give $e$ that color.
	The coloring stays feasible, and the colors of~all skeleton edges stay fixed, so~the same argument applies to~the next edge.
	After finitely many insertions every edge is~colored with one of~the same $3n$ colors.
	Therefore
	\[
	\chi'_w(G)\le3n=3\lceil\Delta_w(G)\rceil.
	\]

The triangle construction in \Cref{app:triangle-lower} proves sharpness.
\end{proof}

\Needspace{14\baselineskip}
\begin{algorithm}[The exact general total-weight construction]\label{alg:general-total}
Input: a weighted multigraph $G$ without loops. Put
$n=\lceil\Delta_w(G)\rceil$.
\begin{enumerate}[leftmargin=*,nosep]
\item Sort the edges by decreasing weight. Starting with $F=\varnothing$,
add each edge whenever doing so leaves the degree of $F$ at most $2n$
at both endpoints.
\item Apply the constructive proof of \Cref{lem:three-per-two} to color
$F$ properly with at most $3n$ colors. Add unused colors to obtain a
palette of exactly $3n$ colors.
\item Inspect each edge outside $F$ and assign it a color with room at
both endpoints. Find such a color by scanning the palette.
Return the resulting coloring.
\end{enumerate}
\Cref{lem:three-common-room} guarantees that every scan succeeds.
The construction is deterministic and polynomial and uses at most
$3\lceil\Delta_w(G)\rceil$ colors.
\end{algorithm}

\subsection{Sharpness}\label{app:triangle-lower}
For an~integer $t\ge2$, take three vertices and put $t$ parallel edges between each pair.
Give every edge weight $(t+1)/(2t)>1/2$.
Any two edges meet, and their weights add up~to~more than one.
Thus every color contains at~most one edge, and
\begin{equation}
	\chi'_w(G_t)=3t,\qquad
	\Delta_w(G_t)=t+1,\qquad \bp(G_t)=2t.
	\label{eq:general-triangle}
\end{equation}
This example is~also noted by~Feige and Singh~\cite[Section~1.2]{FS}.
It~shows that $\chi'_w(G_t)/\Delta_w(G_t)\to3$ and gives the lower bound $3/2$ for the coefficient of~$\bp(G)$ in~general multigraphs.

\begin{figure}[htbp]
	\centering
\begin{tikzpicture}[x=1cm,y=1cm,font=\small]
\node[vertex,label=below left:$a$] (a) at (0,0) {};
\node[vertex,label=below right:$b$] (b) at (4,0) {};
\node[vertex,label=above:$c$] (c) at (2,2.6) {};
\foreach \bend in {-9,0,9}{
 \draw[inkblue,thick] (a) to[bend left=\bend] (b);
 \draw[inkblue,thick] (b) to[bend left=\bend] (c);
 \draw[inkblue,thick] (c) to[bend left=\bend] (a);
}
\node[below=9pt] at (2,0) {$t$ parallel edges};
\node[rotate=52,above=7pt] at (1,1.3) {$t$ parallel edges};
\node[rotate=-52,above=7pt] at (3,1.3) {$t$ parallel edges};
\node[align=center] at (2,.85) {Each weight\\$\dfrac{t+1}{2t}>\dfrac12$};
\node[align=left,anchor=west] at (5.1,1.4)
 {Every pair of edges meets.\\[5pt]
  Every color holds one edge.\\[8pt]
  $\chi'_w=3t,\quad\Delta_w=t+1,\quad\bp=2t.$};
\end{tikzpicture}
	\caption{The obstruction for general multigraphs. Each set of three
	drawn curves represents $t$ parallel edges. The required number of
	colors approaches three times the maximum total weight, while it is
	exactly $3/2$ times the local bin-packing bound.}
	\label{fig:general-triangle}
\end{figure}
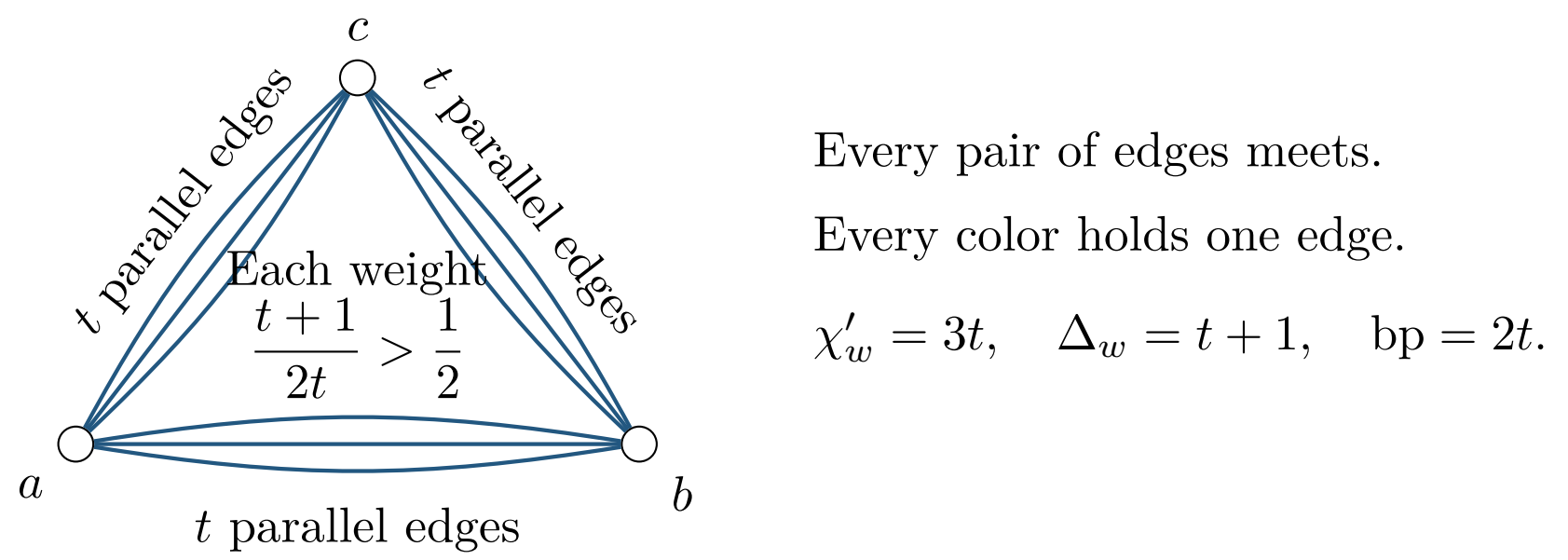

\clearpage
\section{The Simple-Graph Bound and Local-to-Global Rounding}
\label{sec:round-overview}
This appendix gives the full rounding argument used in
Section~\ref{sec:overview-simple}. We first describe local bins by profiles,
then prove the rounding theorem, and finally apply it to the local bin packings.
The remaining subsections give the implementation details.

For a~vertex~$v$, put
\[
\mathcal K_v=\operatorname{conv}
\left\{\one_S:S\subseteq E(v),\ \sum_{e\in S}w_e\le1\right\}.
\]
Thus, $x\in\mathcal K_v$ if~and only if~there is~a~random feasible subset of~$E(v)$ that contains each edge~$e$ with probability~$x_e$.
The rounding theorem uses this local feasibility condition.

\noindent\textbf{Remark.}
	Two edges of~weight $0.6$, each given fractional value $0.8$, have total fractional load $0.96$.
	However, at~most one of~them fits in~a~bin, so~their inclusion probabilities must sum to~at~most one.
	The load constraint alone therefore does not imply membership in~$\mathcal K_v$.

\subsection{Describing Local Bin Packings}\label{sec:local}
We first consider only one vertex, or equivalently one finite set of items
with weights $w_e\le1$. A subset is feasible if its total weight is at
most one. We say a vector $x=(x_e)$ is an \emph{average of feasible sets}
if there are nonnegative numbers $p_S$ such that
\[
 \sum_Sp_S=1,\qquad x_e=\sum_{S\ni e}p_S,
\]
where only feasible sets $S$ are used. In probability terms, choose $S$
with probability $p_S$; then item $e$ is chosen with probability $x_e$.
Write $\mathcal K$ for all such vectors. At a graph vertex we write $\mathcal K_v$.

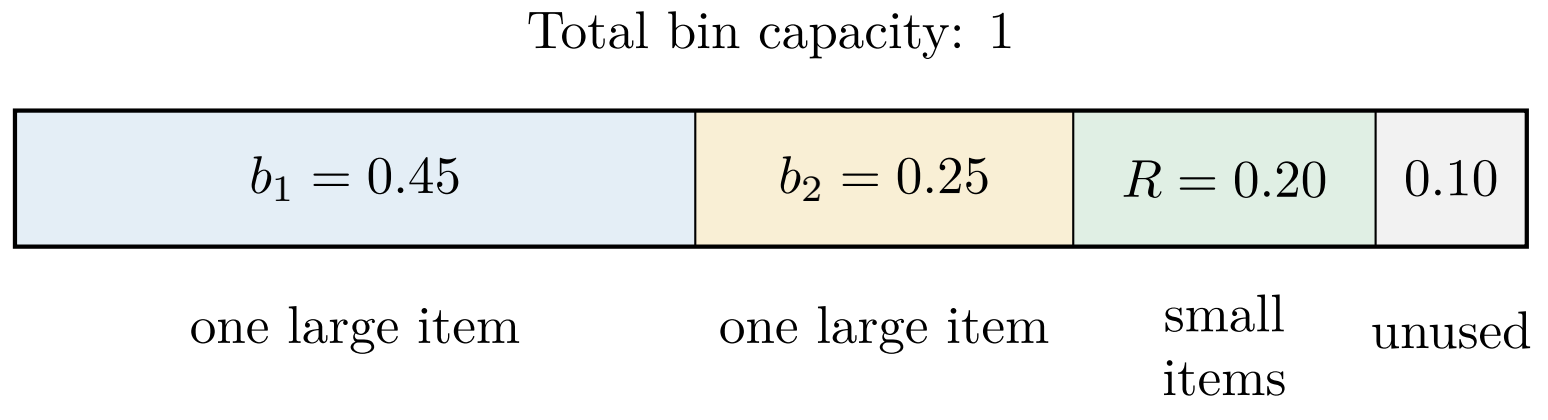
\begin{figure}[htbp]
	\centering
\begin{tikzpicture}[x=1cm,y=1cm,font=\small]
\fill[softblue] (0,0) rectangle (4.5,.9);
\fill[softgold] (4.5,0) rectangle (7,.9);
\fill[softgreen] (7,0) rectangle (9,.9);
\fill[black!5] (9,0) rectangle (10,.9);
\draw[thick] (0,0) rectangle (10,.9);
\foreach \xx in {4.5,7,9}{\draw (\xx,0)--(\xx,.9);}
\node at (2.25,.45) {$b_1=0.45$};
\node at (5.75,.45) {$b_2=0.25$};
\node at (8,.45) {$R=0.20$};
\node at (9.5,.45) {$0.10$};
\node[align=center] at (2.25,-.55) {one large item};
\node[align=center] at (5.75,-.55) {one large item};
\node[align=center] at (8,-.65) {small\\items};
\node[align=center] at (9.5,-.55) {unused};
\node at (5,1.4) {Total bin capacity: $1$};
\end{tikzpicture}
	\caption{A profile records places for large items and space for
	small items. Each large-item slot holds at most one item whose weight
	fits. Several small items may share the reserved space. Unused capacity
	is allowed.}
	\label{fig:bin-description}
\end{figure}

Fix a small positive number $\beta$. Call an item \emph{large} if its
weight is at least $\beta$, and \emph{small} otherwise. A feasible set
has at most $K=\lfloor1/\beta\rfloor$ large items. We describe a bin by
slots of sizes $b_1,\ldots,b_K$, allowing zero-size slots, and by a
space of size $R$ reserved for small items. The description is valid if
\[
 \sum_i b_i+R\le1.
\]
A large item may occupy a slot only if its weight is no greater than
the slot size. Several small items may share the reserved space as long
as their total weight is at most $R$. See \Cref{fig:bin-description}.

The number of~different bins may depend on~the input size.
We~replace them by~a~bounded list of~profiles, losing only a~small fraction of~each item's inclusion probability.
At~this stage, an~item may be~assigned fractionally to~several profiles; the later rounding argument makes the choices integral.

\begin{lemma}\label{lem:profiles}
Fix $0<\delta<1$ and $0<\beta\le r_0=\delta/4$. If $x\in\mathcal K$,
there are at most $M=M(\beta,\delta)$ valid profiles, with
probabilities $\lambda_P\ge0$ satisfying $\sum_P\lambda_P=1$, and nonnegative assigned fractions $a(e;P,i)$ with the following properties.
Put $y_e=(1-\delta)x_e$. Then
\begin{align}
 \sum_{P,i}a(e;P,i)&=y_e,\label{eq:marginal}\\
 \sum_e a(e;P,i)&\le\lambda_P
       &&\text{for each large-item slot }i,\label{eq:slot}\\
 \sum_{e:w_e<\beta}w_ea(e;P,R)&\le\lambda_PR
       &&\text{for the space reserved for small items}.\label{eq:reservoir}
\end{align}
Large items are assigned only to slots that can hold their weight. Small
items are assigned only to positive reserved spaces, each of which has
size at least $r_0$. Every individual assigned fraction is at most $y_e$.
\end{lemma}

\noindent\textbf{Remark.}
	Equation~\eqref{eq:marginal} preserves a~$1-\delta$ fraction of~each item's original probability.
	The capacity bounds use two different units: a~large-item slot holds one item, whereas a~reserved space for small items bounds their total weight.
	Keeping these two constraints separate will allow us~to round large and small edges by~different methods.

\begin{proof}
Choose a random feasible set whose inclusion probabilities are $x_e$.
Give each large item in the chosen set a slot equal to its weight, and
put $R=1-$ the sum of those slot sizes. The small items in the set fit
in this remaining space.

For each distinct positive slot size $b$, let $m(b)$ be the expected
number of slots of that size. This is a finite list of nonnegative
numbers with $\sum_bm(b)\le K$. Put $q=\delta/4$. We keep a short
list of slot sizes as follows. Keep the largest size. Scan toward
smaller sizes, and keep the first size for which the sum of the $m(b)$
at smaller sizes than the previous kept size reaches $q$. Repeat, and keep the
smallest size at the end. If there are no slots, the list is empty.

Each completed step uses at least $q$ of the total sum $\sum_bm(b)$.
Thus at most $L=\lfloor K/q\rfloor+2$ sizes are kept. The sum of the
$m(b)$ strictly between two consecutive kept sizes is less than $q$.
Round every slot size down to the largest kept size no greater than it.
For any threshold $t$, this lowers the expected number of slots of size
at least $t$ by less than $q$: only sizes between the two kept values
around $t$ can cross below $t$.

Also round each reserved space down to $r_0\lfloor R/r_0\rfloor$.
Every resulting profile is valid, and every positive reserved
space has size at least $r_0$. Sort the slots and pad with zeros.
There are now only boundedly many descriptions, since there are at most
$L+1$ choices for each of $K$ slot sizes and at most
$\lfloor1/r_0\rfloor+1$ choices for the reserved space.

Give this rounded distribution total probability $s=1-\delta$.
Add probability $sq$ on a profile consisting of one slot of size
one, and probability $sr_0$ on a profile consisting of reserved
space of size one. These probabilities sum to at most
\[
 s(1+q+r_0)=(1-\delta)(1+\delta/2)\le1.
\]
Put any probability left over on the empty description. In particular,
one may use
\[
 M\le(L+1)^K\bigl(\lfloor1/r_0\rfloor+1\bigr)+3.
\]

The extra size-one slot repairs every lost threshold count. The total
slot capacity, counted as fractions of slots, for sizes at least $t$
is at least
\[
 s\sum_{b\ge t}m(b)
 =\sum_{e:w_e\ge\max\{t,\beta\}}sx_e
 \qquad(0<t\le1).
\]
Assign each large item its required fraction $sx_e$, processing items
from largest weight to smallest. Split fractions across compatible
slots as needed. This process cannot fail. If it stopped at an item
of weight $t$, all fractions already placed came from items of weight
at least $t$, and all slots of size at least $t$ would be full. That
would contradict the displayed inequality. This proves
\eqref{eq:marginal} and \eqref{eq:slot} for large items.

The total expected reserved space is
$Q=\sum_P\lambda_PR(P)$. Rounding a space down loses less than $r_0$,
and the added full-space description repairs that loss. Hence
\[
 Q\ge s(\E R-r_0)+sr_0
 \ge\sum_{e:w_e<\beta}w_e y_e.
\]
Assign a small item's fraction among the reserved spaces in proportion
to their expected capacities:
\[
 a(e;P,R)=y_e\frac{\lambda_PR(P)}{Q}\quad\text{if }Q>0.
\]
These fractions sum to $y_e$. Their total weight in description $P$
is at most $\lambda_PR(P)$, giving \eqref{eq:reservoir}. If $Q=0$,
all small-item demands are zero. Finally, since all fractions are
nonnegative and sum to $y_e$, each one is at most $y_e$.
\end{proof}

The fractions satisfy the required marginals and capacities, but need not form a~distribution on~sets of~distinct items within each profile.
The rounding step selects at most one option for each edge.

\subsection{Choosing Compatible Profiles and Whole Edges}\label{sec:round-proof}
We~first choose one profile for every vertex and color.
These choices must preserve most of~the proposed weight while leaving slack in~every slot and reserved space.
We~then select disjoint slot pairs for the large edges and round the small-edge options by~linear algebra.
The slack created in~the first step pays for the losses in~the second.

\restatedstatement{\theorem}{\localglobalrounding}

\begin{proof}[Proof of~\Cref{thm:round}]
All constants hidden in $O(\cdot)$ below may depend on the fixed parameters,
but not on $n$ or on the number of vertices.

\paragraph*{Dividing the vertices into groups}
Choose a fixed integer $r\ge\max\{2,4/\eps\}$. Assign each vertex
independently and uniformly to one of $r$ groups. For a vertex $v$ and
group $j$, let
\[
 L_{v,j}=\sum_{u:uv\in E(G)}w_{uv}\one_{\{u\text{ is in group }j\}}.
\]
This is the weight from $v$ to group $j$, and its expectation is
$W_G(v)/r$. Apply \Cref{prop:HS} to all these sums. One
choice at a source vertex $u$ contributes at most one to each sum,
because $G$ is simple, and its contributions over all sums add up to
$W_G(u)\le n$. Thus the proposition has $R\le n/r$ and $D\le n$.
It gives a division into groups $V_1,\ldots,V_r$ for which
\[
 L_{v,j}\le \frac{W_G(v)}r+O(\sqrt{n\log n})
 \quad\text{for every }v,j.
\]
This bound allows arbitrarily many edges of very small weight.

Let $H$ contain just the edges joining different groups. Set aside the
other edges for now. Their total weight at every vertex is at most
$n/r+O(\sqrt{n\log n})$. Restrict the proposed fractions to $H$.
They remain averages of feasible sets: intersect each original feasible
set with $E(H)$. If
\[
 T_v^G=\sum_{e\in E_G(v)}w_e\sum_cx_{e,c},\qquad
 T_v^H=\sum_{e\in E_H(v)}w_e\sum_cx_{e,c},
\]
then $\sum_cx_{e,c}\le1$ gives
\begin{equation}\label{eq:deleted-weight}
 T_v^G-T_v^H\le \frac nr+O(\sqrt{n\log n}).
\end{equation}
Until the final error calculation, work only with $H$: every edge,
incident-edge set, and sum below belongs to $H$.

\paragraph*{Profiles and lost weight}
Choose small constants $0<\delta,\theta,\sigma<1$, put $r_0=\delta/4$,
and choose a sufficiently small positive constant $\beta\le r_0$.
The precise requirements on these constants will be given at the end.
Call an edge large if $w_e\ge\beta$, and small otherwise.

Apply \Cref{lem:profiles} to every vertex $v$ and color $c$.
Write $\lambda_{v,c,P}$ for the probability of profile $P$, and
$a_{v,c}(e;P,i)$ for the fraction assigned to edge $e$ in its slot $i$.
For a small edge, $i$ instead denotes the space for small edges.
In either case these fractions sum, over all descriptions and places, to
\[
 y_{e,c}=(1-\delta)x_{e,c}.
\]
Use the same number $M$ of profile labels everywhere, adding empty
descriptions of probability zero if necessary. The bound on $M$ depends
only on $\beta,\delta$.

It is useful to give every label a positive probability. Set
\[
 \lambda'_{v,c,P}=(1-\theta)\lambda_{v,c,P}+\theta/M,
 \qquad \rho=(1-\theta)(1-\sigma).
\]
Initially, choose all descriptions $P_{v,c}$ independently according to
$\lambda'$. For an edge $e=uv$, a color $c$, and a place $i$ at $u$ and
a place $j$ at $v$, define
\[
 X(e,c,i,j)=
 \frac{\rho\,a_{u,c}(e;P_{u,c},i)a_{v,c}(e;P_{v,c},j)}
 {y_{e,c}\lambda'_{u,c,P_{u,c}}\lambda'_{v,c,P_{v,c}}}.
\]
If $y_{e,c}=0$, the value is zero. A positive value is a fractional option
to give $e$ color $c$ and use the stated places at its two ends. Large
edges use ordinary slots at both ends; small edges use the space for small
edges at both ends. Put
\[
 S_e=\sum_{c,i,j}X(e,c,i,j),\qquad
 \mu_e=\rho\sum_c y_{e,c}\le1.
\]
The factor $\lambda'$ in each denominator cancels the probability of
choosing that description. Summing the assigned fractions at each end
therefore gives $\E S_e=\mu_e$. More explicitly, for each color the two endpoint sums give $y_{e,c}^2$; division by $y_{e,c}$ leaves $y_{e,c}$.

For one color, the sum of the $X$ values is $O(1/n)$. Indeed,
$\sum_i a_{v,c}(e;P,i)\le y_{e,c}$,
$\lambda'_{v,c,P}\ge\theta/M$, and $y_{e,c}\le C_0/n$.
There are $O(n)$ colors, so $S_e=O(1)$ for every possible set of choices.
The terms from distinct colors are independent. Adding their variances
gives $\operatorname{Var}(S_e)=O(1/n)$.

At a vertex, measure the total error by adding the absolute errors of its
incident edges, multiplied by their weights:
\[
 Z_v=\sum_{e\ni v}w_e|S_e-\mu_e|.
\]
The Cauchy--Schwarz inequality gives
$\E|S_e-\mu_e|\le\sqrt{\operatorname{Var}(S_e)}$.
Since $\sum_{e\ni v}w_e\le n$, we obtain
\begin{equation}\label{eq:L1}
 \operatorname{Var}(S_e)=O(1/n),\qquad \E Z_v=O(\sqrt n).
\end{equation}

We must also respect the slots and the space for small edges. For an
ordinary slot, let $d(e)=1$ and let its capacity be $h(P,i)=1$; it can
hold at most one edge. For the space for small edges, let $d(e)=w_e$ and
$h(P,i)=R(P)$; its capacity measures total weight. Its load is the sum
of $d(e)X(e,c,i,j)$ over all options using it.
Fix a profile $P$ at $(v,c)$ and one of its places $i$. Averaging over
the choices at the other endpoints gives load at most
\begin{equation}\label{eq:meanload}
 \frac{\rho}{\lambda'_{v,c,P}}
       \sum_{e\ni v}d(e)a_{v,c}(e;P,i)
 \le (1-\sigma)h(P,i).
\end{equation}
Here we used the capacity bounds of \Cref{lem:profiles} and
$\lambda'\ge(1-\theta)\lambda$. We omit spaces of capacity zero.
Every remaining capacity is at least $r_0$.

\paragraph*{Choosing profiles group by~group}
The groups $V_1,\ldots,V_r$ have no edges within them in $H$.
We choose the lists $P_v=(P_{v,c}:c\in\cC)$ one group at a time,
as shown in \Cref{fig:group-passes}. All references to the original distribution of a profile list mean
the independent description choices with probabilities $\lambda'$.

\begin{figure}[htbp]
\centering
\begin{tikzpicture}[x=1cm,y=1cm,font=\small]
\foreach \xx/\name in {0/1,3.5/2,7/3}{
 \draw[fill=softblue,draw=inkblue,rounded corners] (\xx-.95,-.4) rectangle (\xx+.95,1.6);
 \node at (\xx,1.25) {$V_{\name}$};
 \node[vertex] (a\name) at (\xx-.45,.55) {};
 \node[vertex] (b\name) at (\xx+.45,.05) {};
 \draw[black!45,thick,dashed] (a\name)--(b\name);
}
\draw[inkblue,thick] (a1)--(a2);
\draw[inkblue,thick] (b1)--(b2);
\draw[inkblue,thick] (a2)--(b3);
\draw[inkblue,thick] (b2)--(a3);
\draw[inkblue,thick] (a1) to[bend left=24] (a3);
\node at (9,.6) {$\cdots\quad V_r$};
\node[align=center] at (3.5,2.6) {Keep edges between groups; set aside edges within groups.};
\node[panel,fill=softgreen,minimum width=1.9cm] (s1) at (0,-1.45) {choose $V_1$};
\node[panel,fill=softgreen,minimum width=1.9cm] (s2) at (3.5,-1.45) {choose $V_2$};
\node[panel,fill=softgreen,minimum width=1.9cm] (s3) at (7,-1.45) {choose $V_3$};
\draw[-{Latex},thick] (s1)--(s2);
\draw[-{Latex},thick] (s2)--(s3);
\node at (9,-1.45) {$\cdots\quad V_r$};
\end{tikzpicture}
\caption{Divide the vertices into a fixed number of groups. The dashed
edges have small total weight at each vertex and will be colored later.
Choose descriptions one group at a time, keeping all earlier choices.
There are no edges within a group in the remaining graph, so choices
made in one step can be treated as independent.}
\label{fig:group-passes}
\end{figure}
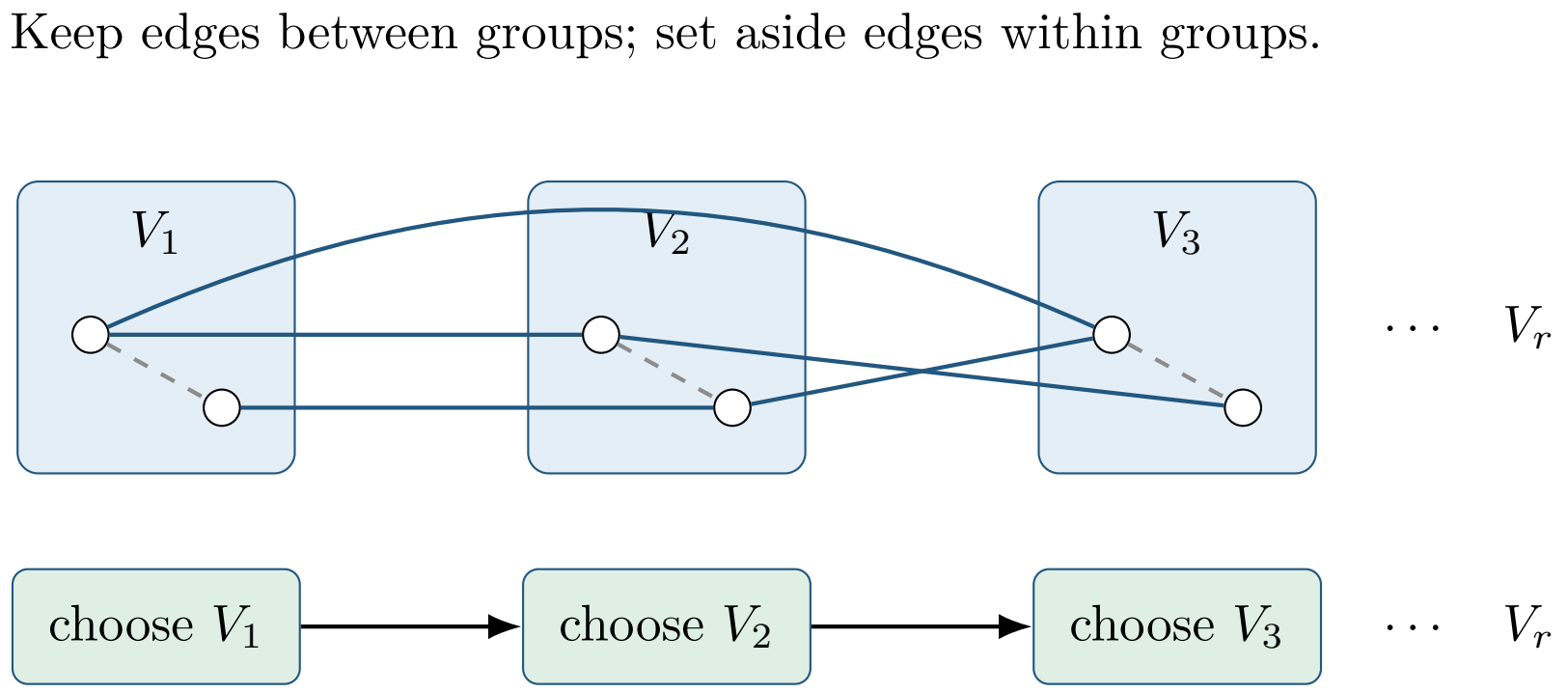

We use second moments to compute an upper bound on the remaining
absolute error. Choose a power of two $\kappa$ with $(2n^2)^{-1}<\kappa\le n^{-2}$,
and put
\begin{equation}\label{eq:rational-square-root}
 \phi(z)=\kappa\left\lceil\frac{\sqrt z}{\kappa}\right\rceil
 \qquad(z\ge0).
\end{equation}
Thus $\sqrt z\le\phi(z)<\sqrt z+\kappa$. After groups
$V_1,\ldots,V_k$ have been chosen, let $\mathcal F_k$ record their lists
and define
\begin{equation}\label{eq:computable-progress}
 U_e^{(k)}=\E[(S_e-\mu_e)^2\mid\mathcal F_k],\qquad
 Q_v^{(k)}=\sum_{e\ni v}w_e\phi(U_e^{(k)}).
\end{equation}
Every unchosen list in this expectation has its original distribution;
this is a way to measure progress, not a rule for the later choices.
The moments are computed by adding the first and second moments of the
independent color contributions. \Cref{sec:computable-progress}
gives the formulas and the rational implementation of $\phi$.
Initially $Q_v^{(0)}=O(\sqrt n)$, by~\eqref{eq:L1} and
$\sum_{e\ni v}w_e\kappa\le n\kappa\le1/n$.

Under the original distributions for the lists in the next group,
the tower property and the concavity of the square root give
\[
 \E[\phi(U_e^{(k)})\mid\mathcal F_{k-1}]
 \le\sqrt{U_e^{(k-1)}}+\kappa
 \le\phi(U_e^{(k-1)})+\kappa.
\]
Consequently
\begin{equation}\label{eq:progress-one-step}
 \E[Q_v^{(k)}\mid\mathcal F_{k-1}]
 \le Q_v^{(k-1)}+1/n.
\end{equation}
These inequalities replace the need to compute conditional expectations
of absolute values.

Put $t=\lceil n^{3/4}\rceil$. At step $k$, for each $u\in V_k$, retain only lists
for which $Q_u^{(k)}\le t$. Every edge incident to $u$ has its other
endpoint outside $V_k$: that endpoint's list is either already fixed
or averaged over its original distribution. Thus $Q_u^{(k)}$ depends on
$\mathcal F_k$ only through the list at $u$ and the earlier choices
$\mathcal F_{k-1}$, not on any other list in $V_k$.
We choose independently from the retained
lists, with their original probabilities rescaled to sum to one.

We show inductively that, just before step $k$, every unchosen vertex
satisfies
\[
 Q_v^{(k-1)}=O_r(\sqrt{n\log n}).
\]
Thus Markov's inequality, $\Prb(Y>t)\le\E Y/t$ for $Y\ge0$,
bounds the probability removed at each vertex of $V_k$ by the computable
number
\[
 p_k=\max_{u\in V_k}\frac{Q_u^{(k-1)}+1/n}{t}
     =O_r(n^{-1/4}\sqrt{\log n})=o(1).
\]
For an empty group take $p_k=0$. The constant can be chosen to cover all
$r$ steps. In particular, every retained distribution has positive total
probability for sufficiently large $n$. To sample it, repeatedly draw an
original list until its computed value of $Q_u^{(k)}$ is at most $t$.
The expected number of draws is at most $1/(1-p_k)$.

Apply \Cref{prop:HS} to the independent retained lists in
$V_k$, imposing the following bounds.
\begin{itemize}
\item For every vertex $v$, color $c$, possible profile $P$ at
$(v,c)$, and place $i$ in $P$, bound the load contributed by neighbors
of $v$ in $V_k$, on the assumption that $v$ uses $P$.
\item For every $v\notin V_k$, bound the part of $Q_v^{(k)}$ contributed
by neighbors in $V_k$.
\end{itemize}
The first kind of row is zero when $v\in V_k$. Every nonzero row is a
sum of nonnegative functions, one function of each neighboring list
being chosen. For an error row, if $v$ was chosen earlier, its list is
fixed; if it will be chosen later, average over its original list.
Thus the contribution of a neighbor $u\in V_k$ is respectively
\[
 w_{uv}\phi\bigl((S_{uv}(P_u,P_v)-\mu_{uv})^2\bigr)
 \quad\text{or}\quad
 w_{uv}\phi\bigl(\E_{P_v}[(S_{uv}(P_u,P_v)-\mu_{uv})^2]\bigr).
\]
Each value is computable from first and second moments of the color
contributions. The contributions of all other neighbors are unchanged
during this step.

We give the size estimates needed for the proposition. In the following
table, a ``source'' is one vertex whose list we are choosing. A row is one
of the quantities to be bounded. The last column adds the contributions
of one fixed source choice over all rows. The scaling makes every
individual contribution at most one. Here $B_0$ is a sufficiently large
fixed constant.
\begin{equation}\label{eq:APtable}
\begin{array}{c|c|c|c}
 \text{quantity bounded}&\text{one source's contribution}
 &\text{scaling}&\text{sum over all rows}\\ \hline
 \text{load}&O(1/n)&\text{multiply by }n/B_0&O(n^2)\\
 \text{error bound}&O(w_{uv})&\text{divide by }B_0&O(n)
\end{array}
\end{equation}
All expected row sums after scaling are $O(n)$.
These estimates hold at every step.

We~verify that the bounds in~\eqref{eq:APtable} do~not depend on~the number of~vertices.
Since $H$ is~simple, a~source has at~most one edge to~each target.
Its contribution to one ordinary-slot load is therefore $O(1/n)$; its
contribution to a small-edge load is $O(w_e/n)$. For a fixed source list
and color, add the contributions over all possible profiles and places
at all neighboring endpoints. There are only boundedly many description
labels and places at an endpoint. Summing the factors for the other
endpoint in the formula for $X$ cancels $y_{e,c}$, up to a fixed factor,
because those assigned fractions sum to $y_{e,c}$ and $\lambda'\ge\theta/M$.
The remaining sum at the source is bounded by its total assigned
ordinary-slot fraction plus its total assigned small-edge weight, divided
by its profile probability. \Cref{lem:profiles} bounds these by
$O(1/\beta)$ and $O(1)$, respectively, since
$\lambda/\lambda'\le1/(1-\theta)$. Thus the sum is $O(1)$ per color,
or $O(n)$ over all colors before scaling. This gives $O(n^2)$
in the first row of the table. For the second row,
$|S_e-\mu_e|=O(1)$ implies $\phi(U_e^{(k)})=O(1)$ for every
choice of fixed lists. Thus one edge contributes $O(w_e)$, and the sum
of the weights at a source is at most $n$.

To check the expected row sums, first use the original list distributions.
A load row has mean at most the full load in~\eqref{eq:meanload}, hence
at most one. Restricting a list to an event of probability at least
$1-p_k$ increases the expectation of each nonnegative contribution by
at most $1/(1-p_k)$. An error row is $O(n)$ for every possible set of
lists, since $\phi(U_e^{(k)})=O(1)$ and the incident weights sum to at
most $n$.
After scaling as in the table, \Cref{prop:HS} therefore has
$R=O(n)$ and $D=O(n^2)$ at every step.

For an error row, let $C$ be the total contribution from neighbors
outside $V_k$. This is a fixed nonnegative quantity at this step.
The original expected contribution from neighbors in $V_k$ is at most
$Q_v^{(k-1)}-C+1/n$, by the square-root inequality above.
Use this computable upper bound, divided by $1-p_k$, for the retained
lists. \Cref{prop:HS} therefore gives
\[
 Q_v^{(k)}
 \le C+\frac{Q_v^{(k-1)}-C+1/n}{1-p_k}
       +O(\sqrt{n\log n})
 \le\frac{Q_v^{(k-1)}}{1-p_k}+O(\sqrt{n\log n}).
\]
Only upper bounds on row means are needed; the sampling procedure does
not evaluate expectations under the retained distributions.
\Cref{sec:computable-progress} explains why the
Harris--Srinivasan algorithm can use these implicitly described lists.
There are only $r$ steps. For a vertex not yet chosen, induction from
$Q_v^{(0)}=O(\sqrt n)$ proves the claimed
$O_r(\sqrt{n\log n})$ bound, justifying the removed-probability bound
used above. More explicitly, at each induction step we may first use
$1/(1-p_k)\le2$. The constants then grow by at most a fixed factor
per step. Since there are only $r$ steps, choosing $n$ sufficiently
large makes every removed probability less than $1/2$, as required.
At the step when $v$ itself is chosen, the retained-list
condition gives $Q_v^{(k)}\le t$. At each later step the displayed
inequality applies again. Consequently, after the last step,
\[
 Z_v\le Q_v^{(r)}=O_r(n^{3/4})\qquad\text{for every }v.
\]

To check the load bounds, fix a vertex $v$, a color $c$, a possible
description $P$, and one of its places $i$. Under the original lists,
the expected contribution from neighbors in $V_k$ is exactly
\[
 m_k=\frac{\rho}{\lambda'_{v,c,P}}
       \sum_{\substack{e=uv\ni v\\u\in V_k}}
                d(e)a_{v,c}(e;P,i).
\]
In particular, $\sum_{k=1}^r m_k\le(1-\sigma)h(P,i)$ by
\eqref{eq:meanload}. \Cref{prop:HS}, followed by undoing
the load scaling, bounds the actual contribution from group $V_k$ by
\[
 \frac{m_k}{1-p_k}
 +O\left(\sqrt{\frac{\log n}{n}}+\frac{\log n}{n}\right).
\]
We imposed this bound for every possible $P$, so it holds for the
description eventually chosen at $v$, whether that choice was made
before, during, or after step $k$. Put $p=\max_{1\le k\le r}p_k=o(1)$.
Adding the $r$ bounds gives total load at most
\[
 \frac{(1-\sigma)h(P,i)}{1-p}
 +O_r\left(\sqrt{\frac{\log n}{n}}+\frac{\log n}{n}\right)
 \le(1-\sigma/2)h(P,i)
\]
for sufficiently large $n$, since every positive capacity is at least
$r_0$. We~have therefore bounded the lost weight at~every vertex and retained slack in~every capacity.
It~remains to~choose whole edges from the fractional options.

An edge may currently have total assigned fraction $S_e>1$.
Divide all its option fractions by $\max\{1,S_e\}$, and call the new
fractions $\widehat X(e,c,i,j)$. Their sum and the unused fraction are
\[
 s_e=\min\{1,S_e\},\qquad r_e=1-s_e.
\]
This only decreases loads. Since $\mu_e\le1$, we have
$s_e\ge\mu_e-|S_e-\mu_e|$. Hence
\begin{equation}\label{eq:capture}
 \sum_{e\ni v}w_es_e
 \ge\sum_{e\ni v}w_e\mu_e-Z_v
 \ge\rho(1-\delta)\sum_{e\ni v}w_e\sum_cx_{e,c}
       -O(n^{3/4}).
\end{equation}

\paragraph*{Rounding the large edges}
For large edges, each option uses one ordinary slot at each end. We use
\Cref{prop:DP} to choose options without using a slot twice.
We also provide places for edges left uncolored, so that the number of
such edges can be controlled at every vertex.

Divide $[\beta,1]$ into finitely many weight intervals, assigning an edge
on a shared boundary to just one interval. Let $\ell_h,u_h$ be the lower
and upper endpoints of interval $h$, with
$u_h\le(1+\xi)\ell_h$, where $\xi>0$ is a small fixed constant.
Let $E(v,h)$ be the edges at $v$ in that interval, and put
\[
 R_{v,h}=\sum_{e\in E(v,h)}r_e,
 \qquad
 B_{v,h}=\left\lceil\frac{R_{v,h}}{1-\gamma}\right\rceil
                   +\lceil n^{2/3}\rceil,
\]
where $0<\gamma<1$ is another small fixed constant.
Create $B_{v,h}$ extra slots at $(v,h)$, reserved for uncolored edges.
Each such slot may receive at most one edge. These are counting devices;
they are separate from the slots belonging to colors.

Construct a family of three-element sets, allowing repeated copies of a
set. Introduce one new element for each large edge; call it the job
element for that edge. Let $A$ be all job elements. The remaining elements,
forming a set $B$, are the ordinary slots in the selected profiles
and the extra slots just created. Each three-element set will consist of
one job element and one slot at each endpoint of that edge.

Put $N=\lceil n\rceil^6$. For each color option $(e,c,i,j)$, add
\[
 \lfloor N\widehat X(e,c,i,j)\rfloor
\]
copies of the set containing the job element of $e$ and its two ordinary
slots. For an edge $e=uv$ in interval $h$, divide the unused fraction
$r_e$ equally among all pairs of extra slots, one at $(u,h)$ and one at
$(v,h)$. For each such pair add
\[
 \left\lfloor\frac{Nr_e}{B_{u,h}B_{v,h}}\right\rfloor
\]
copies of the set containing that pair and the job element of $e$.
Selecting a set of the first kind will color the edge; selecting a set
of the second kind will leave it uncolored.

\begin{figure}[htbp]
\centering
\begin{tikzpicture}[x=1cm,y=1cm,font=\small]
\node[panel,minimum width=2.8cm] (uc) at (-3.5,1) {slot at $u$\\color $c$};
\node[panel,minimum width=2.8cm] (vc) at (3.5,1) {slot at $v$\\color $c$};
\node[panel,fill=black!4,draw=black!50,minimum width=2.8cm] (ud) at (-3.5,-1) {extra slot at $u$\\leave uncolored};
\node[panel,fill=black!4,draw=black!50,minimum width=2.8cm] (vd) at (3.5,-1) {extra slot at $v$\\leave uncolored};
\node[draw,rounded corners,fill=white,minimum width=1.3cm,minimum height=.7cm] (e) at (0,0) {edge $e$};
\draw[inkblue,very thick] (uc)--(e)--(vc);
\draw[black!55,thick,dashed] (ud)--(e)--(vd);
\node[align=center] at (0,1.9) {two possible choices for the same edge};
\end{tikzpicture}
\caption{Each option is a set containing an edge label and one slot at each
endpoint. The upper option gives color $c$; the lower option leaves the
edge uncolored and uses an extra slot at each end. The two options share
the edge label, so at most one can be selected. Disjoint selected sets
also prevent a slot from being used by two edges.}
\label{fig:matching-options}
\end{figure}
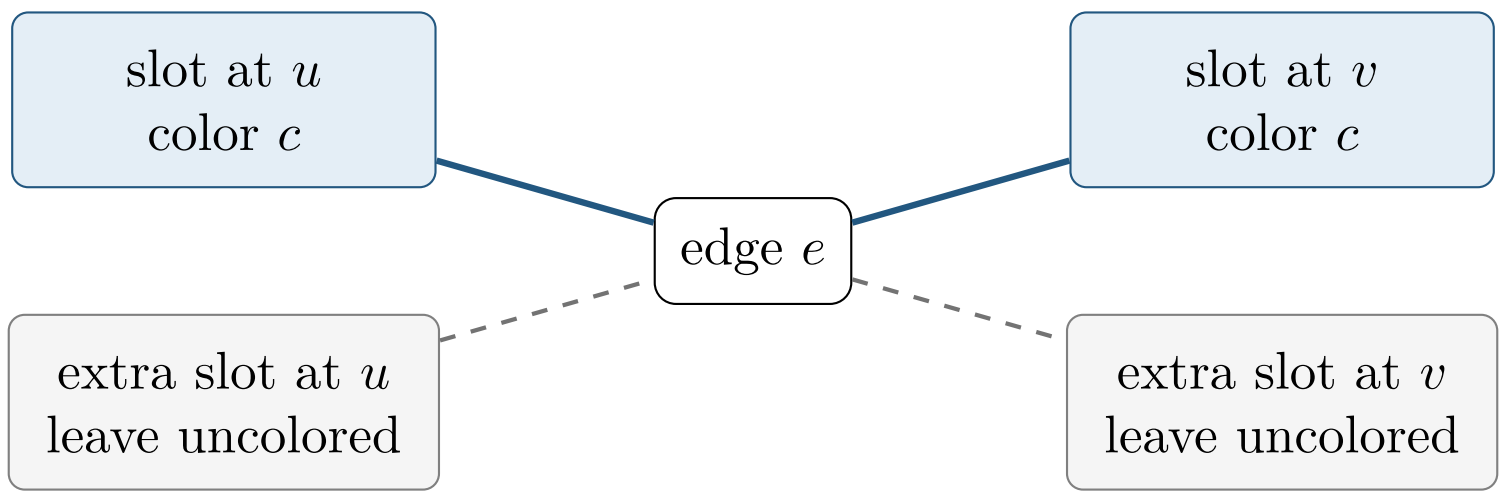

For an element $z$, let $d(z)$ be the number of sets containing it,
counting copies. Let $\Delta_2$ be the largest number of sets containing
any given pair of distinct elements, again counting copies. We claim
\begin{equation}\label{eq:degrees}
 \begin{aligned}
 d(a)&=N-O(n^2) &&(a\in A),\\
 d(b)&\le(1-\eta)N &&(b\in B),\\
 \Delta_2&=O(Nn^{-2/3}),
 \end{aligned}
 \qquad \eta=\min\{\sigma/2,\gamma\}>0.
\end{equation}
To check the first line, a job has $O(n)$ color options, since there
are $O(n)$ colors and boundedly many slots per color. Also,
$R_{v,h}\le n/\beta$, so $B_{v,h}=O(n)$ and the job has $O(n^2)$
pairs of extra slots. Before taking integer parts, the numbers of copies
for that job sum to $N(s_e+r_e)=N$. Each integer part loses less than
one, proving the first line.

For an ordinary slot, the number of sets using it is at most its
fractional load multiplied by $N$, which is at most $(1-\sigma/2)N$.
For an extra slot at $(v,h)$, the corresponding bound is
$NR_{v,h}/B_{v,h}\le(1-\gamma)N$. This proves the second line.

For the last line, consider the possible pairs. A job element and an
ordinary slot lie together in $O(N/n)$ sets: the slot fixes the color,
and the total fraction for one edge and color is $O(1/n)$.
A job element and an extra slot lie together in at most
$N/B_{v,h}=O(Nn^{-2/3})$ sets. Two ordinary slots at opposite endpoints
can occur together only for the unique original edge joining those
endpoints, and again the bound is $O(N/n)$. Two extra slots at opposite
endpoints occur together in at most $N/(B_{u,h}B_{v,h})$ sets, for the
same reason. All other pairs of slot elements occur together in no set,
as do pairs of job elements. This proves \eqref{eq:degrees}.

Take $D=\lfloor(1-\eta/2)N\rfloor$ and any fixed $0<\zeta<1/9$.
For large $n$, the hypotheses of \Cref{prop:DP} hold.
Indeed, every job element belongs to more than $D$ sets by a fixed
positive fraction of $N$, whereas every slot element belongs to fewer
than $D$ sets. The extra factor $D^{-a}$ required by that proposition
tends to zero. Finally, since $N$ has order $n^6$,
\[
 Nn^{-2/3}=o(D^{1-\zeta}).
\]
The proposition supplies pairwise disjoint sets covering all job
elements. Thus it chooses exactly one option for each large edge and
uses each ordinary or extra slot at most once.

At most $B_{v,h}$ edges of $E(v,h)$ are left uncolored, because each
uses a different extra slot. Their total weight is therefore at most
\begin{equation}\label{eq:residual}
 u_hB_{v,h}
 \le\frac{1+\xi}{1-\gamma}
                \sum_{e\in E(v,h)}w_er_e+O(n^{2/3}).
\end{equation}
The number of intervals depends only on the fixed constants, so summing
this bound over all intervals still gives an $O(n^{2/3})$ error.

\paragraph*{Rounding the small edges}
For a small edge $e$, use its color options of fractions $\widehat X$
and one additional option of fraction $r_e$ that leaves it uncolored.
These fractions sum to one. For the space for small edges of size $R>0$,
divide every load by $R$, so the capacity becomes one.
A color option then adds $w_e/R\le\beta/r_0$ to each of two capacity
constraints, one at each endpoint. The current fractional load on each
such space is at most $1-\sigma/2$.

At each vertex add one more constraint that records the total weight of
small edges left uncolored. The uncolored option for $e=uv$ adds $w_e$
to this sum at $u$ and at $v$. Thus every option contributes to at most
two constraints, and every contribution is at most
$\tau=\beta/r_0$; here $r_0<1$.
Choose $\beta$ small enough that
\[
 2\beta/r_0\le\sigma/2.
\]
\Cref{lem:two} now chooses one whole option for every small edge.
It increases each load by at most $2\beta/r_0$.
Consequently all the spaces for small edges obey their capacities, and
the small uncolored weight at $v$ is at most
\[
 \sum_{\substack{e\ni v\\w_e<\beta}}w_er_e+2\beta/r_0.
\]

For each vertex and color, a large edge occupies a distinct ordinary
slot and has weight at most that slot's size. The total weight of the
small edges fits the space set aside for them. A profile has
total size at most one, so the resulting partial coloring is feasible.
Every color option used had positive assigned fraction, and therefore
had $x_{e,c}>0$.

It remains to add the errors. Write
$T_v=T_v^H=\sum_{e\ni v}w_e\sum_cx_{e,c}$ and
$\Lambda=(1+\xi)/(1-\gamma)$. We have $T_v\le W_H(v)\le n$.
The bounds for large and small uncolored edges give total uncolored
weight at most
\[
 \Lambda\sum_{e\ni v}w_er_e+O(n^{2/3}).
\]
Subtract this from $W_H(v)$, use $r_e=1-s_e$, and then apply
\eqref{eq:capture}. The colored weight is at least
\[
 \begin{split}
 W_H(v)-\Lambda\sum_{e\ni v}w_er_e-O(n^{2/3})
 &\ge \sum_{e\ni v}w_es_e-(\Lambda-1)n-O(n^{2/3})\\
 &\ge T_v-\bigl[1-\rho(1-\delta)+(\Lambda-1)\bigr]n
                 -O(n^{3/4}).
 \end{split}
\]
Choose $\delta,\theta,\sigma,\xi,\gamma$ positive and small enough that
the quantity in square brackets is less than $\eps/2$.
Then set $r_0=\delta/4$ and choose $\beta$ as required above.
By~\eqref{eq:deleted-weight}, the colored weight in the original graph
is therefore at least
\[
 T_v^G-\left(\frac1r+\frac\eps2\right)n-O(n^{3/4})
 \ge T_v^G-\eps n-O(n^{3/4}).
\]
Here $1/r\le\eps/4$, and $\sqrt{n\log n}=O(n^{3/4})$.
The integer $r$ and all other choices depend only on $C_0,C_1,\eps$.
Taking $n$ sufficiently large proves the theorem, with no dependence
on the number of vertices.
\end{proof}

\begin{algorithm}[Rounding the proposed colors]\label{alg:round}
	Input: the graph, palette, fractions, and local feasible-bin distributions in~\Cref{thm:round}; fix the desired error.
	\begin{enumerate}[leftmargin=*,nosep]
		\item Partition the vertices into a fixed number of groups. Set aside
		the small total weight of~edges lying within groups.
		\item Build the profile menus and the paired endpoint fractions.
		\item Choose profile lists one group at a time, passing the load and
		error tests in~\Cref{sec:round-proof}.
		\item Clip an edge's total proposal to one if necessary. Create ordinary
		and extra slots, then find disjoint options for the large edges.
		\item Round one option for each small edge by \Cref{lem:two}.
		Return the resulting feasible partial coloring.
	\end{enumerate}
	Procedures that run in~expected polynomial time for each choice are specified in~\Cref{sec:rounding-implementation}.
\end{algorithm}

\subsection{The Bound for Simple Graphs}\label{sec:simple}
We~apply \Cref{thm:round} to~prove the simple-graph bound in~\Cref{thm:main}.
The empty graph is~immediate, so~assume that $G$ has an~edge and put $n=\bp(G)$.
At~each vertex $v$, choose a~packing into $n$ bins, allowing empty bins, and denote the edge sets in~these bins by~$B_{v,1},\ldots,B_{v,n}$.

Use the colors $1,\ldots,n$.
For every edge and color set
\[
x_{e,c}=1/n.
\]
For a~fixed $v$ and $c$, the vector of~these values is
\[
(x_{e,c}:e\ni v)=\frac1n\sum_{j=1}^n\one_{B_{v,j}}.
\]
It~is~an~average of~feasible sets, so~it~belongs to~$\mathcal K_v$.
Also $\sum_cx_{e,c}=1$ for every edge, and $W_G(v)\le n$.
Thus \Cref{thm:round} applies with $C_0=C_1=1$.

Fix $\eps>0$.
The theorem gives a~partial coloring with the $n$ colors in~which the weight left uncolored at~any vertex is~at~most
\[
R=\eps n+O_\eps(n^{3/4}).
\]
\Cref{lem:general-cleanup} colors those remaining edges with at~most $4R+1$ additional colors.
Therefore
\[
\chi'_w(G)\le n+4\eps n+O_\eps(n^{3/4}).
\]
After division by~$n$, the last error term tends to~zero.
Since $\eps$ may be~chosen arbitrarily small, this proves $\chi'_w(G)\le(1+o(1))n$.
Every threshold used depends only on~$\eps$, not on~the number of~vertices.
For the finitely many smaller values of~$n$, \Cref{lem:general-cleanup} gives the bound $4n+1$, so~they can be~included by~increasing the error function at~those values.

\noindent\textbf{Remark.}
	A~triangle with unit-weight edges has $\bp(G)=2$ and needs three colors.
	Thus, the asymptotic bound for simple graphs does not assert equality for every value of~$\bp(G)$.
	The proof above applies to~simple graphs with odd cycles as~well as~to~bipartite graphs.

\begin{algorithm}[Coloring a simple graph]\label{alg:simple}
Input: a weighted simple graph and a desired accuracy.
Choose the local-packing and rounding accuracies sufficiently small
for the desired final error.
\begin{enumerate}[leftmargin=*,nosep]
\item Compute the approximate local packings in
\Cref{sec:local-input-algorithms}. Let $n$ be the largest number of bins
used, and pad the other packings with empty bins.
\item Set $x_{e,c}=1/n$ for $c\in[n]$, using a uniform choice among the
local bins as each feasible distribution.
\item Run \Cref{alg:round} with these fractions and distributions.
\item Color the remaining edges in decreasing weight order with fresh
colors, as in \Cref{lem:general-cleanup}. Return the completed coloring.
\end{enumerate}
The extra factor from the approximate local packings is absorbed in
the chosen error. This randomized algorithm has expected polynomial
running time for each fixed accuracy.
\end{algorithm}

\subsection{Computing Progress and Sampling Vertex Lists}\label{sec:rounding-implementation}
\label{sec:computable-progress}

A~vertex has exponentially many possible profile lists, so~the algorithm must evaluate progress and sample retained lists without enumerating them.
The conditional second moments in~\eqref{eq:computable-progress} provide the required computations.
All probabilities, weights, and assigned fractions are rational.
There are $O(n)$ colors, and the number $M$ of~profile labels depends only on~the fixed accuracy parameters.

The algorithm needs three operations: evaluate the progress value for a
proposed list, draw a list that passes its local test, and resample a
small set of lists when a load constraint fails. We describe them in
that order. An entire list is one variable in the resampling procedure,
but drawing or testing that variable only requires going through its
colors. We never store the set of all possible lists.

Fix an edge $e=uv$ and a color $c$. For a profile $P$, define
\[
 A_{v,c}(P)=\sum_i a_{v,c}(e;P,i),\qquad
 b_{v,c}(P)=\frac{A_{v,c}(P)}{\lambda'_{v,c,P}}.
\]
The sum includes all ordinary slots, or the reserved space when $e$ is
small. If the list at $v$ has already been fixed, put
\[
 m_{v,c}=b_{v,c}(P_{v,c}),\qquad
 q_{v,c}=b_{v,c}(P_{v,c})^2.
\]
If that list is still unchosen, its original distribution instead gives
\[
 m_{v,c}=\sum_P A_{v,c}(P)=y_{e,c},\qquad
 q_{v,c}=\sum_P\frac{A_{v,c}(P)^2}{\lambda'_{v,c,P}}.
\]
Define $m_{u,c},q_{u,c}$ in the same way. The contribution of this color
to $S_e$ is
$\rho b_{u,c}(P_{u,c})b_{v,c}(P_{v,c})/y_{e,c}$.
Thus its first and second moments are
\[
 a_c=\frac{\rho m_{u,c}m_{v,c}}{y_{e,c}},\qquad
 d_c=\frac{\rho^2q_{u,c}q_{v,c}}{y_{e,c}^2}.
\]
When $y_{e,c}=0$, set $a_c=d_c=0$. The unchosen endpoint lists are
independent, and their choices at different colors are independent.
Fixing any complete endpoint lists preserves the latter independence.
It follows that
\begin{equation}\label{eq:progress-moment-formula}
 U_e^{(k)}
 =\left(\sum_c a_c-\mu_e\right)^2
   +\sum_c(d_c-a_c^2).
\end{equation}
This formula evaluates the conditional second moment using only
$O(M|\mathcal C|)$ operations per edge. In particular, it does not sum
over all possible complete lists.
The two terms have different roles. The first is the square of the
difference between the current mean and the target $\mu_e$; the second
is the variance still contributed by unchosen profiles. After the final
group no lists remain unchosen, so this variance term is zero and only
the actual squared error remains.
Here the expectations for unchosen lists are always taken under their
original distributions. The later rejection step is handled by an upper
bound on conditional means, rather than by asserting independence
between colors after rejection.

Use $\kappa$ and $\phi$ from \eqref{eq:rational-square-root}.
The number of bits needed to write $\kappa$ is $O(\log(n+2))$.
For a nonnegative rational $z$, compute $\phi(z)$ by finding the
least nonnegative integer $b$ with
$b^2\kappa^2\ge z$. Integer comparisons and binary search compute it in
polynomial time. The same method computes the integer threshold
$t=\lceil n^{3/4}\rceil$, using $t^4\ge n^3$. Every progress value and
every test that accepts or rejects a list is therefore rational and
polynomially computable.

By \eqref{eq:progress-one-step}, this upward rounding adds at most
$1/n$ to the expected potential at each step, under the original
list distributions. At the end, $U_e^{(r)}=(S_e-\mu_e)^2$, so
$Z_v\le Q_v^{(r)}$. Thus the original $O(n^{3/4})$ absolute-error bound
is retained.

At group $V_k$, first compute
\[
 p_k=\max_{u\in V_k}\frac{Q_u^{(k-1)}+1/n}{t},
\]
with value zero for an empty group. To draw a retained list at $u$, draw
its $O(n)$ profile labels independently from the original rational
probabilities and test $Q_u^{(k)}\le t$. Repeat on rejection. Markov's
inequality bounds the rejection probability by $p_k$. In the range where
the proof is applied, $p_k=o(1)$ and in particular $p_k\le1/2$, so the
expected number of attempts is at most two. There is no edge within
$V_k$, so the acceptance test for $u$ depends only on its own list and
on earlier groups. It remains the same throughout this group's
resampling procedure. Consequently this procedure samples exactly the
required independent retained distributions, even though their
normalizing probabilities are not computed.
The choices in different colors of a retained list may now be dependent.
Rejection sampling automatically gives each
accepted list its original probability divided by the total acceptance
probability. The crucial independence is between the lists at distinct
vertices of the current group: their acceptance tests use disjoint
random choices. This is the independence needed by partial resampling.

We~apply the Harris--Srinivasan procedure directly to~these sampled lists. In the proof of their
Theorem~4.7, a fractional product distribution is first found and then
rounded by their Algorithm~2. Here that distribution is already supplied
by the sampling procedure just described, so the first step is
unnecessary. More explicitly, their Theorem~3.8 uses resampling
probabilities $p_{i,j}=\lambda_{i,j}/\sum_j\lambda_{i,j}$. The proof
of their Theorem~4.7 takes $\lambda_{i,j}=(1+\varepsilon)z_{i,j}$,
where $z$ is the supplied distribution, so these resampling
probabilities are exactly $z_{i,j}$. No new distribution has to be
computed. The rounding algorithm only needs to draw independent
values of its variables, find a violated row, and sample a subset of the
assignments currently contributing to that row. The row family here has
polynomial size, and its coefficients are computable by
\eqref{eq:progress-moment-formula} and the displayed formula for $\phi$.

The required subset sampling has a~direct dynamic program.
If~a~violated row currently has coefficients $a_1,\ldots,a_N\in[0,1]$, the Harris--Srinivasan rule chooses a~set $I$ of~a~specified size $d$ with probability proportional to~$\prod_{i\in I}a_i$. Put
\[
 E(j,\ell)=\sum_{\substack{I\subseteq\{1,\ldots,j\}\\|I|=\ell}}
                  \prod_{i\in I}a_i.
\]
The recurrence
\[
 E(j,\ell)=E(j-1,\ell)+a_jE(j-1,\ell-1),\qquad E(0,0)=1
\]
computes the table in $O(Nd)$ operations, taking impossible entries to
be zero. Starting with $j=N$ and $\ell=d$, include $j$ with probability
$a_jE(j-1,\ell-1)/E(j,\ell)$ and decrease $\ell$ when it is included;
then decrease $j$ and continue. This samples the required distribution.
Only the current assignments are used, not the set of all possible
values of any variable. This is the implementation principle in
Appendix~C of~\cite{HS}.
The recurrence separates subsets according to whether they contain
$j$. In the backward sampling step, the numerator is exactly the total
weight of subsets that do contain $j$, and the denominator is the total
weight of all remaining choices. Thus the table provides both the
normalizing constant and the conditional probabilities needed to draw
the subset. All its entries are obtained by rational additions and
multiplications of polynomially many input coefficients.

The proof of Theorem~4.7 of~\cite{HS} bounds the expected number of
resampling steps by a constant times the number of variables. It also
allows supplied upper bounds on the row means. For our load rows use
$m_k/(1-p_k)$; for an error row with unchanged contribution $C$ use
$(Q_v^{(k-1)}-C+1/n)/(1-p_k)$, with the scaling from
\eqref{eq:APtable}.
Thus the implementation uses computable upper bounds on~the retained means.
The thresholds in~the Harris--Srinivasan procedure may be~rounded upward to~rational numbers with additive error at~most one; this preserves its guarantees and changes our final error only by~a~constant.
Sampling a~list, evaluating the rows, and selecting a~resampling subset all take expected polynomial time.
The same is~therefore true of~the complete procedure for choosing retained lists, for fixed accuracy.

\subsection{Finding Disjoint Options}\label{sec:matching-algorithm}
The large-edge rounding step requires an~algorithm for the matching problem supplied existentially by~\Cref{prop:DP}.
Here the degree gap is~a~fixed fraction, which allows an~expected-polynomial implementation.
We~first regularize the instance without an~excessive increase in~its size, and then implement the random matching steps.

A~job belongs to~$A$, a~slot belongs to~$B$, and every option contains one job and two distinct slots.
Degrees and pair counts include repeated options.
The aim is to choose one option for every job without using any slot
twice. The degree gap says that each job has more choices than any one
slot has competing uses. We first reduce pair counts to a logarithmic
scale. We then repeatedly choose a small matching and remove its
vertices. The useful quantity is the ratio between the minimum number
of choices for a remaining job and the maximum number of uses of a
remaining slot. This ratio increases at each step. Once it is large,
independent choices followed by local resampling finish the construction.

The regularization lemma below is a tool for one such step. It fills
deficits at slots using new jobs, so that the random-step estimates can
use exact degrees. It is important that no new option is offered to an
old job: otherwise a matching in the enlarged instance might cover an
original job with an option that did not exist in our problem.

\Needspace{14\baselineskip}
\begin{lemma}
\label{lem:rank-three-regularization}
Let $G=(A,B)$ be a family of triples with $A\ne\varnothing$: every edge consists of one vertex in $A$ and two distinct vertices in $B$.
Suppose every vertex of $A$ has degree $a$, every vertex of $B$ has
degree at most $b$, and $a\ge b\ge1$ are integers. In polynomial
time, one can construct a family of triples
$G'=(A',B')$ containing $G$ as an induced subhypergraph such that
\begin{itemize}
\item every vertex of $A'$ has degree $a$ and every vertex of $B'$
has degree $b$;
\item no new edge meets an old vertex of $A$;
\item the maximum pair count is at most
$\max\{1,\Delta_2(G)\}$;
\item fewer than $4a$ vertices are added to $B$, and at most
$|B|/2+2b$ vertices are added to $A$.
\end{itemize}
The polynomial bound refers to the explicitly listed hypergraph.
\end{lemma}

\begin{proof}
The total degree deficit in $B$ is
\[
 T=\sum_{v\in B}(b-d_G(v))=b|B|-2a|A|.
\]
Let $g=\gcd(b,2a)$. Choose an integer $N$ with
\[
 2a\le N<4a,\qquad N\equiv-|B|\pmod{2a/g}.
\]
Such an $N$ exists because $2a/g\le2a$. Add $N$ new vertices to $B$,
each with deficit $b$. Then the total deficit is divisible by $2a$.
Put
\[
 m=\frac{T+Nb}{2a}.
\]
We have $b\le m\le |B|/2+2b$. Add $m$ new vertices to $A$.
Each new job must receive $a$ options, and each option uses two slot
incidences. Thus a new job absorbs exactly $2a$ units of deficit. The
choice of $N$ makes this division integral and also supplies enough new
jobs that each individual slot deficit can be spread over distinct jobs.

First find a simple bipartite graph $J$ between the new $A$ vertices
and all the $B$ vertices. Every new $A$ vertex must have degree $2a$;
each $B$ vertex must have degree equal to its deficit, denoted $d_v$.
These requirements are feasible: $0\le d_v\le b\le m$ and
$\sum_vd_v=2am$. For completeness, use a network with capacity $2a$
from the source to each new $A$ vertex, capacity one on every possible
new-$A$--$B$ edge, and capacity $d_v$ from $v$ to the sink. For
$0\le k\le m$,
\[
 \sum_v\min\{d_v,k\}\ge\frac{k}{m}\sum_vd_v=2ak.
\]
These inequalities verify every cut condition for flow $2am$.
An integer maximum flow therefore supplies $J$ in polynomial time.

Process the new $A$ vertices one at a time. A new vertex $u$ has
$2a$ distinct neighbors in $J$. Pair those neighbors, forbidding
any pair already contained in an edge of the growing hypergraph.
Before processing $u$, each of its neighbors still needs the one
incident edge supplied at this step. Its present degree is therefore
at most $b-1$, so it has at most $b-1$ forbidden partners. The allowed
graph on the $2a$ neighbors has minimum degree at least
\[
 2a-1-(b-1)=2a-b\ge a.
\]
It has a perfect matching. One elementary proof takes a largest
matching: if $x,y$ are unmatched, they are nonadjacent, and the sum
of their degrees is at least $2a$. But unless a three-edge augmenting
path exists, each matched pair accounts for at most two of these
incidences, giving a sum at most $2a-2$, a contradiction. A standard
polynomial matching algorithm finds the required pairing.

For every paired pair $v,w$, add the hyperedge $\{u,v,w\}$. This
gives $u$ exactly $a$ hyperedges and uses each of its incidence edges
in $J$ once. After all new vertices have been processed, every deficit
is filled. A pair consisting of a new $A$ vertex and a $B$ vertex
occurs at most once, since $J$ is simple. A pair of $B$ vertices is
never reused. Pairs meeting old $A$ vertices are unchanged. This proves
the pair count claim, and the construction adds only polynomially many
vertices and edges.
\end{proof}

\Needspace{10\baselineskip}
\begin{proposition}
\label{prop:algorithmic-rank-three-matching}
Fix $\eta,\zeta>0$. For sufficiently large integer $D$, an explicitly
listed family of triples $G=(A,B)$ satisfying
\[
 d(a)\ge(1+\eta)D\quad(a\in A),\qquad
 d(b)\le D\quad(b\in B),\qquad
 \Delta_2(G)\le D^{1-\zeta}
\]
has a matching covering every element of $A$ that can be found by a randomized algorithm in expected polynomial time. The degree threshold and running-time
constants may depend on $\eta,\zeta$. Each edge has exactly one vertex
in $A$ and two distinct vertices in $B$.
\end{proposition}

\begin{proof}
If $A$ is empty, return the empty matching. If no jobs remain at any later stage, return the matching already found. We henceforth assume the current job set is nonempty, so the initial $D$ is at most the number of explicitly listed input edges.
We describe how to make the no-forbidden-configuration case of the
Delcourt--Postle argument algorithmic under the fixed slack used here.
The only concentration input is their iterative bipartite lemma and
its variable-model local-lemma proof~\cite[Lemma~3.3 and Sections~5.5--6]{DP}.
We use the standard Moser--Tardos algorithm~\cite{MT}: sample all independent
variables and, while a bad event holds, resample its variables. Under
the symmetric local-lemma bound $ep(d+1)<1$, its expected number of
resamplings is polynomial in the number of explicitly checkable events.

The proof keeps a matching already chosen and a residual instance for
the uncovered jobs. Every accepted random step adds disjoint original
options to that matching. Its degree estimates ensure that the residual
instance still has many choices per job, while the competition at each
slot decreases faster. Resampling changes only the proposed step; we
accept its matching only after all of that step's tests pass.

Choose fixed rational $\theta$ with
$0<\theta<\min\{\eta/4,1/8\}$. Delete edges at vertices of $A$
until each has degree $a_*:=\lceil(1+4\theta)D\rceil$. This is
possible for large $D$, and all degrees are then at most $2D$.

\paragraph*{Reducing the degrees once}
The small-matching estimates below require logarithmic pair counts,
whereas the input only gives $D^{1-\zeta}$. A preliminary independent
sampling achieves this without spending a substantial part of the
degree gap. Ordinary degrees retain a growing expectation, but the
expected number of retained copies through a fixed pair tends to zero.
Choose a fixed rational $s$ with $0<s<\min\{\zeta/2,1/2\}$ and
put $p=2^{-\lceil(1-s)\log_2D\rceil}$. This is an exactly
computable dyadic probability; equivalently its exponent can be found
using integer powers, since $s$ is rational. Write $d=pD$.
Independently retain each edge with probability $p$.

Let $\delta=\theta/10$. The bad events are a retained degree below
$(1-\delta)pa_*$ at a vertex of $A$, a retained degree above
$(1+\delta)d$ at a vertex of $B$, and a pair count above
\[
 t=\left\lfloor\frac{\lfloor\log_2d\rfloor^2}{4}\right\rfloor.
\]
For large $D$, $t\le0.6(\log d)^2$. Chernoff bounds give degree-event
probabilities $\exp(-\Omega_{\theta}(d))$. A pair has expected
pair count at most $pD^{1-\zeta}\le D^{s-\zeta}$, which tends to zero;
the standard binomial upper-tail bound gives probability at most
$\exp(-\Omega_{s,\zeta}((\log D)^3))$ for a pair event.

Each event involves at most $2D$ edge variables. One edge variable
belongs to at most three degree events and three pair events.
Consequently the dependency degree is $O(D)$, and the symmetric
Moser--Tardos condition holds for large $D$. All events are checked
by counting retained incidences. This finds a retained graph in
expected polynomial time. Set
\[
 a_0=\lfloor(1-\delta)pa_*\rfloor,
 \qquad b_0=\lceil(1+\delta)d\rceil.
\]
For large $D$, its $A$ degrees are at least $a_0$, its $B$ degrees
are at most $b_0$, and $a_0/b_0\ge1+2\theta$. Its pair counts
are at most $0.6(\log d)^2$.

\paragraph*{Choosing a small random matching}
Suppose the current degree bounds are integers $a\ge b$, with
$1+\theta\le R:=a/b\le2000$, and pair counts at most
$(\log b)^2$. Fix rational constants
\[
 0<\gamma\le\min\{\theta/8,1/1000\},\qquad
 \varepsilon=\gamma/10^{12},\qquad w=1.
\]
These satisfy the inequalities in the Delcourt--Postle iterative
lemma for $r=3$, $g=2$, and $R\le2000$, once $b$ is large enough.
In particular,
\[
 \gamma\ge 96R(15R+3)\varepsilon,
 \qquad \varepsilon\ge b^{-1/18},\qquad R\le\log b.
\]
The last two conditions hold for all sufficiently large $b$ because the parameters are fixed.
The remaining hypotheses also hold: $\gamma\le1/144$, the pair count is~at~most $4\log^2b$, and the forbidden configuration hypergraph is~empty.

Trim the degrees in $A$ to $a$, and apply
\Cref{lem:rank-three-regularization} to make the $B$ degrees
exactly $b$. This extension has polynomial size and preserves the
pair-count bound. Independently select each edge with probability
$\varepsilon/b$, and independently mark it for deletion with probability
$\varepsilon$. These latter coins are called equalizing coins in
\cite[Section~5.1]{DP}: with an empty forbidden family and $w=1$, each
edge is deleted with the same probability $\varepsilon$.
Keep the selected edges that
meet no other selected edge and were not marked for deletion. They
form a matching. Delete its vertices and the edges marked for deletion.

The degree conclusions of the cited iterative lemma are
\[
 a'=a e^{-\varepsilon(3+\gamma)},\qquad
 b'=b e^{-\varepsilon(R+2-\gamma)}.
\]
The exponents explain the gain. An option at a job must survive at its
two slots, whereas an option at a slot must also survive competition at
its job, whose degree is $a=Rb$. The latter competition is stronger
because $R>1$. The deletion coins contribute equally to both sides.
Ignoring the small allowance $\gamma$, the ratio therefore gains a
factor $e^{\varepsilon(R-1)}$ in one step.

We briefly specify why the probabilistic step is algorithmic, rather
than relying only on its positive probability. Let $F_0$ be the marked
edges and $V'$ the vertices not covered by the proposed matching. For
each vertex $v$ of the regularized graph, define
\[
 D_v=\{e\ni v:e\notin F_0,\ e\setminus\{v\}\subseteq V'\}.
\]
Thus $|D_v|$ counts options that remain available at their other
vertices. It equals the residual degree if $v$ itself survives; its
definition does not require conditioning on that survival. The bad
events are
\[
 \mathcal A_u=\{|D_u|<a'\}\quad(u\in A),\qquad
 \mathcal A_v=\{|D_v|>b'\}\quad(v\in B).
\]
These are the one-sided events denoted $A_a,A_b$ in
\cite[Section~6]{DP}. Only surviving vertices need the resulting degree
guarantees, but testing the events at all vertices is convenient.
Each event can be checked from selection and deletion coins in a
fixed neighborhood in the incidence graph.
The concentration calculation in~\cite[Sections~5.4.3, 5.5, and~6]{DP} bounds their probabilities by
$\exp(- (\log b)^{3/2})$, while their dependency degrees are bounded
by a fixed power of $b$. With the forbidden family empty, no other
configuration events or probability-evaluation subroutines are needed.
Moser--Tardos therefore finds acceptable coins in expected polynomial
time. We may use weaker, rational degree thresholds, as specified next;
their bad events are subsets of these events and satisfy the same bounds.
The neighborhood bound is what makes this useful when the total graph
is very large: a fixed event does not depend on coins at distant jobs.
Since the regularized maximum degree is at most $2000b$, a neighborhood
of fixed radius contains only a fixed power of $b$ incidences. The
displayed failure probability is eventually smaller than the reciprocal
of any such power, giving the required local-lemma inequality.

\paragraph*{Exact thresholds and growth of the degree ratio}
The exponential expressions are convenient for analysis, but the
algorithm must make exact accept-or-reject decisions. We therefore
relax them slightly to rational thresholds before testing any event.
Put $\tau=\varepsilon\theta/100$. Compute rational numbers $L,U$
such that
\begin{align*}
 (1-2\tau)e^{-\varepsilon(3+\gamma)}
 &\le L\le(1-\tau)e^{-\varepsilon(3+\gamma)},\\
 (1+\tau)e^{-\varepsilon(R+2-\gamma)}
 &\le U\le(1+2\tau)e^{-\varepsilon(R+2-\gamma)}.
\end{align*}
Use the integer thresholds $a_+=\lfloor aL\rfloor$ and
$b_+=\lceil bU\rceil$ in the event tests. These relax the preceding
degree bounds. Taylor expansion with a certified remainder computes
the required exponential intervals in polynomial time; only a fixed
relative precision depending on $\theta$ is required. All actual
coin probabilities and event comparisons remain rational.

For large $b$, the integer rounding changes the relative bounds by
at most another $\tau$. Hence
\[
 \frac{a_+}{b_+}
 \ge R e^{\varepsilon(R-1-2\gamma)}
          \frac{1-3\tau}{1+3\tau}
 \ge R e^{\varepsilon\theta/2}.
\]
The ratio thus reaches $1000$ after at most
$\lceil2\log(1000)/(\varepsilon\theta)\rceil+1$ steps, a constant
depending only on $\eta$. Before and immediately after the last step
it remains below $2000$. Throughout these constantly many steps,
$b\ge c_{\eta}d$ for a positive constant $c_{\eta}$. Thus the
original pair-count bound $0.6(\log d)^2$ is at most
$(\log b)^2$ for large $D$, and every step's hypotheses hold.

\paragraph*{Discarding the regularization after each step}
To~keep the instance size polynomial, discard the auxiliary extension after each random matching step.
Retain only the matching edges and residual vertices belonging to~the graph before regularization. No regularization edge
meets an old $A$ vertex. Consequently an old $A$ vertex is removed
only when an original edge matches it, and the degree of each surviving
old $A$ vertex is unchanged by discarding dummy vertices. Degrees
in old $B$ can only decrease. Old $B$ vertices used by a discarded
dummy matching edge may simply be left unavailable in later steps.
Such a slot was already removed by the random step, so keeping it
unavailable cannot create a new loss beyond the degree estimates just
proved. It also cannot cover an old job: dummy edges contain only new
jobs. Thus every old job that disappears has genuinely been matched
by an original option.
Therefore the next current graph is always a subhypergraph of the
initial reduced graph. The extensions never accumulate or multiply.
Each individual extension has $O(|V(G)|+D)$ vertices and polynomially
many edges, and there are only constantly many extensions.

\paragraph*{Finishing}
Once $a\ge1000b$, trim each remaining job to exactly $900b$ options.
Independently choose one uniformly for each job. A bad event specifies
two options from different jobs that share a slot and are both chosen.
Its probability is $(900b)^{-2}$. A job belongs to at most
$(900b)(2b)$ such events, so an event has dependency degree at most
$3600b^2$. The symmetric local-lemma inequality holds:
\[
 e(900b)^{-2}(3600b^2+1)<1.
\]
Another Moser--Tardos run finds a collision-free choice. Together with
the original matching edges retained from the previous steps, it
covers every original job exactly once.

All input hyperedges, including repeated copies, are explicitly listed.
After the first trimming their number, the numbers of bad events,
the sizes of the regularized graphs, and the number of coin variables
are polynomial in that input size. Each event test and resampling uses
polynomial time. The parameters are fixed, all actual coin probabilities
are rational with polynomial encoding length, and the number of stages
is independent of input size. The standard expected resampling bounds
therefore give expected polynomial total time. Every returned object is
checked to be a matching covering every element of $A$, so every output is correct.
\end{proof}

\subsection{Combining the Procedures}
\label{subsec:simpalg}
For the simple-graph algorithm, use the approximate local packings in \Cref{sec:local-input-algorithms}.
The compression in~\Cref{lem:profiles} scans explicit slot sizes and assigns fractions greedily, and hence takes polynomial time for fixed accuracy.

The profile choices use the resampling procedure in \Cref{sec:computable-progress}.
The matching of large-edge options uses
\Cref{prop:algorithmic-rank-three-matching}: its job degrees
exceed its slot degrees by a fixed fraction, exactly as in
\eqref{eq:degrees}. The list of option triples has polynomial size.
Indeed its multiplicity scale is $\lceil n\rceil^6$, and the parameters
$n=\lceil\Delta_w(G)\rceil$ or an approximate local bin count are
bounded by a constant times the number of input edges, apart from a
constant depending on accuracy.

The small-edge choices use \Cref{lem:two}, with exact rational
linear algebra. Sorting, grouping parallel edges, finding ordinary
matchings, proper bipartite edge coloring, and the Euler-tour
construction also take polynomial time. Each random procedure has an
expected polynomial number of resampling steps and can check its
output. Thus these algorithms have expected polynomial running
time for every fixed accuracy.

Choose a fixed threshold large enough for the inequalities in these
proofs. Below that threshold use decreasing-weight greedy cleanup
directly, before invoking any resampling procedure. Its color count is
then bounded by a constant depending only on accuracy. This gives an
algorithm on every input and accounts for the additive constants in
the fixed-accuracy statements.

These are algorithms for each fixed accuracy.
Their running-time constants may depend on~that accuracy, as~allowed by~the uniform asymptotic guarantees.

\Needspace{10\baselineskip}

\clearpage
\section{The Bipartite Total-Weight Bound}
\label{sec:bip-multi}
This appendix proves the result in Section~\ref{sec:overview-bipartite}.
We group parallel edges, prove the missing-color estimates, and use the
rounding theorem from \Cref{sec:round-overview} to complete the coloring.

\subsection{Grouping Parallel Edges}\label{sec:parallel-grouping}
We first group parallel edges without changing any color load.
Let~$G$ be~a~weighted multigraph with $W_G(v)\le n$ at~every vertex.
Initially, each edge forms its own group.
For each fixed pair of~endpoints, repeatedly combine two groups whose weights are both at~most $1/2$.
The weight of~a~group is~the sum of~its edge weights.
Every group produced this way has weight at~most one.
When the process ends, at~most one group of~weight at~most $1/2$ remains for each pair of~endpoints.
Regard each group as~a~single edge.
Let $H$ consist of~the groups of~weight greater than $1/2$, and let $R$ consist of~the other groups.
Thus $H$ is~a~multigraph, whereas $R$ is~simple and has weights at~most $1/2$.
A~color assigned to~a~group will be~assigned to~all the original edges in~that group.
This keeps the load of~every color at~every vertex unchanged.

\begin{figure}[htbp]
	\centering
\begin{tikzpicture}[x=1cm,y=1cm,font=\small]
\node[vertex,label=left:$u$] (u) at (0,0) {};
\node[vertex,label=right:$v$] (v) at (3.8,0) {};
\foreach \yy/\ww/\cc in {1.6/0.35/inkblue,.8/0.30/inkblue,0/0.20/warmred,-.8/0.15/warmred,-1.6/0.10/warmred}{
 \draw[thick,\cc] (u)..controls (1.1,\yy) and (2.7,\yy)..node[pos=.5,fill=white,inner sep=1pt] {$\ww$}(v);
}
\node at (1.9,1.8) {Five parallel edges};
\draw[-{Latex},thick] (4.55,0)--(5.5,0);
\node[vertex,label=left:$u$] (uu) at (6.2,0) {};
\node[vertex,label=right:$v$] (vv) at (10,0) {};
\draw[very thick,inkblue] (uu)..controls (7.3,1.1) and (8.9,1.1)..node[pos=.5,above=3pt] {$0.65>1/2$}(vv);
\draw[very thick,warmred,dashed] (uu)..controls (7.3,-1.1) and (8.9,-1.1)..node[pos=.5,below=3pt] {$0.45\le1/2$}(vv);
\node at (8.1,1.8) {Two groups};
\end{tikzpicture}
	\caption{Parallel edges may be grouped and then treated as single weighted
	edges. Here $0.35+0.30=0.65$ and $0.20+0.15+0.10=0.45$.
	All original edges in a group receive the same color. The grouping rule
	leaves at most one group of weight at most $1/2$ for each pair of vertices.}
	\label{fig:packets}
\end{figure}
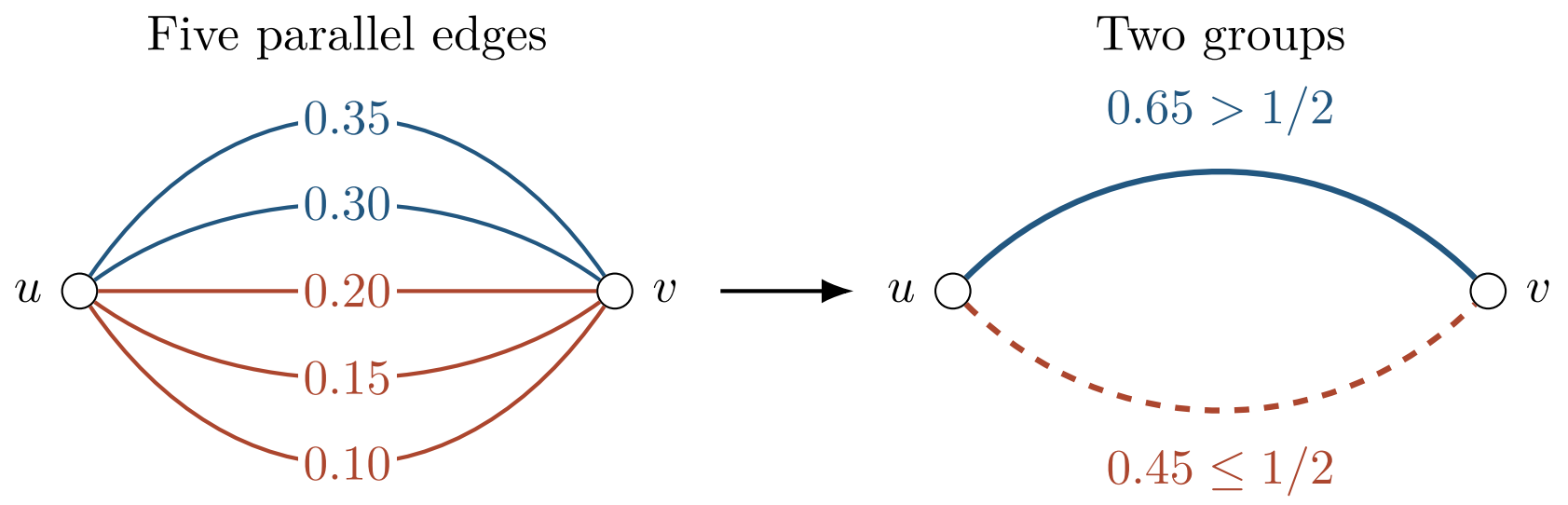

\subsection{Coordinating Missing Colors}\label{sec:missing-colors}\label{sec:mix-proof}
We~next study a~proper edge coloring of~a~multigraph $H$ with colors $1,\ldots,K$.
A~color is~\emph{missing at~$v$} if~no~edge at~$v$ has that color.
Write $F_v$ for the set of~missing colors and $f_v=|F_v|$.
We~will change the coloring while keeping every number $f_v$ fixed.

\begin{figure}[htbp]
	\centering
\begin{tikzpicture}[x=1cm,y=1cm,font=\small]
\foreach \idx/\shift/\title/\missing in {0/0/Before/b,1/6.3/After/a}{
 \begin{scope}[xshift=\shift cm]
 \node at (2.1,1.1) {\title};
 \foreach \xx/\name in {0/u,1.4/x,2.8/y,4.2/v}{
  \node[vertex,label=below:$\name$] (\name\idx) at (\xx,0) {};
 }
 \ifnum\idx=0
  \draw[inkblue,very thick] (u\idx)--node[above] {$a$} (x\idx);
  \draw[warmred,very thick,dashed] (x\idx)--node[above] {$b$} (y\idx);
  \draw[inkblue,very thick] (y\idx)--node[above] {$a$} (v\idx);
 \else
  \draw[warmred,very thick,dashed] (u\idx)--node[above] {$b$} (x\idx);
  \draw[inkblue,very thick] (x\idx)--node[above] {$a$} (y\idx);
  \draw[warmred,very thick,dashed] (y\idx)--node[above] {$b$} (v\idx);
 \fi
 \node at (2.1,-1) {Both ends miss $\missing$.};
 \end{scope}
}
\draw[-{Latex},thick] (4.75,0)--(5.7,0);
\end{tikzpicture}
	\caption{Exchanging two colors on a whole path component preserves a
	proper edge coloring. Only the missing colors at the two ends change. In the path shown,
	both ends miss the same color before the switch and the same color
	afterward. The proof in \Cref{sec:mix-proof} also permits even-length paths.}
	\label{fig:path-switch}
\end{figure}
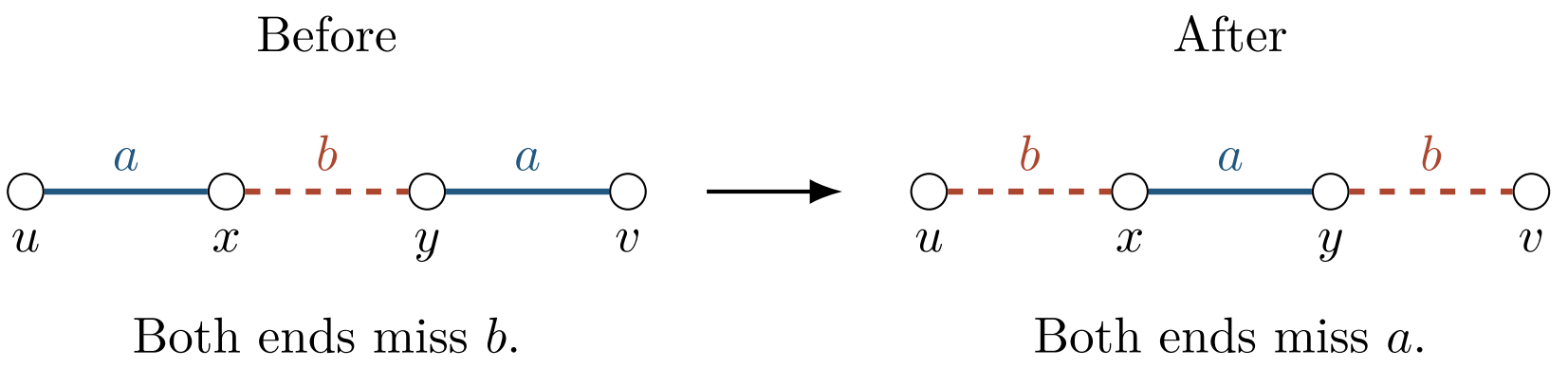

Fix two colors $a,b$ and consider the subgraph formed by~their edges.
Since the coloring is~proper, every vertex has degree at~most two in~the resulting graph.
Its components are paths, even cycles, and isolated vertices.
On~any whole path component, we~may exchange $a$ and $b$ on~every edge.
The coloring stays proper.
Only the missing colors at~the two ends of~the path change; see \Cref{fig:path-switch}.
This operation is~usually called a~\emph{Kempe switch}.

Let~$Q$ be~an~additional graph whose edges remain to~be~colored.
The next lemma assigns each edge fractions of~colors missing at~both endpoints.
It~controls the load of~each color and the total weight lost through edges whose assigned fractions sum to~less than one.
These estimates use first and second moments and do~not require a~positive lower bound on~the edge weights.

\Needspace{22\baselineskip}
\begin{lemma}\label{lem:weighted-moment-mix}
	Fix positive constants $\alpha,C,L$.
	Suppose a~multigraph without loops $H$ has a~proper coloring with an~even number $K\ge4$ of~colors.
	Let $F_v$ be~the colors missing at~$v$, and write $f_v=|F_v|$.
	Suppose
	\[
	f_v\ge\alpha n\quad\text{for every }v,\qquad K\le Cn.
	\]
	Let $Q$ be~another multigraph without loops on~the same vertex set, with nonnegative weights $q_e$, such that
	\[
	W_Q(v):=\sum_{e\ni v}q_e\le n,\qquad
	\sum_{e=uv}q_e\le L\quad\text{for every pair }u,v.
	\]
	For sufficiently large $n$, path switches in~$H$ give missing-color sets and fractions $x_{e,c}$ with the following properties:
	\begin{align}
		x_{e,c}>0&\ \Longrightarrow\ c\in F_u\cap F_v
		&&(e=uv),\label{eq:moment-support}\\
		0\le x_{e,c}&\le \frac{C}{\alpha^2n},\qquad
		\sum_cx_{e,c}\le1,\label{eq:moment-small}\\
		\sum_{e\ni v}q_ex_{e,c}
		&\le\frac{W_Q(v)}{f_v}
		+O\bigl(n^{-1/2}(\log n)^{3/2}\bigr),
		\label{eq:moment-load}\\
		\sum_{e\ni v}q_e\Bigl(1-\sum_cx_{e,c}\Bigr)
		&\le O\bigl(n^{3/4}(\log n)^{3/4}\bigr).
		\label{eq:moment-loss}
	\end{align}
	The constants and the required lower bound on~$n$ depend only on~$\alpha,C,L$, not on~the number of~vertices.
	No~positive lower bound on~$q_e$ is~required.
\end{lemma}

\begin{proof}[Proof of \Cref{lem:weighted-moment-mix}]
Write $a_{u,c}=\one_{\{c\in F_u\}}$. We choose a short list of pairings
of the colors and switch paths for each pairing. In addition
to averages of $a_{u,c}$, we keep track of averages of $a_{u,c}a_{u,d}$
for distinct colors $c,d$.

The first averages will bound the load proposed for each color.
The second averages have a different purpose: they bound the weighted
squared error in the number of common missing colors on an edge.
An average overlap alone would allow many edges to have no common
missing color, compensated by a few edges with large overlap.
The squared-error estimate rules out a large total weight of such
poorly served edges. We first find color pairings that mix these two
kinds of averages, then choose actual path switches that stay close
to the averages, and finally turn the overlaps into fractions.

For a list $h$, write $\|h\|^2=\sum_i h_i^2$; $h^T$ denotes the
same list written as a row rather than a column.

\paragraph*{Mixing colors and color pairs}
For a pairing $M$ of the $K$ colors, let $A_M$ average within each pair.
Let $B_M$ act on lists indexed by unordered pairs $\{c,d\}$ of distinct
colors. To apply $B_M$, independently exchange or keep the two colors
of every pair in $M$, apply the resulting permutation to the indices of
the list, and average over all these choices. Both matrices are symmetric,
have nonnegative entries, and have every row sum equal to one. Applying
either averaging operation twice has the same effect as applying it once.
Consequently
\[
 \|A_Mh\|^2=h^TA_Mh,\qquad
 \|B_Mh\|^2=h^TB_Mh.
\]
These identities let us measure mixing by a single quadratic expression.
The constant lists are unchanged, so the quantity to reduce is the
part of a list left after subtracting its average. The contraction
below is an expectation over the random pairing $M$; it need not
hold for each individual pairing.

Choose $M$ uniformly at random. Direct expansion gives
\[
 \E\|A_Mh\|^2=\frac{K-2}{2(K-1)}\|h\|^2
 \le\tfrac12\|h\|^2
 \quad\text{when }\sum_i h_i=0.
\]
We prove the same statement for $B_M$. An entry of $\E B_M$ depends only
on whether its two unordered-pair indices are equal, meet in one color,
or are disjoint. The respective entries are
\[
 p=\frac{K+2}{4(K-1)},\qquad
 q=\frac1{4(K-1)},\qquad
 r=\frac1{2(K-1)(K-3)}.
\]
For example, a pair $\{i,j\}$ is itself paired in $M$ with probability
$1/(K-1)$ and is then unchanged. Otherwise it is unchanged with
probability $1/4$, giving $p$. To change just $j$ to a specified third
color, those two colors must be paired, one swap must be made, and one
must be omitted, giving $q$. To reach a disjoint pair, there are two
possible pairings of the four colors, and both swaps must be made,
giving $r$.

Every symmetric array $(h_{ij})_{i\ne j}$ has a decomposition
\[
 h_{ij}=\bar h+(b_i+b_j)+r_{ij},
 \qquad \sum_i b_i=0,
 \qquad \sum_{j\ne i}r_{ij}=0\quad\text{for every }i,
\]
where
\[
 \bar h=\frac{1}{\binom K2}\sum_{i<j}h_{ij},
 \qquad b_i=\frac{1}{K-2}\sum_{j\ne i}(h_{ij}-\bar h).
\]
Indeed, the row sums after subtracting $\bar h$ add to zero, so
$\sum_i b_i=0$. The row sum of the array $b_i+b_j$ at $i$ is
$(K-2)b_i$, leaving every row sum of $r$ equal to zero.
The three arrays are mutually perpendicular under the inner product
$\langle x,y\rangle=\sum_{i<j}x_{ij}y_{ij}$:
the last two have total sum zero, and
\[
 \sum_{i<j}(b_i+b_j)r_{ij}
 =\sum_i b_i\sum_{j\ne i}r_{ij}=0.
\]

The matrix $\E B_M$ multiplies these three arrays, respectively, by
\[
 1,\qquad
 p+(K-4)q-(K-3)r=\frac{K-2}{2(K-1)},\qquad
 p-2q+r=\frac{K-2}{4(K-3)}.
\]
Here is the calculation behind the last two expressions. For a fixed
index $\{i,j\}$, separate the other indices into those meeting it in
one color and those disjoint from it. If $h_{ij}=b_i+b_j$ and
$\sum_i b_i=0$, the sums over these two classes are
$(K-4)h_{ij}$ and $-(K-3)h_{ij}$, respectively.
If every row sum of $h$ is zero, the two sums are instead
$-2h_{ij}$ and $h_{ij}$. Multiplying these class sums by $q$ and $r$
and adding the diagonal term $ph_{ij}$ gives the displayed factors.
Thus there are just two kinds of nonconstant behavior to control,
and both contract at the required rate.
The last two numbers are at most $1/2$ for even $K\ge4$. Therefore
\[
 \E\|B_Mh\|^2=h^T(\E B_M)h\le\tfrac12\|h\|^2
 \quad\text{when }\sum_{i<j}h_{ij}=0.
\]

Take $s=14\lceil\log_2K\rceil$ independent uniform pairings. Apply these
two contraction bounds successively to all unit coordinate lists after
subtracting their averages. The expected sum of the squared errors for
the two products is at most
\[
 \bigl(K+\tbinom K2\bigr)2^{-s}<K^{-10}.
\]
Thus some fixed list of pairings has the following property for both
products: it sends any list to within $K^{-4}$ times its original norm
of the constant list with the same average. Indeed, if $P$ is either product and $J$ replaces every entry by the
average, the Cauchy--Schwarz inequality gives
\[
 \|(P-J)x\|^2\le\|x\|^2\sum_j\|(P-J)e_j\|^2,
\]
where the $e_j$ are the unit coordinate lists. Fix such a list of pairings.

\paragraph*{Controlling the moments}
We now transfer this averaging calculation to the actual coloring.
We intend to propose fraction $K/(f_uf_v)$ for each common missing
color of an edge $uv$. The lists below record the neighbor contributions
to the resulting load and to the square of an edge's total proposal.
All three lists below are attached to a target vertex $v$ and sum
information from its neighbors $u$. The list $z_v$ measures how much
weighted demand can request a given color from those neighbors.
The lists $d_v$ and $Z_v$ will be the diagonal and off-diagonal terms
when we square the number of common missing colors: the diagonal
uses $a_{u,c}^2=a_{u,c}$, whereas the off-diagonal uses
$a_{u,c}a_{u,d}$.
During the switches define
\begin{align*}
 z_{v,c}&=\sum_{e=uv}\frac{q_e}{f_u}a_{u,c},\qquad
 d_{v,c}=\sum_{e=uv}\frac{q_e}{f_u^2}a_{u,c},\\
 Z_{v,\{c,d\}}&=
       \sum_{e=uv}\frac{q_e}{f_u^2}a_{u,c}a_{u,d}\quad(c\ne d).
\end{align*}
The numbers $f_u$ do not change during path switches. For a fixed
color pairing $M$, flip an independent fair coin on every nontrivial
two-color path and switch exactly when its coin says so. At one vertex,
the missing indicators for a paired pair of colors are either equal,
and unchanged, or unequal, and fairly exchanged. The coins for different
color pairs are independent. Hence, conditional on the coloring before
the stage, the means of the three displayed lists after the stage are
\[
 A_Mz_v,\qquad A_Md_v,\qquad B_MZ_v.
\]
Only the distributions at~one vertex enter this calculation; missing-color indicators at~different vertices may be~dependent.
In particular, the two ends of one switched path share the same coin.
We do not replace these dependent endpoint choices by independent ones.
Linearity of expectation is enough for the three displayed means.

We now choose the switches while keeping every coordinate close to this
mean. Process the $K/2$ color pairs one at a time. At one step the
independent choices are the path coins for the current pair of colors.
Switching within one pair does not change the edges or path components
of another pair, since their colors are disjoint. Thus the independent
coins for the entire stage are fixed before this process begins.
For every coordinate of $z,d$ or $Z$, track its conditional mean: average
only over the color pairs that have not yet been processed. A coordinate
$z_{v,c}$ or $d_{v,c}$ can change only when the pair containing $c$ is processed. A
coordinate $Z_{v,\{c,d\}}$ can change only when the pair containing $c$
or the pair containing $d$ is processed. Thus it changes at most twice.
If $c,d$ are paired together, their product at every vertex is unchanged.

For example, suppose $c$ belongs to the current color pair and $d$
belongs to an unprocessed pair $\{d,d'\}$. After the current choice,
the conditional contribution at a neighbor $u$ is
\[
 a_{u,c}\,\frac{a_{u,d}+a_{u,d'}}2.
\]
The second factor is already a fixed number: later switching the
$\{d,d'\}$-path cannot change its sum. Thus this expression is a
nonnegative function of the current path coin alone. When the pair
containing $d$ is eventually processed, the first factor has become
fixed. This is why processing the pairs one at a time converts the
products into sums to which \Cref{prop:HS} applies.

Whenever such a conditional mean changes, its new value is a constant
plus a sum of nonnegative functions of the independent path coins at
that step. For a product involving a color from an unprocessed pair,
replace that missing indicator by its mean, a number in $[0,1]$. If its
pair was processed earlier, use its fixed indicator instead. A path coin
changes missing-color indicators at no more than two vertices. Its contribution is therefore
at most
\[
 \frac{2L}{\alpha n}\quad\text{to a coordinate of }z,
 \qquad
 \frac{2L}{\alpha^2n^2}\quad\text{to a coordinate of }d\text{ or }Z.
\]
To control lower as well as upper errors, also use the sum obtained by
replacing each summand by its maximum possible value minus that summand.
All of these new summands are nonnegative.

The bound $L$ is used here to add together all edges of $Q$ joining
a target vertex to one end of the path. It bounds their total weight,
even if their number is large. This is the point at which a bound on
the weight of each parallel bundle replaces an ordinary-degree bound.

Apply \Cref{prop:HS}, multiplying the $z$ sums by a sufficiently
small fixed constant times $n$, and the $d,Z$ sums by a sufficiently small
fixed constant times $n^2$. Every summand then lies in $[0,1]$. The
expected size of each scaled sum is $O(n)$, since before scaling
\[
 z_{v,c}\le 1/\alpha,\qquad
 d_{v,c},Z_{v,\{c,d\}}\le1/(\alpha^2n).
\]
The same estimates apply to the sums used for lower errors.

The total contribution of any one value of a path coin to all scaled
sums is $O(n^2)$. Each of the coin's at most two ends
$u$ contributes to target vertices $v$ with total weight at most
$W_Q(u)\le n$. For the $z,d$ sums there are only the two current colors.
For the $Z$ sums one color must be current, and the other has at most
$K$ possibilities. The factors $1/f_u$ and $1/f_u^2$ cancel the respective
scalings up to fixed constants. Thus the $z,d$ contribution is $O(n)$ and
the $Z$ contribution is $O(Kn)=O(n^2)$. The sums used for lower errors
only double these bounds.

\Cref{prop:HS} now bounds all scaled errors by
$O(\sqrt{n\log n}+\log n)$. Each coordinate is affected at most twice
as the color pairs are processed. We can therefore complete one whole
stage with coordinate errors
\[
 O\bigl(n^{-1/2}\sqrt{\log n}\bigr)\quad\text{for }z,
 \qquad
 O\bigl(n^{-3/2}\sqrt{\log n}\bigr)\quad\text{for }d,Z.
\]
All constants are independent of the number of vertices. Notice that
no ordinary-degree bound was used.
There are $K/2$ decisions in a stage, but a fixed coordinate does not
pay for all of them. Its conditional mean changes only at the one or
two color pairs identified above. This is why the displayed errors
do not acquire an additional factor of $K$.

Apply this procedure for all $s=O(\log n)$ stages. Averaging does not
increase the largest coordinate error. Put
\[
 E_n=\sqrt n(\log n)^{3/2},\qquad
 D_v=\sum_{e=uv}\frac{q_e}{f_u}.
\]
The fixed sums of the lists are
\[
 \sum_cz_{v,c}=W_Q(v),\qquad
 \sum_cd_{v,c}=D_v,\qquad
 \sum_{c<d}Z_{v,\{c,d\}}=\frac{W_Q(v)-D_v}{2}.
\]
The mixing property of the chosen pairings, together with the accumulated
coordinate errors, consequently gives
\begin{align}
 z_{v,c}&=\frac{W_Q(v)}K+O(E_n/n),\label{eq:moment-first}\\
 d_{v,c}&=\frac{D_v}K+O(E_n/n^2),\label{eq:moment-diagonal}\\
 Z_{v,\{c,d\}}
 &=\frac{W_Q(v)-D_v}{K(K-1)}+O(E_n/n^2).
 \label{eq:moment-second}
\end{align}
The errors from imperfect averaging are smaller than the displayed errors:
the initial norms are $O(\sqrt n)$, $O(n^{-1/2})$, and $O(1)$,
respectively.

\paragraph*{Constructing the color fractions}
For an edge $e=uv$, write
\[
 \lambda_e=\frac{K|F_u\cap F_v|}{f_uf_v}.
\]
The reference overlap $f_uf_v/K$ is the intersection size one would expect
from two independent uniform subsets of sizes $f_u$ and $f_v$.
Thus $\lambda_e$ compares the actual overlap with this reference value.
The switched sets need not be independent; the moment estimates, rather
than such an independence assumption, justify what follows.
We do not require every $\lambda_e$ to be close to one. Instead we bound
the weighted sum of their squared errors at each vertex.
By \eqref{eq:moment-first},
\[
 \sum_{e\ni v}q_e\lambda_e
 =\frac K{f_v}\sum_{c\in F_v}z_{v,c}
 =W_Q(v)+O(E_n).
\]
Although $F_v$ was itself determined by the switches, this substitution
is valid: \eqref{eq:moment-first} holds for every color, so it can be
summed over any resulting set $F_v$. The same observation applies
to the pairs of colors in the second-moment calculation.
For the second moment, expand the square of $|F_u\cap F_v|$ and use
\eqref{eq:moment-diagonal} and \eqref{eq:moment-second}. We get
\begin{align*}
 \sum_{e\ni v}q_e\lambda_e^2
 &=\frac{K^2}{f_v^2}\left(
       \sum_{c\in F_v}d_{v,c}
       +2\sum_{\{c,d\}\subseteq F_v}Z_{v,\{c,d\}}\right)\\
 &=\frac{K D_v}{f_v}
   +\frac{K(f_v-1)}{f_v(K-1)}\bigl(W_Q(v)-D_v\bigr)+O(E_n)\\
 &\le W_Q(v)+O(E_n).
\end{align*}
For the last inequality, the first term is $O(1)$ because
$D_v\le1/\alpha$ and $K/f_v\le C/\alpha$; the coefficient of
$W_Q(v)-D_v$ is at most one because $f_v\le K$.
Combining the first and
second moments gives
\[
 \sum_{e\ni v}q_e(\lambda_e-1)^2=O(E_n).
\]
The Cauchy--Schwarz inequality therefore gives
\begin{equation}
 \sum_{e\ni v}q_e|\lambda_e-1|
 \le\sqrt{W_Q(v)\sum_{e\ni v}q_e(\lambda_e-1)^2}
 =O\bigl(n^{3/4}(\log n)^{3/4}\bigr).
 \label{eq:moment-total-error}
\end{equation}

Finally set
\[
 x_{e,c}=\frac{K}{f_uf_v}\,
       \frac{\one_{\{c\in F_u\cap F_v\}}}{\max\{1,\lambda_e\}}
       \qquad(e=uv).
\]
Before division by $\max\{1,\lambda_e\}$, each common missing color
receives the same proposed fraction $K/(f_uf_v)$, whose total on the
edge is $\lambda_e$. If this exceeds one, we scale all its fractions
down to total one. If it is below one, we keep them as they are and
charge the deficit $1-\lambda_e$ to the final cleanup. In either case
the deficit is at most $|\lambda_e-1|$, already controlled above.
These fractions have the required support and individual upper bound,
and their sum on an edge is $\min\{1,\lambda_e\}$. Equation
\eqref{eq:moment-total-error} proves \eqref{eq:moment-loss}. For a color
missing at $v$, division by $\max\{1,\lambda_e\}$ only decreases the load,
so \eqref{eq:moment-first} gives
\[
 \sum_{e\ni v}q_ex_{e,c}
 \le\frac K{f_v}z_{v,c}
 \le\frac{W_Q(v)}{f_v}+O(E_n/n).
\]
For a color present at $v$ the load is zero. This proves
\eqref{eq:moment-load} and completes the proof.
\end{proof}

Now the proof of \Cref{lem:overview-missing} follows easily. 
Apply \Cref{lem:weighted-moment-mix} with $q_e=2w_e$, size parameter $2n$, and $\alpha=\varepsilon/2$, $C=2$, $L=1$. Indeed, $W_Q(v)=2W_R(v)\le2n$, every pair has total charge at most one because $R$ is simple, $f_v\ge\varepsilon n$, and $K\le4n$. Dividing the charge estimates \eqref{eq:moment-load}--\eqref{eq:moment-loss} by two yields the claimed weight bounds; \eqref{eq:moment-support}--\eqref{eq:moment-small} give the support and small-proposal properties.

\subsection{Computing the Color Switches}\label{sec:mixing-implementation}
There are two separate searches in the proof. The first chooses a list
of pairings of the color names; this list is independent of the graph.
The second chooses which paths to switch in the current coloring.
We describe both searches because simply trying all path switches
would not give the required running time.

The list of color pairings in \Cref{lem:weighted-moment-mix}
can also be found explicitly. For each candidate list, form its
single-color and unordered-pair averaging matrices. If $J_d$ denotes
the $d$ by $d$ matrix all of whose entries are $1/d$, accept a list
when
\[
 \left\|A_s\cdots A_1-J_K\right\|_F^2
 +\left\|B_s\cdots B_1-J_{\binom K2}\right\|_F^2
 \le K^{-8}.
\]
Here the squared matrix norm means the sum of its squared entries,
so this is an exact rational test. It implies both operator bounds
used in the proof. Enumerating lists eventually finds one, by the
averaging argument already proved.

Alternatively, generate independent uniformly random pairings and
repeat this test. The expected squared error is less than $K^{-10}$,
so Markov's inequality gives success probability at least $1-K^{-2}$
on each trial. The matrices have polynomial dimension; each row of a
pair-averaging matrix has at most four possible images. Thus this
particular search has an expected polynomial-time implementation.

For a fixed pairing, the two-color paths and cycles are found by
following edges. Each path switch is a binary choice, and its moment
coefficients are explicit finite sums. The Harris--Srinivasan resampling
algorithm therefore chooses the switches directly.
More explicitly, process the color pairs in the order used in the proof.
For the current pair, use one binary variable per nontrivial path and
the scaled upper- and lower-error sums as the constraints of
\Cref{prop:HS}. Their summands are computable from the present coloring,
the weights $q_e$, and the fixed numbers $f_u$. Each successful step
fixes that pair's switches; it does not change the path components of
the remaining pairs in the stage. There are $O(K\log n)$ such steps
and $O(|V|K^2)$ moment coordinates. The contribution bounds proved
above are polynomial in $n$, so the cited resampling procedure has
expected polynomial running time for fixed $\alpha,C,L$.
Finally compute the sets $F_v$, their intersections on edges of $Q$,
and the fractions $x_{e,c}$ by the displayed formula. No search remains
in this last step.

\subsection{Completing the Coloring of the Groups}
\label{sec:finish-packets}
The missing-color and rounding lemmas now complete the coloring of~the grouped graph.
The argument uses only the number of~colors in~a~proper coloring of~$H$, so~it~applies to~both bipartite and general graphs.

\begin{lemma}\label{lem:finish-packets}
Fix $\tau,C>0$, and let $G$ be a weighted multigraph with $W_G(v)\le n$ at every vertex.
Apply the grouping in \Cref{sec:parallel-grouping}; let $H$ consist of groups heavier than $1/2$, and let $R$ consist of the remaining groups.
Suppose $H$ has a proper coloring with an even number $K$ of colors, where
\[
 (2+\tau)n\le K\le Cn.
\]
Then the original graph has a feasible coloring with $K+o(n)$ colors.
The error depends only on $\tau,C$, not on the number of vertices.
\end{lemma}
\begin{proof}
For every vertex,
\[
 d_H(v)\le2W_H(v)\le2(n-W_R(v)).
\]
Thus its number of missing colors satisfies
\[
 f_v=K-d_H(v)\ge\tau n+2W_R(v).
\]
Apply \Cref{lem:weighted-moment-mix} to the simple graph $R$, giving
an edge of weight $w_e$ charge $q_e=2w_e$. Its weighted degree in these
charges is at most $2n$, and the charge between any pair is at most one.
The scale parameter in that lemma can therefore be $2n$, with
$\alpha=\tau/2$ and color-bound constant $C$.
It supplies fractions supported on common missing colors, each
$O_{\tau,C}(1/n)$, whose sum on each edge is at most one. At every vertex,
\[
 \sum_{e\ni v}w_e\left(1-\sum_cx_{e,c}\right)=o(n).
\]
For each color, their local weighted load is at most
\[
 \sum_{e\ni v}w_ex_{e,c}
 \le\frac{W_R(v)}{f_v}+o(1)
 \le\frac12-\frac{\tau}{2C}+o(1)<\frac12
\]
for sufficiently large $n$. Here we used $f_v\le K\le Cn$.
\Cref{lem:half-load} makes each local vector an average of feasible
edge sets.

Fix $\delta>0$. \Cref{thm:round}, applied to the simple graph
$R$, gives a feasible partial coloring using the same $K$ colors and
leaving weight at most $\delta n+o(n)$ at every vertex. Each color used
on an edge of $R$ is missing at both endpoints in $H$, so the colorings
fit together. \Cref{lem:general-cleanup} colors the uncolored
edges with at most $4\delta n+o(n)$ new colors. Expand the groups back
into their original edges; this preserves every color load. Taking
$\delta$ arbitrarily small proves the lemma.
\end{proof}

\subsection{The Tight Leading Coefficient}\label{sec:bipartite-total-proof}
\begin{theorem}\label{thm:bipartite-total-weight}
Every bipartite weighted multigraph satisfies
\[
\chi'_w(G)\le(2+o(1))\Delta_w(G).
\]
The error is uniform in graph order, multiplicity, and minimum positive weight.
For each fixed accuracy and rational input, the construction has expected
polynomial running time. The leading coefficient two is best possible.
\end{theorem}
\begin{proof}
The empty graph needs no colors. Otherwise let
$n=\lceil\Delta_w(G)\rceil\ge1$.
Make the parallel-edge groups from \Cref{sec:parallel-grouping}.
	Every edge of~$H$ has weight greater than $1/2$, so~its ordinary maximum degree is~at~most $2n$.
	The graph $H$ is~bipartite.
	Fix $0<\varepsilon<1$, and let $K$ be~the least even integer at~least $(2+\varepsilon)n$.
	K\H{o}nig's theorem properly colors $H$ with at~most $2n\le K$ colors.
	For large $n$, $K\le4n$.
	\Cref{lem:finish-packets}, with $\tau=\varepsilon$ and $C=4$, gives $\chi'_w(G)\le(2+\varepsilon+o_\varepsilon(1))n$.
	Since $n=\Delta_w(G)+O(1)$, taking $\varepsilon$ arbitrarily small proves the bound.

The star example at the start of \Cref{sec:lower-bounds} proves sharpness.
\end{proof}

\begin{algorithm}[The bipartite total-weight construction]\label{alg:bipartite-total}
Input: a weighted bipartite multigraph and a desired accuracy.
Put $n=\lceil\Delta_w(G)\rceil$, and choose the palette and rounding
accuracies sufficiently small for the desired final error.
\begin{enumerate}[leftmargin=*,nosep]
\item Group parallel edges as in \Cref{sec:parallel-grouping}.
Let $H$ contain the groups heavier than $1/2$ and let $R$ contain
the remaining groups.
\item Properly color $H$ using K\H{o}nig's theorem. Add unused colors
to obtain the even palette size $K$ chosen in the proof.
\item Apply \Cref{lem:weighted-moment-mix}, as in the proof of
\Cref{lem:finish-packets}, to switch the coloring of $H$ and obtain
fractions for $R$. Construct the local feasible distributions using
\Cref{lem:algorithmic-local-load} with $b=1/2$ and $r=1$.
\item Run \Cref{alg:round} on $R$ and color its remaining edges with
fresh colors using \Cref{lem:general-cleanup}.
\item Give every original edge the color of its group.
Return the resulting coloring.
\end{enumerate}
\end{algorithm}

\clearpage
\section{Bin-Packing Bounds for Bipartite Multigraphs}
\label{sec:bp-bounds}\label{sec:bipartite-bounds}
This appendix gives the proofs for Section~\ref{sec:overview-bipartite-bp}.
We extend the rounding theorem to bounded multiplicity, obtain the harmonic
bound, and then improve it using compatible selections. The final charge
calculation and the $1.62$ deduction appear together below.

\subsection{Rounding with Bounded Multiplicity}
\label{sec:bounded-rounding}
Packet formation leaves a~bounded number of~edges between each pair of~vertices.
The next lemma extends the rounding theorem to~this setting.
Its proof tracks the two places where simplicity was used: the influence of~one vertex list and the pair counts of~the matching instance.

\begin{lemma}\label{lem:bp-mult-round}
Fix $L\in\mathbb N$ and $C_0,C_1,\delta>0$. Let $Q$ be a weighted multigraph without loops of multiplicity at most $L$, with $W_Q(v)\le n$,
and let $\mathcal C$ contain at most $C_1n$ colors. Suppose
\[
 0\le x_{e,c}\le C_0/n,\qquad \sum_cx_{e,c}\le1,
 \qquad (x_{e,c}:e\ni v)\in\mathcal K_v
 \quad\text{for every }v,c.
\]
Then, for sufficiently large $n$, a feasible partial coloring uses
color $c$ on $e$ only when $x_{e,c}>0$, and colors total weight at least
\[
 \sum_{e\ni v}w_e\sum_cx_{e,c}
 -\delta n-O_{L,C_0,C_1,\delta}(n^{3/4})
\]
at every vertex. The threshold and error constant are independent of
the number of vertices.
\end{lemma}
\begin{proof}
Use the proof of \Cref{thm:round}, keeping every parallel edge
as a separate edge. We give all changes where simplicity was used.

In the division into vertex groups, one source vertex contributes at
most $L$ to a row, while its sum over rows is still at most $n$.
Divide the rows by $L$ before applying \Cref{prop:HS} and
then multiply back. The error is still
$O_L(\sqrt{n\log n})$.

The list at a source $u$ now contributes the sum over its at most $L$
edges to a given target $v$. Its ordinary-slot contribution is
$O(L/n)$; its contribution to an error row is
$O(\sum_{e: e=uv}w_e)$. Increase the fixed scaling constant $B_0$ in
\eqref{eq:APtable} by a factor of $L$. The total contributions over
all target rows are still $O_L(n^2)$ for the scaled load rows and
$O_L(n)$ for the error rows, because sums over incident edges still
use the bound $W_Q(u)\le n$. Thus the Harris--Srinivasan bounds,
the induction on the vertex groups, and the $O(n^{3/4})$ total-error
bound hold with constants depending on $L$.

The variance calculation is for one edge at a time, using independent
colors, and is unchanged. The deletion of edges inside vertex groups
also remains valid: there are no loops, and all edges joining the same
two vertices are kept or deleted together.

In the construction of three-element sets, job--slot pairs have the
same bounds as before. A pair of ordinary slots can now belong to
options for at most $L$ original edges, giving $O(LN/n)$ copies.
A pair of extra slots can likewise belong to at most $L$ original
edges, giving $O(LNn^{-4/3})$ copies. Consequently the largest pair
count is still $O_L(Nn^{-2/3})$. With $N=\lceil n\rceil^6$, this is
$o(D^{1-\zeta})$ for every fixed $0<\zeta<1/9$, exactly as required
in the original proof. The job degrees, the slot degrees, the rounding
of small edges, and the final weight estimates do not use simplicity.
This proves the stated extension.
\end{proof}

\subsection{Coloring from Harmonic Charge}\label{sec:bp-details}
An~edge of~weight $1/(k+1)<w\le1/k$ has harmonic charge $1/k$ and can be~grouped with $k-1$ edges from the same class.
Each full packet then has total charge one and total weight at~most one.
We~color the full packets using K\H{o}nig's theorem and insert the remaining edges using \Cref{lem:weighted-moment-mix} and the bounded-multiplicity rounding theorem from \Cref{sec:bounded-rounding}.

\Needspace{17\baselineskip}
\begin{lemma}
	\label{lem:harmonic-load}
	For an~edge weight $0<w\le1$, put
	\[
	a(w)=\frac1{\lfloor1/w\rfloor},\qquad
	A(G)=\max_v\sum_{e\ni v}a(w_e).
	\]
	Every bipartite weighted multigraph $G$ with $W_G(v)\le n$ at~every vertex satisfies
	\[
	\chi'_w(G)\le A(G)+o(n).
	\]
	More precisely, for every $\eta>0$ there is~$N(\eta)$ such that the number of~colors is~at~most $A(G)+\eta n$ whenever $n\ge N(\eta)$.
	The threshold is~independent of~the number of~vertices, the number of~parallel edges, and the smallest weight.
\end{lemma}
\begin{proof}
	\textit{Group the very small edges.} The empty graph is~immediate.
	Notice that $w\le a(w)\le2w$, so~$A(G)\le2n$.
	Fix $0<\varepsilon<1$.
	Choose an~integer $T$ large enough that
	\[
	b=1/T\le\min\{1/4,\varepsilon/8\}.
	\]
	Both $b$ and $T$ stay fixed while $n$ tends to~infinity.

	For each endpoint pair, group its original edges of~weight below $b$ into packets of~weight in~$[b,2b)$, leaving at~most one tiny packet of~weight below $b$.
	Keep every original edge of~weight at~least $b$ separately.
	Give each non-tiny object of~weight $w$ charge $a(w)$; give each tiny packet charge $w/(1-b)$.
	Grouping preserves the weighted degree and can later be~undone without changing any color load.
	\Cref{fig:tiny-packets} shows this step.

	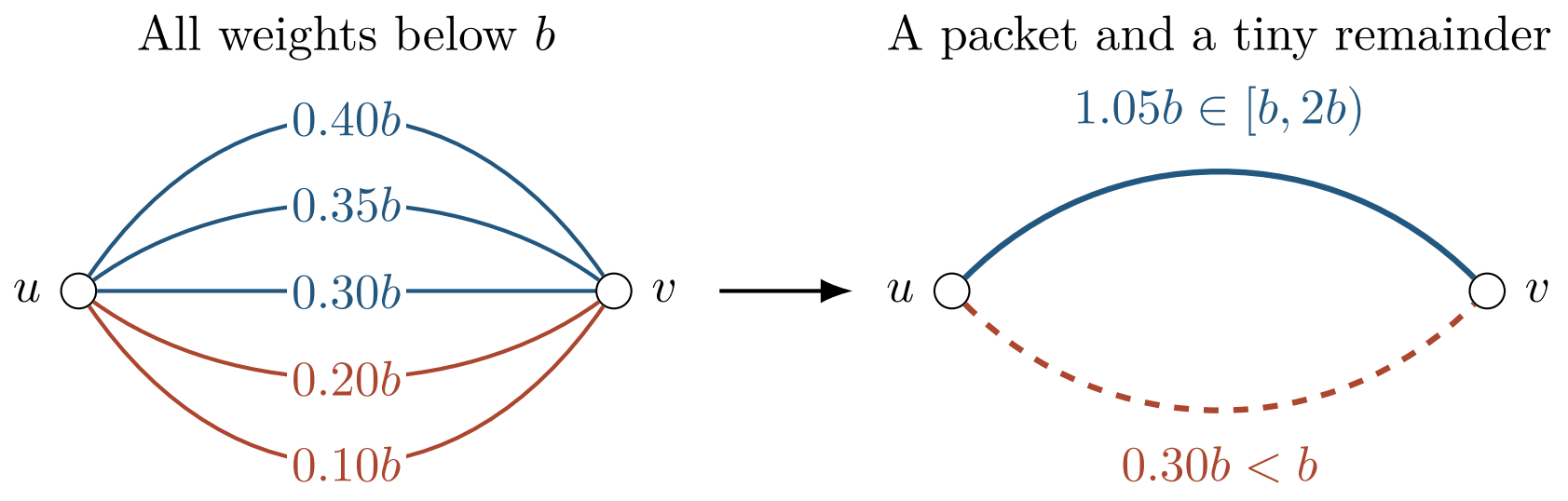
\begin{figure}[htbp]
		\centering
\begin{tikzpicture}[x=1cm,y=1cm,font=\small]
\node[vertex,label=left:$u$] (u) at (0,0) {};
\node[vertex,label=right:$v$] (v) at (3.8,0) {};
\foreach \yy/\ww/\cc in {1.6/0.40b/inkblue,.8/0.35b/inkblue,0/0.30b/inkblue,-.8/0.20b/warmred,-1.6/0.10b/warmred}{
 \draw[thick,\cc] (u)..controls (1.1,\yy) and (2.7,\yy)..node[pos=.5,fill=white,inner sep=1pt] {$\ww$}(v);
}
\node at (1.9,1.8) {All weights below $b$};
\draw[-{Latex},thick] (4.55,0)--(5.5,0);
\node[vertex,label=left:$u$] (uu) at (6.2,0) {};
\node[vertex,label=right:$v$] (vv) at (10,0) {};
\draw[very thick,inkblue] (uu)..controls (7.3,1.1) and (8.9,1.1)..node[pos=.5,above=3pt] {$1.05b\in[b,2b)$}(vv);
\draw[very thick,warmred,dashed] (uu)..controls (7.3,-1.1) and (8.9,-1.1)..node[pos=.5,below=3pt] {$0.30b<b$}(vv);
\node at (8.1,1.8) {A packet and a tiny remainder};
\end{tikzpicture}
		\caption{Collect small parallel edges until their total first reaches
		$b$. Every completed packet has weight below $2b$, and at most one
		tiny remainder is left between each pair. The original edge weights
		may be arbitrarily small; $b$ is a fixed grouping threshold.}
		\label{fig:tiny-packets}
	\end{figure}

	At~each vertex, the increase in~total charge is~at~most $4bn$.
	Indeed, a~completed small-edge packet of~weight $w<2b$ has
	\[
	a(w)<\frac{w}{1-w}\le(1+4b)w,
	\]
	where the last inequality uses $b\le1/4$.
	A~tiny packet also has charge at~most $(1+4b)w$.
	The original edges forming these packets had total charge at~least their total weight.
	Thus the increase is~at~most $4b$ times their original total weight, which is~at~most $4bn$.
	The charges of~the other original edges stay unchanged.

	\textit{Form full packets and color them.} Within each endpoint pair and each non-tiny class $k=\lfloor1/w\rfloor\in\{1,\ldots,T\}$, combine exactly $k$ objects at~a~time.
	Regard every full group as~one edge of~a~multigraph $H$.
	Its weight is~at~most one, and its charge is~exactly one.
	Let $R$ consist of~the ungrouped non-tiny objects and the tiny packets.
	It~has multiplicity at~most
	\[
	L=1+T(T-1)/2.
	\]
	There is~no~class-$1$ remainder, so~every edge of~$R$ has weight at~most $1/2$.
	Retain its charge
	\[
	q_e=\begin{cases}
		a(w_e),&w_e\ge b,\\
		w_e/(1-b),&w_e<b,
	\end{cases}
	\qquad Q_v=\sum_{e\in E_R(v)}q_e.
	\]
	In~particular, $w_e\le q_e\le1/2$, and
	\begin{equation}\label{eq:hload-charge}
		d_H(v)+Q_v\le A(G)+4bn\le A(G)+\varepsilon n/2.
	\end{equation}

	Let $K$ be~the least even integer at~least $A(G)+\varepsilon n$.
	For large $n$, $K\le4n$.
	K\H{o}nig's theorem gives a~proper $K$-coloring of~$H$.
	Its missing-color counts satisfy
	\[
	f_v=K-d_H(v)\ge Q_v+\varepsilon n/2.
	\]
	\textit{Make room for the remaining edges.} Apply \Cref{lem:weighted-moment-mix} to~$R$ with charges $q_e$, using size parameter $3n$.
	Indeed, $Q_v\le A(G)+4bn\le3n$, $f_v\ge(\varepsilon/6)(3n)$, $K\le2(3n)$, and the total charge between any pair is~at~most $L/2$.
	The lemma gives fractions supported on~the common missing colors, with $x_{e,c}=O_\varepsilon(1/n)$, $\sum_cx_{e,c}\le1$, and
	\[
	\sum_{e\in E_R(v)}q_ex_{e,c}\le\frac{Q_v}{f_v}+o(1),
	\qquad
	\sum_{e\in E_R(v)}q_e\left(1-\sum_cx_{e,c}\right)=o(n).
	\]
	All constants may depend on~$\varepsilon$ and its chosen $T$, but not on~the graph or~its smallest weight.
	Also $Q_v/f_v\le1-\varepsilon/8$, because $f_v\le K\le4n$.

	\textit{Check the local bin condition.} We~verify that these fractions are averages of~feasible sets at~one vertex and one color.
	For $2\le k\le T$, let $R_k(v)$ be~the non-tiny edges at~$v$ with $\lfloor1/w_e\rfloor=k$, and put
	\[
	p_k=\max\left\{\frac1k\sum_{e\in R_k(v)}x_{e,c},
	\max_{e\in R_k(v)}x_{e,c}\right\}.
	\]
	For tiny edges let
	\[
	p_0=\max\left\{\sum_{e\in E_R(v):\,w_e<b}\frac{w_e}{1-b}x_{e,c},
	\max_{e\in E_R(v):\,w_e<b}x_{e,c}\right\}.
	\]
	Maxima over empty sets are zero.
	The sum of~these probabilities is~at~most the local charge load plus $T\max_e x_{e,c}$, hence is~at~most one for sufficiently large $n$.
	Choose case $k$ with probability $p_k$, case zero with probability $p_0$, and the empty set with the remaining probability.
	In~case $k$, \Cref{lem:bp-classchoice} supplies a~set of~at~most $k$ edges with conditional inclusion probabilities $x_{e,c}/p_k$.
	In~case zero, these conditional probabilities have weighted sum at~most $1-b$ and every item has weight below $b$.
	\Cref{lem:algorithmic-local-load}, with $r=1$, therefore supplies
	a distribution on sets of weight at most one with these conditional
	inclusion probabilities.
	This proves the local condition required by~\Cref{lem:bp-mult-round}.

	\textit{Choose whole edges and finish.} Because $q_e\ge w_e$, the total original weight not proposed at~each vertex is~$o(n)$.
	Apply \Cref{lem:bp-mult-round} with any fixed error parameter $\delta>0$.
	It~colors all but weight at~most $\delta n+o(n)$ at~every vertex.
	Its colors are missing from $H$ at~both ends, so~the two colorings fit together.
	All remaining edges have weight at~most $1/2$.
	\Cref{lem:cleanup} colors them with at~most $2\delta n+o(n)+1$ additional colors.
	Expanding the groups and packets gives
	\[
	\chi'_w(G)\le A(G)+(\varepsilon+2\delta)n+o(n).
	\]
	Choose $\varepsilon,\delta$ sufficiently small compared with the given $\eta$, then choose $T$, and finally take $n$ sufficiently large.
	Every choice depends only on~$\eta$.
	This proves the stated uniform bound.
\end{proof}

\subsection{Constructing the Harmonic Local Distributions}\label{sec:harmonic-inputs}
Fix a vertex $v$ and a color, and use the classes, threshold $b$, and
probabilities $p_k,p_0$ from the proof of \Cref{lem:harmonic-load}.
Abbreviate $R_k(v)$ by $R_k$ and the proposed fractions by $x_e$.
The class probabilities sum to at most one; the remaining probability
chooses the empty set. We construct the conditional distributions with
polynomial support.

For $p_k>0$, the vector $(x_e/p_k:e\in R_k)$ has coordinates at most
one and sum at most $k$. The second part of
\Cref{lem:algorithmic-local-load} supplies the required distribution
on sets of at most $k$ items; their weights are at most $1/k$, so they fit.
For $p_0>0$, the conditional tiny-item vector has weight at most $1-b$.
Use the first part of the lemma with $r=1$. In each case the unconditional
inclusion probability of $e$ is exactly $x_e$. There are only a fixed
number of classes, and the total support size is polynomial.
The same bounded-load construction applies whenever the proof uses
fractions of load at most $r-b$ for items of weight at most $b$.
The resulting local distributions have explicit polynomial support and can be~constructed from the proposed fractions.
For the packet constructions, use these local distributions and those in
\Cref{sec:feasible-choices}.

\subsection{The Harmonic Bound}
Let
\[
r_1=1,\qquad r_{j+1}=r_j(r_j+1),\qquad
h_\infty=\sum_{j\ge1}\frac1{r_j}=1.691030206757\ldots.
\]
The constant $h_\infty$ is~the supremum of~the harmonic charge that fits in~one unit bin.
This is the classical harmonic weight bound of Lee and Lee
\cite[Theorem~1 and Corollary~1]{LeeLee85}; see also~\cite{Epstein}.
For completeness, we~prove the upper bound needed here.

\begin{lemma}\label{lem:general-bin-charge}
	Put $a(w)=1/\lfloor1/w\rfloor$.
	Every finite list of~positive weights with sum at~most one has total $a$-charge at~most $h_\infty$.
\end{lemma}
\begin{proof}
	For every positive integer $t$,
	\[
	a(w)\le\frac{t+1}{t}w\qquad(0<w\le1/t).
	\]
	Indeed, $k=\lfloor1/w\rfloor\ge t$ and $1/k<((k+1)/k)w\le((t+1)/t)w$.

	Start with remaining total weight at~most $1/r_1=1$.
	At~step $j$, suppose the remaining total is~at~most $1/r_j$.
	If~an~item exceeds $1/(r_j+1)$, it~has charge exactly $1/r_j$.
	Remove it.
	The new total is~less than $1/r_j-1/(r_j+1)=1/r_{j+1}$, so~the next step applies.
	Otherwise all remaining items have weight at~most $1/(r_j+1)$.
	Their total charge is~at~most
	\[
	\frac{r_j+2}{r_j+1}\frac1{r_j}
	=\frac1{r_j}+\frac1{r_{j+1}}.
	\]
	Including the earlier removed items, this is~at~most $h_\infty$.
	The procedure must stop, since each step that continues removes an~item from a~finite list.
	An~empty remainder also satisfies the bound.
\end{proof}

\begin{theorem}
	\label{thm:bp-harmonic-first}
	Every bipartite weighted multigraph $G$ satisfies
	\[
	\chi'_w(G)\le (h_\infty+o(1))\bp(G).
	\]
	The error is~independent of~graph order, multiplicity, and minimum weight.
	For rational weights and each fixed accuracy, the coloring can be~found in~expected polynomial time.
\end{theorem}
\begin{proof}
	Put $n=\bp(G)$.
	Pack each vertex's incident edges into $n$ bins.
	\Cref{lem:general-bin-charge} gives incident charge at~most $h_\infty n$.
	The harmonic-charge coloring lemma (\Cref{lem:harmonic-load}) colors the graph using at~most its largest incident charge plus $o(n)$ colors.
	Its error is~uniform even for arbitrarily small weights, and its construction is~algorithmic.
	This proves the result.
\end{proof}

\subsection{The Five-Thirds Bound}\label{sec:short-five-thirds}

\begin{theorem}
	\label{thm:bp-five-thirds-short}
	Every bipartite weighted multigraph $G$ satisfies
	\[
	\chi'_w(G)\le\bigl(\tfrac53+o(1)\bigr)\bp(G).
	\]
	The error is~independent of~the number of~vertices, the multiplicities, and the smallest edge weight.
\end{theorem}

The improvement comes from coloring a~small fraction of~two compatible weight ranges separately.
This reduces the harmonic charge of~the remaining edges to~at~most $3/2$ per original bin, while the selected edges use only $1/6$ of~a~color per bin asymptotically.

\begin{proof}
	Put $n=\bp(G)$ and consider the two ranges
	\[
	A=(1/2,5/8],\qquad B=(1/3,3/8].
	\]
	An~$A$-edge and a~$B$-edge fit together, as~do~two $B$-edges.
	At~each vertex, $n$ local bins contain at~most $n$ edges from $A$ and at~most $2n$ edges from $A\cup B$.

	\textit{Selecting a~sixth of~each range.} Choose one sixth of~each range, to~within one edge at~every vertex.
	In~the graph of~one range, divide the incidences at~each vertex into groups of~six and one possible smaller group.
	Split the vertex into one copy per group.
	The resulting graph is~bipartite and has maximum degree six.
	Properly color it~with six colors by~K\H{o}nig's theorem, and keep any one color.
	Every full group contributes one selected edge; the smaller group contributes zero or~one.
	Apply this construction separately to~$A$ and $B$.

	Let $F$ be~the selected union and put $K=\lceil n/6\rceil+1$.
	At~every vertex, $F$ has at~most $K$ incident $A$-edges and at~most $2K$ incident edges in~total.
	Split each vertex into two copies, put all its $A$-edges at~the first copy, and distribute its $B$-edges so~that each copy has degree at~most $K$.
	Properly $K$-color the graph of~copies.
	After joining the copies back together, each color has at~most two edges at~a~vertex, with at~most one from $A$.
	Every such pair fits, so~$F$ needs at~most $n/6+O(1)$ colors.

	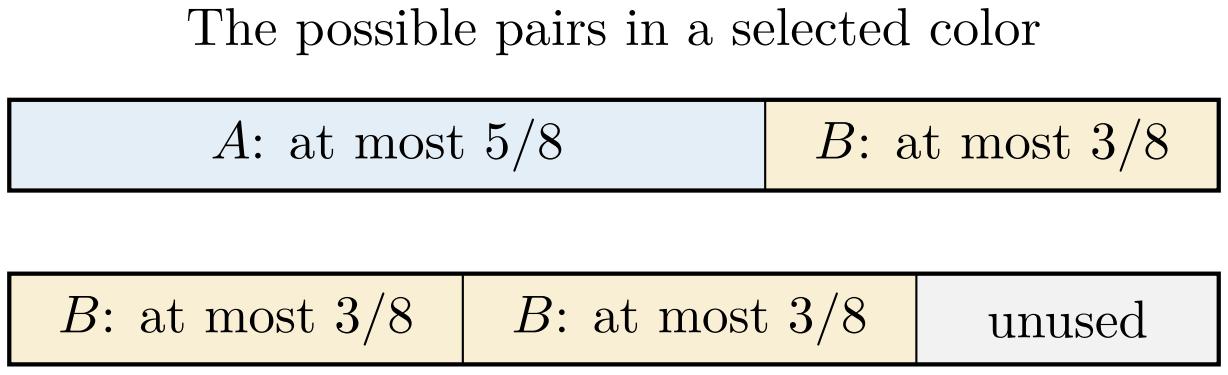
\begin{figure}[htbp]
		\centering
\begin{tikzpicture}[x=1cm,y=1cm,font=\small]
\fill[softblue] (0,.8) rectangle (5,1.4);
\fill[softgold] (5,.8) rectangle (8,1.4);
\draw[thick] (0,.8) rectangle (8,1.4);
\draw (5,.8)--(5,1.4);
\node at (2.5,1.1) {$A$: at most $5/8$};
\node at (6.5,1.1) {$B$: at most $3/8$};
\fill[softgold] (0,-.35) rectangle (6,.25);
\fill[black!5] (6,-.35) rectangle (8,.25);
\draw[thick] (0,-.35) rectangle (8,.25);
\draw (3,-.35)--(3,.25);
\draw (6,-.35)--(6,.25);
\node at (1.5,-.05) {$B$: at most $3/8$};
\node at (4.5,-.05) {$B$: at most $3/8$};
\node at (7,-.05) {unused};
\node at (4,1.85) {The possible pairs in a selected color};
\end{tikzpicture}
		\caption{Splitting each vertex into two copies allows these pairs.
		Both have total weight at most one; singleton and empty colors are
		also allowed.}
		\label{fig:short-five-thirds}
	\end{figure}

	\textit{Charge the remaining edges.} Define
	\[
	q(w)=
	\begin{cases}
		\frac56 a(w),&w\in A\cup B,\\
		a(w),&\text{otherwise}.
	\end{cases}
	\]
	The discounted charge satisfies
	\begin{equation}\label{eq:short-five-thirds-light}
		0<w\le\frac12\quad\Longrightarrow\quad q(w)\le\frac43w.
	\end{equation}
	For $w\le1/3$, put $k=\lfloor1/w\rfloor\ge3$; then $a(w)=1/k\le(4/3)w$.
	For $1/3<w\le3/8$, the discounted charge is~$5/12\le(4/3)w$.
	For $3/8<w\le1/2$, it~is~$1/2\le(4/3)w$.

	Every original unit bin now has charge at~most $3/2$.
	If~it~has no~edge heavier than $1/2$, equation~\eqref{eq:short-five-thirds-light} gives charge at~most $4/3$.
	Otherwise let $h$ be~the weight of~its unique edge heavier than $1/2$.
	The two possibilities give
	\[
	\begin{array}{ll}
		h\le5/8:&\displaystyle
		q(h)+\sum_{w\ne h}q(w)
		\le\frac56+\frac43(1-h)\le\frac32,\\[5pt]
		h>5/8:&\displaystyle
		q(h)+\sum_{w\ne h}q(w)
		\le1+\frac43(1-h)\le\frac32.
	\end{array}
	\]
	Here the sum means the other items in~this bin, with multiplicities.

	Let $R=G-F$.
	Selecting a~sixth of~$A$ removes a~sixth of~its charge, to~within one at~each vertex.
	Selecting a~sixth of~$B$ does the same, to~within $1/2$.
	Consequently,
	\[
	\sum_{e\in E(R),\ e\ni v}a(w_e)
	\le\sum_{e\ni v}q(w_e)+\frac32
	\le\frac32n+\frac32.
	\]
	The harmonic-charge lemma colors $R$ with $(3/2+o(1))n$ colors, since $W_R(v)\le n$.
	Use fresh colors for $R$ and combine the two colorings:
	\[
	\chi'_w(G)\le\frac n6+\frac{3n}{2}+o(n)
	=\bigl(\tfrac53+o(1)\bigr)n.
	\]
	All constant losses above are independent of~the graph, so~the uniform error bound follows from that of~the harmonic-charge lemma.
\end{proof}

\subsection{The Bound with Coefficient 1.62}\label{sec:bipartite-improved}

We~now combine several compatible selections of~edges.
Write $a(w)=1/\lfloor1/w\rfloor$ and $A(F)=\max_v\sum_{e\ni v,\ e\in F}a(w_e)$ for an~edge set~$F$.
By~\Cref{lem:harmonic-load}, every bipartite weighted multigraph~$F$ of~weighted degree at~most~$n$ satisfies
\begin{equation}\label{eq:bip-harmonic-input}
	\cw(F)\le A(F)+o(n),
\end{equation}
with a~uniform error term.
We~will apply this bound to~the remainder after the selections.
The first lemma ensures that the selections can be~made simultaneously on~whole edges while controlling their counts and the remainder charge at~each vertex.
With error at most four, it is a special case of
\cite[Theorem~10]{BLRT26}: the elements are allowed edge-label pairs,
a basis of the partition matroid chooses one label per edge, and each
element contributes at most one to each of two endpoint constraints.
The direct proof below gives the strict error bound stated here.
The same constraint-dropping method appears for edge-count decompositions
in~\cite[Section~3]{FS}.

\begin{lemma}
	\label{lem:bip-balanced-labels}
	Let $G$ be~a~weighted multigraph, and suppose each edge~$e$ has finitely many possible labels, including a~remainder label~$0$, with probabilities $p_{e,j}$ summing to~one.
	For each positive label $j$, its eligible edges are divided into sets $H_j,S_j$.
	One can assign exactly one allowed label to~each physical edge so~that, at~every vertex, each selected $H_j$ count, each selected $S_j$ count, and the remainder harmonic charge differs from its prescribed expectation by~less than four.
	The assignment is~computable in~polynomial time when the edge weights and prescribed probabilities are rational.
\end{lemma}
\begin{proof}
	Use variables $x_{e,j}$, initially equal to~$p_{e,j}$, and keep $\sum_jx_{e,j}=1$ for each edge.
	Also keep one equation for each quantity in~the statement.
	Repeatedly move to~a~vertex of~the current feasible polytope and fix every variable that is~zero or~one.
	An~unfixed edge has at~least two fractional variables.
	If~there are $M$ fractional variables, at~most $M/2$ edge equations remain.
	Every fractional variable occurs in~at~most two endpoint equations besides its edge equation, and its coefficient there is~at~most one: it~is~either a~count coefficient or~$a(w_e)$.
	Hence the endpoint equations contain at~most $2M$ fractional-variable appearances.
	If~every endpoint equation contained at~least five such variables, there would be~at~most $2M/5$ endpoint equations, for fewer than $M$ equations altogether.
	This contradicts extremality of~a~point with $M$ fractional variables.
	Drop an~endpoint equation containing at~most four fractional variables and repeat.
	Once dropped, its final value can change by~less than four; equations retained to~the end hold exactly.
	The edge equations are never dropped, so~every edge receives one label.
	Each iteration fixes a~variable or~drops an~equation, so~the process terminates after polynomially many iterations.
	The existence argument applies to~real probabilities.
	For rational input, every linear program has rational data of~polynomial encoding length, which gives the claimed running time.
\end{proof}

The next lemma uses the vertex-splitting and K\H{o}nig-coloring method
from~\Cref{lem:cleanup} for two compatible weight ranges.
\begin{lemma}
	\label{lem:bip-two-copy}
	Let the edge set of~a~bipartite multigraph $J$ be~partitioned into $H,S$, with $H$-weights in~$(1/2,h]$ and $S$-weights in~$(0,b]$.
	If~$2b\le1$ and $h+b\le1$, then $J$ has a~feasible coloring with
	\[
	L=\max_v\max\left\{d_H(v),
	\left\lceil\frac{d_H(v)+d_S(v)}2\right\rceil\right\}
	\]
	colors, constructible in~polynomial time.
\end{lemma}
\begin{proof}
	Make two copies of~each vertex.
	Put all $H$ incidences at~the first copy and distribute $S$ incidences so~each copy has degree at~most $L$; the first has room since $d_H(v)\le L$, and total room is~$2L$.
	Make the choices independently at~the two endpoints of~each edge.
	The copied graph remains bipartite and has maximum degree $L$.
	K\H{o}nig's edge-coloring theorem gives a~proper $L$-coloring.
	After identifying copies, a~color has at~most two edges at~a~vertex.
	If~one is~in~$H$, the other is~in~$S$; otherwise both are in~$S$.
	The two compatibility inequalities make both cases feasible.
\end{proof}

\begin{lemma}
	\label{lem:bip-compatible-selections}
	Let $G$ be~a~bipartite weighted multigraph, and fix finitely many selections, indexed by $j$.
	Selection $j$ has a~positive rational cost $f_j$, a~heavy range $H_j\subseteq(1/2,h_j]$, a~small range $S_j\subseteq(0,b_j]$, and a~nonnegative function $\rho_j$ on~$S_j$, where $2b_j\le1$ and $h_j+b_j\le1$.
	Suppose every unit bin $B$ satisfies
	\begin{equation}\label{eq:bip-mode-budget}
		\bigl|\{e\in B:w_e\in H_j\}\bigr|+\sum_{\substack{e\in B\\w_e\in S_j}}\rho_j(w_e)\le2.
	\end{equation}
	An~edge of~weight $w$ is~assigned probability $f_j$ for selection $j$ if~$w\in H_j$, probability $f_j\rho_j(w)$ if~$w\in S_j$, and zero otherwise; write the resulting probability as~$p_j(w)$.
	Assume $\sum_jp_j(w)\le1$ for all $w$.
	Set
	\[
	q(w)=a(w)\bigl(1-\sum_jp_j(w)\bigr).
	\]
	If~every unit bin has $q$-charge at~most $Q$, then
	\[
	\cw(G)\le\bigl(Q+\sum_jf_j+o(1)\bigr)\bp(G).
	\]
	The error is~uniform in~graph order, multiplicity, and minimum weight.
\end{lemma}
\begin{proof}
	Put $n=\bp(G)$ and apply \Cref{lem:bip-balanced-labels} to~all edges.
	Since each vertex has an~$n$-bin packing and each bin has at~most one heavy edge, the expected selected $H_j$ count is~at~most $f_jn$.
	Equation~\eqref{eq:bip-mode-budget} similarly bounds the expected number of~edges assigned to~selection $j$ by~$2f_jn$.
	After rounding, let $H'_j$ and $S'_j$ be~the subgraphs of~edges assigned to~selection $j$ whose weights lie in~$H_j$ and $S_j$, respectively.
	Then
	\[
	d_{H'_j}(v)<f_jn+4,\qquad
	d_{H'_j}(v)+d_{S'_j}(v)<2f_jn+8.
	\]
	\Cref{lem:bip-two-copy} colors selection $j$ with at~most $\lceil f_jn\rceil+4$ fresh colors.
	The prescribed remainder charge at~each vertex is~$\sum_{e\ni v}q(w_e)\le Qn$, so~balanced labeling gives $A(R)<Qn+4$ for the remainder $R$.
	Its weighted degree is~at~most $n$, and~\eqref{eq:bip-harmonic-input} colors it~with $Qn+4+o(n)$ fresh colors.
	Summing the disjoint palettes proves the claim.
\end{proof}

For a selection whose small edges lie in harmonic class $k\ge2$, take
$H=(1/2,h]$, $S=(1/(k+1),b]$ with $b\le1/k$, and $\rho(w)=2/k$.
A bin without an $H$ edge has at most $k$ selected small edges.
With an $H$ edge it has at most $\lfloor k/2\rfloor$, so
\eqref{eq:bip-mode-budget} holds in both cases.
The three selections used below are all of this form.

\subsubsection{The Three Selections and Their Charge}
\label{app:bip-three-selections}
The following three compatible selections give the charge estimate for
\Cref{thm:bipartite-improved}.

A row selects a heavy edge with probability $f$ and a small edge with
probability $2f/k$.
\begin{table}[htbp]
\centering
\begin{tabular}{c|c|c|c}
$k$ & Heavy range $H$ & Small range $S$ & Cost $f$\\\hline
2 & $(1/2,5/8]$ & $(1/3,3/8]$ & $1/20$\\
2 & $(1/2,7/12]$ & $(1/3,5/12]$ & $1/10$\\
3 & $(1/2,17/24]$ & $(1/4,7/24]$ & $1/12$
\end{tabular}
\caption{Three compatible selections with total cost $7/30$.}
\label{tab:bip-three-selections}
\end{table}
Every row satisfies $h+b=1$ and $b<1/2$, and the total selection
probability at any weight is at most $1/20+1/10+1/12=7/30$.
Let $p(w)$ be this total probability and put
$q(w)=a(w)(1-p(w))$, where $a(w)=1/\lfloor1/w\rfloor$.

\begin{lemma}\label{lem:bip-three-choice-charge}
Every unit bin has total $q$-charge at most $104/75$.
\end{lemma}
\begin{proof}
We compare the remaining charge with a fixed multiple of the actual
weight. The weight in one bin is at most one, so this part is easy
to sum. Only the positive excess above that multiple requires a
case check. Very small items have no positive excess; among the
larger items, the bin capacity limits how many can appear together.
Put $\tau=14/13$ and define the excess
\[
 e(w)=\max\{q(w)-\tau w,0\}.
\]
For any bin $B$, this definition gives
\[
 \sum_{e\in B}q(w_e)\le\frac{14}{13}\sum_{e\in B}w_e
                          +\sum_{e\in B}e(w_e)
 \le\frac{14}{13}+\sum_{e\in B}e(w_e).
\]
Our task is therefore to bound the last sum. The selected weight
ranges lower the charge of some medium and heavy items, which is
why their excesses must be checked separately.
The following bounds follow by substituting the three selection
probabilities. On each interval where $q$ is constant, its excess above
$\tau w$ decreases with $w$; it therefore suffices to check the left
endpoints of those intervals.
\[
\begin{array}{c|c@{\qquad}c|c}
\text{heavy range}&\text{bound on }e(w)&
\text{light range}&\text{bound on }e(w)\\\hline
(1/2,7/12]&89/390&(1/3,1/2]&\alpha:=103/1560\\
(7/12,5/8]&31/130&(1/4,1/3]&\beta:=16/351\\
(5/8,1]&19/78&(1/5,1/4]&\gamma:=9/260
\end{array}
\]
For smaller weights, the harmonic bound
$a(w)\le(1+1/t)w$ whenever $w\le1/t$ gives
\begin{equation}\label{eq:bip-small-excess}
 e(w)\le\max\{0,1/t-1/13\}\,w
 \qquad(w\le1/t).
\end{equation}
We use this with $t=5,6,12,13$, giving slopes $8/65$, $7/78$,
$1/156$, and zero, respectively. If the bin has no heavy item, the
preceding bounds give $e(w)\le w/5$ for every item. Its total charge
is therefore at most $14/13+1/5<104/75$.

Now suppose the bin has a heavy item. According to its range in the
first column of the table, the remaining light items have total weight
less than $C=1/2$, $5/12$, or $3/8$, respectively. We claim that for
each of these three capacities their total excess is at most
\begin{equation}\label{eq:bip-light-excess}
 E(C):=\frac8{65}C+\frac1{50}.
\end{equation}
There are three possibilities for the largest remaining weight.

If it exceeds $1/3$, that item contributes at most $\alpha$ and the others
have total weight below $C-1/3$. Using \eqref{eq:bip-small-excess}
gives the second column below.

Otherwise, if it exceeds $1/4$, that item contributes at most $\beta$.
For $C=1/2$, if another item exceeds $1/5$, it contributes at most
$\gamma$ and leaves total weight below $1/20<1/13$, whose excess is zero.
If there is no such item, the remaining excess is at most
$(8/65)(1/4)\le\gamma$. This gives $\beta+\gamma$ in either case.
For the other two capacities, the weight left after the largest item
is below $1/6$ or $1/8$, so the slope $7/78$ applies.
These bounds give the third column.
\[
\begin{array}{c|c|c|c}
 C & \text{largest weight}>1/3 &
 1/4<\text{largest weight}\le1/3 & E(C)\\\hline
 1/2 & \alpha+7/468 & \beta+\gamma & 53/650\\
 5/12 & \alpha+1/1872 & \beta+7/468 & 139/1950\\
 3/8 & \alpha & \beta+7/624 & 43/650
\end{array}
\]
Each entry in the middle columns is at most the last entry of its row.

Finally, if every remaining item has weight at most $1/4$, let $m$
be the number exceeding $1/5$. Since $C\le1/2$, we have $m\le2$.
The other items use the slope $8/65$, so their total excess is at most
\[
 m\gamma+\frac8{65}\left(C-\frac m5\right)
 =\frac8{65}C+\frac{m}{100}\le E(C).
\]
This proves \eqref{eq:bip-light-excess} in every case.

Write $h(C)$ for the heavy-item bound corresponding to capacity $C$.
The three rows satisfy
\[
 h(C)+\frac8{65}C=\frac{113}{390}.
\]
The bin's total $q$-charge is consequently at most
\[
 \frac{14}{13}+h(C)+E(C)
 =\frac{14}{13}+\frac{113}{390}+\frac1{50}
 =\frac{104}{75}.
\]
\end{proof}

\subsubsection{The Resulting Coloring Bound}
\begin{theorem}
	\label{thm:bipartite-improved}
	Every bipartite weighted multigraph $G$ satisfies
	\[
	\cw(G)\le(1.62+o(1))\bp(G).
	\]
	Equivalently, for every $\eta>0$ there is~a~constant $C_\eta$, independent of~graph order, multiplicity, and minimum positive weight, such that $\cw(G)\le(1.62+\eta)\bp(G)+C_\eta$.
	For rational weights and fixed $\eta>0$, such a coloring can be found
	by a randomized algorithm in expected polynomial time.
\end{theorem}
\begin{proof}
Put $n=\bp(G)$ and use the three selections in
\Cref{tab:bip-three-selections}. Their costs sum to
$F=1/20+1/10+1/12=7/30$. Each row satisfies $h+b=1$ and $b<1/2$,
and its selection probability at any weight is at most its cost.
Thus the total probability is at most $7/30<1$.
\Cref{lem:bip-three-choice-charge} shows that the remainder charge
of every unit bin is at most $Q=104/75$.
Applying \Cref{lem:bip-compatible-selections} gives
\[
 \cw(G)\le\left(\frac7{30}+\frac{104}{75}+o(1)\right)n
 =\left(\frac{81}{50}+o(1)\right)n.
\]
The selections partition the original edges, and their additive
rounding errors contribute only $O(1)$ colors. The remaining error
comes from \eqref{eq:bip-harmonic-input}.
\end{proof}

\begin{algorithm}[The bipartite $1.62$ construction]\label{alg:bipartite-improved}
Input: a weighted bipartite multigraph and a desired accuracy.
Choose the packing and coloring accuracies small enough for the
desired final error.
\begin{enumerate}[leftmargin=*,nosep]
\item Find approximate local packings as in
\Cref{sec:local-input-algorithms} and let $n$ be the largest number
of bins used.
\item Assign the three selection labels and the remainder label by
\Cref{lem:bip-balanced-labels}, using the probabilities in
\Cref{tab:bip-three-selections}.
\item Color each selected subgraph by \Cref{lem:bip-two-copy},
using a separate palette for each selection.
\item Color the remainder by \Cref{lem:harmonic-load} with a fresh
palette. Return the union of these colorings.
\end{enumerate}
The label assignment and
the three selected colorings take polynomial time; the harmonic
remainder coloring takes expected polynomial time.
\end{algorithm}

\clearpage
\section{Bin-Packing Bounds for General Multigraphs}
\label{sec:general-bounds}
This appendix gives the proofs for Section~\ref{ov:page-12}.
The first packet bound introduces the fractional bundle costs.
We then prove the endpoint-label extension of rounding and the local
inequalities needed for the $1.999$ bound. The construction and its
implementation follow in the same appendix.

\subsection{A First Bound for General Multigraphs}
\label{sec:general-bp}

The bipartite arguments color packets by K\H{o}nig's theorem; for general multigraphs, dense odd vertex sets create an additional obstruction.
We first show that controlling triangles already gives a bound of $\bigl(\tfrac54h_\infty+o(1)\bigr)\bp(G)$.
The construction also provides the preparation and rounding steps used in the stronger bound of \Cref{thm:general-bp-subtwo}.

We use three related terms in these constructions. A \emph{bundle}
contains all edges on one original pair of endpoints. An \emph{object}
is an edge, or a group of parallel edges, treated as one weighted item
during preprocessing. A \emph{packet} is a collection of objects on one
pair whose total current weight is at most one; its objects receive the
same color. Grouping replaces several edges by an object, whereas
packing collects objects into packets.

\begin{theorem}
\label{thm:general-bp-baseline}
Every weighted multigraph $G$ without loops satisfies
\[
 \chi'_w(G)\le \bigl(\tfrac54h_\infty+o(1)\bigr)\bp(G),
 \qquad \frac54h_\infty=2.113787758446\ldots.
\]
More precisely, for each $\varepsilon>0$ there is $C_\varepsilon$ such
that at most $(\frac54h_\infty+\varepsilon)\bp(G)+C_\varepsilon$
colors suffice. The constant is independent of the number of vertices,
parallel multiplicities, and the smallest edge weight. For rational
weights and explicitly listed edges, each fixed accuracy has a
randomized algorithm with expected polynomial running time.
\end{theorem}

Before giving the proof, we explain the coefficient $5/4$.
We pack parallel edges into packets of weight at most one and properly color the resulting multigraph.
Its maximum degree is at most $(h_\infty+o(1))\bp(G)$.
The first lemma below bounds the number of packets inside every triangle by $(2+o(1))\bp(G)$.
For an odd set of $s\ge5$ vertices, degree counting loses at most the factor $s/(s-1)\le5/4$.
These estimates bound the fractional chromatic index of the packet multigraph, and an ordinary edge-coloring theorem converts that bound into a proper coloring.
Finally, the charge left at each vertex is smaller than the number of missing colors there, allowing us to insert almost all remaining edges into the same palette.

\subsubsection{Packing Two Collections}

For a finite collection $A$ of item types with prescribed counts, write
$F(A)$ for its \emph{fractional count packing cost}. A pattern is a
list of nonnegative integer type counts of total weight at most one.
Assign nonnegative masses to these patterns so that their average
counts equal the prescribed counts; $F(A)$ is the least total mass.
Patterns may repeat a type. Actual packets constructed later will use
distinct available edges. An ordinary packing is one allowed choice,
so $F(A)$ is at most the ordinary bin-packing number.

\begin{lemma}\label{lem:general-two-bags}
For any two collections $A,B$, put $a=F(A)$ and $b=F(B)$. Then
\[
 F(A+B)\ge\min\{a/2+b,\ a+b/2\}.
\]
Consequently, for three collections $A,B,C$,
\[
 F(A)+F(B)+F(C)
 \le2\max\{F(A+B),F(A+C),F(B+C)\}.
\]
\end{lemma}
\begin{proof}
Take a fractional packing of $A+B$ of total mass $T$. In each bin,
keep the $A$ items together as one packet and the $B$ items as another.
An $A$-heavy bin has an $A$ packet heavier than $1/2$; define
$B$-heavy similarly. No bin is both. Write $\ell,m,z$ for the masses
of $A$-heavy, $B$-heavy, and remaining bins, interchanging the names
if necessary so that $\ell\ge m$.

Pair all $B$-heavy mass with the same amount of $A$-heavy mass.
In one such pair let the packet weights be $(a_0,b_0)$ and
$(c_0,d_0)$, where $a_0,d_0>1/2$. If $a_0\le d_0$, then
$a_0+c_0\le1$, so the two $A$ packets fit together. If
$a_0\ge d_0$, the two $B$ packets fit together. Thus at least one
collection costs one bin per unit paired mass. For the other
collection, keep its heavy packet in one bin and put two copies of
its small packet in a bin of mass $1/2$. Its cost is at most $3/2$.

Let $p\in[0,m]$ be the paired mass for which the $A$ packets are put
together. An unpaired $A$-heavy bin costs one for $A$ and $1/2$ for
$B$; a remaining bin costs $1/2$ for each, by doubling its packets.
The separate packing costs are therefore at most
\[
 A'=\ell+\frac m2-\frac p2+\frac z2,
 \qquad B'=\frac\ell2+\frac m2+\frac p2+\frac z2.
\]
Empty packets can be omitted. Since
\[
 A'/2+B'=\ell+3m/4+p/4+3z/4\le T,
\]
and the displayed minimum is increasing in both arguments and
symmetric in them, the first assertion follows. For the second,
rename the three costs so that $a\ge b\ge c$. The first assertion
gives $F(A+B)\ge a/2+b\ge(a+b+c)/2$.
\end{proof}

\subsubsection{Coloring Packets and Inserting the Remaining Edges}

We use two standard results about ordinary multigraph edge coloring.
A fractional proper coloring assigns nonnegative masses to matchings,
covering each edge by total mass at least one; its minimum total mass
is denoted by $\chi'_f(H)$. Edmonds' matching theorem~\cite{Edmonds65} gives
\begin{equation}\label{eq:general-fractional-index}
 \chi'_f(H)=\max\left\{\Delta(H),\
   \max_{\substack{S\subseteq V(H)\\ |S|\ge3\text{ odd}}}
       \frac{2|E_H(S)|}{|S|-1}\right\}.
\end{equation}
Here $E_H(S)$ consists of edges with both endpoints in $S$.
The following result converts this fractional bound into a proper coloring with a lower-order additive error.

\begin{proposition}[Sanders and Steurer~\cite{SS}]
\label{prop:general-ordinary-coloring}
A loopless multigraph $H$ has a proper coloring with at most
\[
 \chi'_f(H)+\sqrt{9\chi'_f(H)/2}
\]
colors. For explicitly listed edges, such a coloring can be found
deterministically in time
\[
 O\!\left(|E(H)|\sqrt{\Delta(H)}\bigl(|V(H)|+\Delta(H)\bigr)\right).
\]
\end{proposition}
This is Theorem~2.23 of~\cite{SS}. Its error depends only on the
fractional chromatic index, not on the number of vertices. For each
fixed $\eta>0$ the bound is also
$(1+\eta)\chi'_f(H)+O_\eta(1)$.

\begin{lemma}
\label{lem:general-charge-finish}
Fix an integer $T\ge4$, put $b=1/T$, and fix $L\in\mathbb N$ and $\rho,C>0$.
Let $H$ have a proper coloring with an even number $K\le Cn$ of
colors. Let $R$ be a weighted multigraph on the same vertices, of
multiplicity at most $L$. Give an edge of $R$ charge
\[
 q_e=\begin{cases}
       a(w_e),&w_e\ge b,\\
       w_e/(1-b),&w_e<b,
     \end{cases}
 \qquad Q_R(v)=\sum_{e\in E_R(v)}q_e.
\]
If $K-d_H(v)\ge Q_R(v)+\rho n$ at every vertex, then for each
$\zeta>0$, after two-color path switches in $H$, $R$ can be colored using colors missing in the resulting coloring of $H$, together
with at most $\zeta n+O_{b,L,\rho,C,\zeta}(1)$ new colors.
For rational inputs and fixed constants, the algorithm has expected
polynomial running time.
\end{lemma}
\begin{proof}
Write $f_v=K-d_H(v)$. The charges dominate weights, their sum at a
vertex is at most $Cn$, and their sum on any pair is at most $L$.
After multiplying its size parameter by a fixed constant,
\Cref{lem:weighted-moment-mix} gives fractions supported on
common missing colors such that $x_{e,c}=O(1/n)$,
$\sum_cx_{e,c}\le1$, and
\[
 \sum_{e\in E_R(v)}q_ex_{e,c}
 \le Q_R(v)/f_v+o(1)\le1-\rho/C+o(1).
\]
The total uncovered weight at each vertex is $o(n)$.

To apply the rounding lemma, we must realize these fractions as inclusion probabilities of feasible sets at each vertex.
Use the class probabilities and conditional distributions from
\Cref{sec:harmonic-inputs}, now allowing $k=1$ as well as
$2\le k\le T$. The same cardinality construction applies to $k=1$:
it chooses at most one edge. The tiny class still uses the load bound
$1-b$. The fixed charge margin gives
\[
 \sum_{k=0}^T p_k
 \le\sum_{e\in E_R(v)}q_ex_{e,c}+(T+1)\max_e x_{e,c}\le1
\]
for large $n$. Use the empty set for the remaining probability.

\Cref{lem:bp-mult-round}, again with its size parameter rescaled
by a fixed constant, leaves weighted degree at most
$\delta n+o_\delta(n)$, for any fixed $\delta>0$. Choose
$\delta<\zeta/8$. \Cref{lem:general-cleanup} then uses at most
$\zeta n$ additional colors for large $n$; bounded $n$ is absorbed
into the additive constant. The rational constructions in
\Cref{sec:mixing-implementation,sec:harmonic-inputs,sec:rounding-implementation} implement the mixing, local distributions,
and bounded-multiplicity rounding in expected polynomial time.
All constants are uniform in graph order and minimum weight.
\end{proof}

\subsubsection{Proof of the General Bound}
\begin{proof}[Proof of \Cref{thm:general-bp-baseline}]
Put $n=\bp(G)$; the empty graph is immediate. Choose a sufficiently large integer $T\ge4$ and put $b=1/T$. On each parallel bundle, group edges lighter than $b$
into objects of weight in $[b,2b)$, leaving at most one object of weight
below $b$. Edges of weight at least $b$ stay as individual objects.
Expanding an object after coloring preserves feasibility.

The objects at each vertex fit in at most
\[
 n'=\left\lceil\frac{n}{1-2b}\right\rceil+1
\]
bins. To see this, keep the original edges of weight at least $b$ in
their original $n$ bins and insert the new objects greedily. Each new
object weighs less than $2b$. Whenever a new bin is opened, all earlier
bins have load greater than $1-2b$, and the total load is at most $n$.

Next fix $0<\sigma<1$. In a bundle containing $m$ objects of weight at
least $b$, sort them by decreasing weight and set $q=\lfloor\sigma m\rfloor$.
If $q=0$, leave them unchanged. Otherwise set aside the first $q$
objects, divide the rest into consecutive groups of at most $q$, and
round each weight up to the largest in its group. There are at most
$D=\lceil2/\sigma\rceil+2$ non-tiny types per bundle: if $q>0$, then
$q>\sigma m/2$, while $q=0$ implies $m<1/\sigma$.
The rounded
object at position $j>q$ is no heavier than the original object at
position $j-q$. This is an injection on each bundle, so the retained
graph $P$ still fits in $n'$ bins at each vertex. The set-aside objects
have total weight at most $\sigma n/b$ at a vertex.

Give retained objects the charges in
\Cref{lem:general-charge-finish}. \Cref{lem:general-bin-charge}
bounds the original incident charge by $h_\infty n$. Grouping tiny
edges adds at most $4bn$ to it: a new object of weight $w<2b$ has
charge at most $w/(1-2b)\le(1+4b)w$, while the original charge was
at least its total weight. The tiny remainder obeys the same bound.
The rank injection cannot increase total charge. Thus
\begin{equation}\label{eq:general-charge-budget}
 Q_P(v)\le(h_\infty+4b)n.
\end{equation}

Solve the fractional count-packing problem on each bundle of $P$.
Choose an optimum solution supported on at most $D+1$ patterns, one
per possible type constraint. A tiny type may have a very small
weight; its pattern set is finite, and there is at most one actual
tiny object in that bundle. From a pattern of mass $\lambda$, form
$\lfloor\lambda\rfloor$ packets using distinct available objects.
The demand equalities ensure that enough objects of each type exist.
The resulting ordinary multigraph is $H$; write $R$ for the objects
left over. Since a pattern contains at most $1/b$ non-tiny objects,
\[
 \text{multiplicity}(R)\le\lceil(D+1)/b\rceil+1.
\]

Every packet used has total charge at least one. Indeed, the dual
packing problem assigns a price $p_t$ to each type, with total price
at most one in every pattern. Prices can be chosen nonnegative:
replace negative prices by zero, using the constraint for the pattern
with those items deleted. The pattern containing
$\lfloor1/w_t\rfloor$ copies of type $t$ gives $p_t\le a(w_t)$.
For a tiny type, $a(w_t)\le w_t/(1-b)$ as well. Equality of primal
and dual optima makes the total price exactly one on every positive
pattern. Therefore
\begin{equation}\label{eq:general-packet-charge}
 d_H(v)+Q_R(v)\le Q_P(v)\le(h_\infty+4b)n.
\end{equation}

For a triangle, let $A,B,C$ be its three retained bundles. Each pair
fits into $n'$ local bins. \Cref{lem:general-two-bags} gives
$F(A)+F(B)+F(C)\le2n'$. The number of actual packets on a bundle
is at most its fractional cost, so the triangle has at most $2n'$
packet edges. For every odd vertex set $S$ of size $s\ge5$,
\[
 \frac{2|E_H(S)|}{s-1}
 \le\frac{s}{s-1}\Delta(H)
 \le\frac54(h_\infty+4b)n.
\]
Equation~\eqref{eq:general-fractional-index} now gives
\[
 \chi'_f(H)\le\max\{2n',\tfrac54(h_\infty+4b)n\}.
\]
Since $\frac54h_\infty>2$, \Cref{prop:general-ordinary-coloring}
and an even-palette
adjustment give a proper coloring with
\[
 K=\left(\frac54h_\infty+O(b)+\eta\right)n+O_{b,\eta}(1)
\]
colors for any fixed $\eta>0$. Add unused colors if needed to ensure
$K\ge\frac54(h_\infty+4b)n$. By
\eqref{eq:general-packet-charge}, the number missing at each vertex
exceeds its residual charge by at least $h_\infty n/4$.
\Cref{lem:general-charge-finish} colors $R$ with these colors
and $\zeta n+O_{b,\sigma,\zeta}(1)$ extra colors, for any fixed
$\zeta>0$. The set-aside objects cost at most
$4\sigma n/b+1$ further colors by \Cref{lem:general-cleanup}.

Choose $b$, then $\eta$, then $\sigma$ and $\zeta$ sufficiently small
in terms of the requested $\varepsilon$. The total is
$(\frac54h_\infty+\varepsilon)n+C_\varepsilon$.
All choices depend only on $\varepsilon$. Small $n$ are absorbed
into $C_\varepsilon$ by the same greedy bound. Finally expand all
packets and objects. This proves the theorem.
\end{proof}

\begin{algorithm}[The packet bound for general multigraphs]
\label{alg:general-bp-baseline}
Input: a weighted multigraph without loops and a desired accuracy.
Choose the parameters in the order specified in the proof.
\begin{enumerate}[leftmargin=*,nosep]
\item Group tiny edges, set aside the prescribed largest objects on
each bundle, and round retained weights by rank.
\item Use the polynomial-size bundle programs described below to form
integer packets.
\item Properly color the packet multigraph by the Sanders--Steurer
algorithm and add unused colors to obtain the required even palette.
\item Apply the missing-color and rounding constructions of
\Cref{lem:general-charge-finish}, then greedily color the remaining
objects with fresh colors.
\item Expand all packets and grouped objects. Check each vertex-color
load and return the resulting coloring.
\end{enumerate}
\end{algorithm}

\subsubsection{Computing the Bundle Packings}
For fixed accuracy, the $D$ non-tiny types admit only constantly many
count vectors $a$ of load $\ell(a)\le1$. If a tiny type has weight
$t$, set $k_a=\lfloor(1-\ell(a))/t\rfloor$. It suffices to keep the
two patterns $(a,0)$ and $(a,k_a)$: every $(a,k)$ with
$0\le k\le k_a$ is a convex combination of these two patterns with
the same total mass. The reduced linear program therefore has the
same optimum. Its size is constant for fixed accuracy, and even a
large $k_a$ has polynomial binary encoding length. Flooring a basic
solution forms only packets that can be filled with the listed objects.

The algorithm uses the approximate local packings from \Cref{sec:local-input-algorithms}, so it does not need to compute $\bp(G)$ exactly.
Their error is absorbed into the requested accuracy.
Rational upper approximations to $h_\infty$ follow from its recurrence;
the tail after $1/r_j$ is at most $2/r_{j+1}$. Packet formation and
ordinary coloring are polynomial; the missing-color construction has
expected polynomial running time. Together these give the claimed
fixed-accuracy algorithm for the packet bound.

\subsection{Rounding While Keeping Endpoint Labels}\label{sec:label-rounding}
\begin{lemma}[Endpoint labels in bounded-multiplicity rounding]
\label{lem:bp-labelled-round}
Fix integers $L,h\ge1$, constants $C_0,C_1\ge1$, and
$\delta,\eta>0$.
In Lemma~\ref{lem:bp-mult-round}, give each incidence $(e,z)$ a fixed
label $\ell_z(e)$ from a set $\mathcal L_z$ of at most $h$ labels.
Suppose that, for every vertex $z$ and color $c$, the vector
$(x_{e,c}:e\ni z)$ is the mean of a distribution on feasible sets,
each of which contains edges of at most one label at $z$.
Let $\pi_{z,c,a}$ be the probability that this distribution chooses
label $a$, assigning the empty set a null label if necessary, and put
$\Pi_{z,a}=\sum_c\pi_{z,c,a}$.
Then the partial coloring in Lemma~\ref{lem:bp-mult-round} may be chosen
so that, in addition to its stated weight and support conclusions,
each color used at $z$ contains edges of at most one label, and
\[
 \#\{c:\text{color $c$ uses label $a$ at $z$}\}
 \le \Pi_{z,a}+\eta n+O_{L,C_0,C_1,h,\delta,\eta}(n^{3/4})
 \qquad (a\in\mathcal L_z).
\]
The threshold and error constant do not depend on the number of
vertices. The conclusion also holds when the edge deletion in the
proof of Lemma~\ref{lem:bp-mult-round} is performed first and the
local distributions are then restricted to the surviving edges.
If the local distributions are rational and supplied with polynomial
support, the partial coloring can be found in expected polynomial
time for fixed parameters.
\end{lemma}

\begin{proof}
Use the proof of Lemma~\ref{lem:bp-mult-round}, with the following
decorations and one additional local test. At $(z,c)$, first draw a
label $a$ with probability $\pi_{z,c,a}$ and condition the feasible-set
distribution on that label. Apply Lemma~\ref{lem:profiles} to this
conditional distribution, with a small compression parameter $\alpha>0$ chosen for the requested output error $\delta$. If $\pi_{z,c,a}>0$, multiply all its profile
probabilities and assigned fractions by $\pi_{z,c,a}$, and decorate
every resulting profile, including its two repair profiles and empty
profile, with $a$. Omit zero-probability labels. Add a null-labeled
empty profile if needed. Thus the total profile probability is one,
the unconditional assigned fraction of $e$ is
$(1-\alpha)x_{e,c}$, and an assigned fraction for $e$ appears only in
a profile decorated with $\ell_z(e)$. There are at most
$hM(\beta,\alpha)+1$ decorated profiles, independently of $n$.
The conditional fractions may be large when $\pi_{z,c,a}$ is small;
this causes no change to the proof because every \emph{unconditional}
assigned fraction is at most $(1-\alpha)x_{e,c}=O(1/n)$.

Smooth the decorated-profile probabilities exactly as in the proof of
Theorem~\ref{thm:round}: with $M_*$ the common bounded decorated
alphabet, put $\lambda'=(1-\theta)\lambda+\theta/M_*$. This preserves
the option support: smoothing changes selection probabilities, not the
assigned fractions. For the independent original list at $z$, let
$C_{z,a}$ count colors whose selected profile is decorated with $a$.
Then
\[
 \mathbb E C_{z,a}
 \le (1-\theta)\Pi_{z,a}+\theta|\mathcal C|
 \le \Pi_{z,a}+\theta C_1n.
\]
Put $t=\lceil n^{3/4}\rceil$. Independence across colors gives
\[
 \operatorname{Var}(C_{z,a})\le |\mathcal C|/4,
 \qquad
 \Pr\{C_{z,a}>\mathbb E C_{z,a}+t\}
 \le\frac{|\mathcal C|}{4t^2},
\]
where the second inequality is Chebyshev's inequality.
In the group-by-group choice of vertex lists, retain at each vertex
its usual list with progress value $Q_z^{(k)}\le\lceil n^{3/4}\rceil$
only if it also passes these at most $h$ label-count tests. The tests
depend on that vertex's list alone. Under the original list law the
discarded probability at a vertex is at most
\begin{equation}\label{eq:label-rejection-bound}
 p'_k=p_k+\frac{h|\mathcal C|}{4t^2}
      =p_k+O(n^{-1/2})=o(1).
\end{equation}
Consequently every
nonnegative row expectation in the Harris--Srinivasan step increases
by at most the same $1/(1-o(1))$ factor. The load and progress
calculations, the $O(n^{3/4})$ error, and both large- and small-edge
rounding go through unchanged. This is a per-vertex conditioning
argument; it does not use a union bound over all vertices.

Every colored edge uses an option with positive assigned fraction at
both ends. Its color therefore selects a decorated profile with the
edge's fixed incidence label. All edges of one color at a vertex have
the same label. The later edge normalization and both whole-edge
roundings only delete options, so they cannot add a used label or
increase a label count. The retained-list test gives
$C_{z,a}\le\Pi_{z,a}+\theta C_1n+t$.
Choose $\theta$ small enough that $\theta C_1<\eta$; the other fixed
accuracy parameters can then be chosen as in the original rounding
proof. This proves the assertion. Restricting a local random feasible
set after group deletion preserves its label and the same argument
applies. The implementation is given in \Cref{sec:label-algorithm}.
\end{proof}

\subsection{A Bound Below Two for General Multigraphs}
\label{sec:general-subtwo}

The construction in \Cref{sec:general-bp} treats triangles separately and bounds larger odd sets by degree counting.
To improve the coefficient below two, we must also account for dense sets of five vertices.
We select disjoint dense blocks, contract them, and restore their internal edges after coloring the quotient.

\begin{theorem}\label{thm:general-bp-subtwo}
Every weighted loopless multigraph $G$ satisfies
\[
 \chi'_w(G)\le(1.999+o(1))\bp(G).
\]
Equivalently, for every $\varepsilon>0$ there is a constant $C_\varepsilon$, independent of graph order, multiplicity, and minimum positive weight, such that
\[
 \chi'_w(G)\le(1.999+\varepsilon)\bp(G)+C_\varepsilon.
\]
For rational weights and explicitly listed edges, such a coloring can
be found by a randomized algorithm with expected polynomial running
time for every fixed $\varepsilon>0$.
\end{theorem}

The main difficulty is restoring a contracted set without conflicting with its external colors.
We record the original endpoint used by each external color and require the internal pattern in that color to avoid this endpoint.
A local price inequality shows that the resulting restricted palettes still suffice.
The proof first establishes the local estimates and the contraction rule, then combines them with the packet and rounding constructions of \Cref{sec:general-bp}.

We use the following constants:
\[
 A=\frac{57}{40},\qquad R=\frac7{10},\qquad
 M=\frac{17}{10},\qquad c=\frac{39}{20},\qquad
 \kappa=\frac{1999}{1000}.
\]
Here $A$ bounds the harmonic charge of a bin containing only light items, $R$ bounds the light items sharing a bin with a heavy item, and $M$ bounds an arbitrary bin.
The constant $c$ is the local five-vertex price bound, and $\kappa$ is the coefficient in the theorem.
We will use the strict inequalities
\begin{equation}\label{eq:subtwo-margins}
4\kappa-(2A+3M)=\frac{23}{500}>0,
\qquad \kappa-\frac{7M}{6}=\frac{47}{3000}>0.
\end{equation}

\subsubsection{Three Bin Estimates}
Write $\tau(w)=a(w)=1/\lfloor 1/w\rfloor$ for the harmonic charge. If $w\le1/j$, then
$\tau(w)\le(1+1/j)w$.  The following elementary estimates replace an
estimate involving the full harmonic constant.

\begin{lemma}\label{lem:subtwo-charge}
A unit bin containing only items of weight at most $1/2$ has total
$\tau$-charge at most $A$.  Items of total weight below $1/2$ have
$\tau$-charge at most $25/36<R$. Every unit bin has
$\tau$-charge at most $61/36<M$.
\end{lemma}
\begin{proof}
For a light-only bin, at most two items exceed $1/3$.  With at most
one, the charge is at most $25/18$.  With two, their charge is one and
their complement weighs below $1/3$.  If the complement has no item
above $1/4$, its charge is at most $5/12$; otherwise it has one such
item of charge $1/3$, and its remaining load is below $1/12$ with
charge at most $13/144$.  This gives $1+1/3+13/144=205/144<57/40=A$.

For load below $1/2$, if no item exceeds $1/3$, the charge is at
most $2/3$. Otherwise the unique item above $1/3$ contributes $1/2$,
and the remainder has load below $1/6$. Since
$\tau(w)\le(7/6)w$ for $w\le1/6$, its charge is at most $7/36$.
Thus the total charge is at most $25/36$. A whole bin is either
light-only or contains one heavy item followed by such a remainder.
\end{proof}

The margin $M-61/36=1/180$ also covers tiny external objects charged
at $w/(1-b)$ rather than $\tau(w)$: the increase in one local bin is
at most $b/(1-b)<1/180$ when $b\le1/200$.

\begin{lemma}\label{lem:subtwo-preprocess}
Let $G$ be a weighted multigraph. Fix an integer $T_0\ge4$, put $b=1/T_0$, and fix $0<\sigma<1$. With $n=\bp(G)$, there is a
retained object graph $P$ and a set-aside collection such that every
retained star fits in
\[
N=\left\lceil\frac{n}{1-2b}\right\rceil+1
\]
bins, every original pair has at most
$D=\lceil2/\sigma\rceil+2$ non-tiny types and one tiny object,
and the set-asides have weighted degree at most $\sigma n/b$.
On each original pair the retained objects yield whole packets whose
number is at most the pair's fractional count-packing cost.  Objects
not used by packets have multiplicity at most
$\lceil(D+1)/b\rceil+1$ per original pair.  Every used packet has
harmonic charge at least one.  Expanding objects and packets back to
their original edges preserves feasibility.
\end{lemma}
\begin{proof}
Apply the grouping, rank rounding, and packet construction from the
proof of \Cref{thm:general-bp-baseline}, with $T=T_0$ and $n'=N$.
That construction gives the displayed local-bin bound, the bound $D$
on the number of types, the set-aside weight, and the residual
multiplicity. Flooring an optimal fractional packing gives no more
whole packets than its cost, and expansion preserves feasibility.

For the harmonic-charge assertion, choose optimal nonnegative dual
prices. The price of a type of weight $w$ is at most $\tau(w)$,
because $\lfloor1/w\rfloor$ copies fit in one bin. Every pattern
with positive mass has total dual price one, so every used packet
has harmonic charge at least one.
\end{proof}

\subsubsection{The Five-Vertex Price Inequality}\label{app:subtwo-five-local}
We work with the fractional count-packing model of
\Cref{sec:general-bp}.  A type specifies one original pair and
one positive item weight; patterns may repeat types.  On five vertices
let $x$ be a demand whose star at each vertex is the mean of feasible
integer star patterns.  Write $t_e$ for the separate fractional bin
cost of pair $e$, $T=\sum_e t_e$, and $C=p\cdot x$ for a nonnegative
type price vector $p$, normalized so that every feasible whole pattern
has price at most one.  Let $\beta_e$ be the largest price of a packet
using pair $e$ alone.

These prices are a way to test whether the whole demand fits in a
given amount of palette.  The normalization says that one color can
carry price at most one.  The quantity $t_e$ measures how much palette
pair $e$ needs when packed separately, whereas $\beta_e$ measures the
largest price that one packet on that pair can carry.  They play
different roles: large $T$ will force many heavy items, and the
$\beta_e$ will help us bound their total price.

\begin{lemma}\label{lem:subtwo-five-local}
If $T\ge2\kappa$, then $C\le c=39/20$.
\end{lemma}
The proof has two steps.
First, density forces enough heavy items to occur in the local star distributions.
We then distribute the light-item prices between their endpoints and bound the resulting local prices by three explicit model configurations.

\begin{proof}[Proof of \Cref{lem:subtwo-five-local}]
Write $H_e$ for the mean number of heavy items, of weight above
$1/2$, on pair $e$, and $d_H(v)=\sum_{e\ni v}H_e$.  A star contains at
most one heavy item, so $d_H(v)\le1$.  Packing type $t$ alone in
$\lfloor1/w_t\rfloor$-copy bins gives
$t_e\le\sum_{t\text{ on }e}\tau(w_t)x_t$.  Averaging this charge
over the local bin laws and using \Cref{lem:subtwo-charge} gives,
at vertex $v$, charge at most $A+(M-A)d_H(v)$.  Summing over the
five vertices counts each pair charge twice, so
\[
2T\le5A+2(M-A)\sum_e H_e.
\]
Indeed, a heavy-free star contributes at most $A$, and a star
containing a heavy item contributes at most $M$.  Thus the preceding
inequality is just the average of these two bounds.  If
$\sum_eH_e<3/2$, its right-hand side would be less than
$5A+3(M-A)=2A+3M<4\kappa$, whereas $2T\ge4\kappa$.
Consequently $\sum_eH_e\ge3/2$.  This lower bound on heavy mass is
the saving that we will use throughout the price argument.

\paragraph*{Reducing the pair prices}
Packets on disjoint pairs fit together; hence
$\beta_e+\beta_f\le1$ when $e$ and $f$ are disjoint.  Put
$a=\max_e\beta_e$.  If $a\le6/11$, a light-only star has price at most
$Aa$.  Pair endpoint states containing the same heavy type.
For each heavy type $t$ on $uv$, the states containing $t$ have total
mass $x_t$ at both endpoints, since both local distributions realize
the same demand and a feasible star contains at most one heavy item.
Couple the conditional endpoint distributions for this type and give
the resulting paired states total mass $x_t$; omit types with $x_t=0$.
Keep the
entire packet of smaller weight on its pair and both off-pair parts; at any
third vertex the two off-pair loads are each below $1/2$, so the
merged pattern is feasible.  The paired price is at most $1+a$.
After pairing, the heavy states use at most $(1+a)/2$ price per
endpoint state, averaged over the two endpoints, and a light-only
state uses at most $Aa$.  There is
one unit of state mass at each of the five vertices, and their total
price is $2C$.
Consequently
\[
 C\le\frac52\max\{Aa,(1+a)/2\}
 \le\frac{171}{88}<c.
\]

Assume $a\ge6/11$ and name a maximum pair $01$.  Order the other
three $0$-spokes so that $\beta_{02}\ge\beta_{03}\ge\beta_{04}$,
and put $s=1-a$, $x=\max(s,\beta_{02})$, and
$y=\max(s,\beta_{03})$.  The disjoint-packet inequalities bound the
ten pair prices by
\[
\begin{array}{c|cccccc}
 e&01&02&03,04&12&13,14&23,24,34\\ \hline
 m_e&a&x&y&1-y&1-x&s.
\end{array}
\]
For example, pair $12$ is disjoint from $03$, so
$\beta_{12}\le\min\{a,1-\beta_{03}\}=1-y$.
Pairs $13$ and $14$ are disjoint from $02$, giving
$\beta_{13},\beta_{14}\le\min\{a,1-\beta_{02}\}=1-x$.
The remaining pairs $23,24,34$ are disjoint from $01$ and have
price at most $1-a=s$.
This vector is the convex combination of three simple arrays:
$S_0$ prices all four $0$-spokes at $a$, $S_1$ prices all four
$1$-spokes at $a$, and $Q$ prices triangle $012$ at $a$; every
other pair has price $s$.  The coefficients of $S_0,S_1,Q$ are
\begin{equation}\label{eq:subtwo-model-weights}
 \lambda_0=\frac{y-s}{2a-1},\qquad
 \lambda_1=\frac{a-x}{2a-1},\qquad
 \lambda_Q=\frac{x-y}{2a-1}.
\end{equation}
They are nonnegative and sum to one, because $s\le y\le x\le a$
and $a-s=2a-1>0$.  To check the decomposition, pair $01$ has price
$a$ in all three arrays, pair $02$ has price $a$ in $S_0,Q$ and
$s$ in $S_1$, and pairs $03,04$ have price $a$ only in $S_0$.
Their averages are therefore $a,x,y$, respectively.  The other
entries follow in the same way.  We have reduced a vector with ten
coordinates to two star-shaped arrays and one triangle-shaped array.
The individual arrays need not bound the original prices; it is
their displayed average that does so.

\paragraph*{A common bound for the three models}
We next derive an inequality that applies to all three arrays.  For upper pair
prices $m_e\ge\beta_e$, choose nonnegative endpoint shares with
$z_{ue}+z_{ve}=2m_e$.  Give a heavy item its full price at each end
and a light item on $e$ the multiplier $z_{ve}/m_e$ at end $v$;
zero-price terms are omitted.  This only changes where light-item
price is charged: its two multipliers sum to two.  Heavy-item price
is also counted twice, so the sum of the five local means is $2C$.
We are free to choose the shares so that the local bounds become
easier to satisfy.

There are two kinds of local state to bound: a light-only star at
one vertex, and a pair of endpoint states containing the same heavy
item.  Choose $h_0=h_1=h$ and $h_2=h_3=h_4=0$, where $h=0$ for
$a<3/4$ and $h=1/2$ for $a\ge3/4$, and define
\begin{align*}
 B_v&=h_v+\operatorname{top}_1(z_{ve}-h_vm_e)_+
          +(A-1)\operatorname{top}_2(z_{ve}-h_vm_e)_+,\\
 D_v&=\max_{e\ni v}(z_{ve}-m_e)_+,
 &J_v&=\max_{e\ni v}z_{ve},\\
 q^G_{uv}&=1+m_{uv}+R(D_u+D_v),
 &q^S_{uv}&=2m_{uv}+R(J_u+J_v).
\end{align*}
The $\operatorname{top}_i$ are the $i$th largest positive entries
among the pairs incident to $v$, padded with zero.  The quantity
$B_v$ will bound the price of a light-only state.  The quantity
$D_v$ records the largest positive difference between an endpoint
share and its pair price; $J_v$ records the largest endpoint share.
They give two bounds for a heavy pair: $q^G$ combines the two stars
into a feasible pattern, while $q^S$ bounds their heavy items and
light remainders separately.

Fix a light-only star pattern at $v$, and let
$p_e$ be the total price of its items on pair $e$. Put $r_e=p_e/m_e$
when $m_e>0$ and $r_e=0$ otherwise; if $m_e=0$, both $p_e$ and
$z_{ve}$ are zero.  Thus $r_e$ is the price on one pair, measured in
units of its allowed packet price $m_e$; it is not its total weight.
Each pair packet has price at most $\beta_e\le m_e$.
Also, each type $t$ on $e$ has price at most $m_e\tau(w_t)$, since
$\lfloor1/w_t\rfloor$ copies fit in one packet. Hence
\[
 0\le r_e\le1,\qquad
 \sum_{e\ni v}r_e\le A,\qquad
 \sum_{e\ni v}m_er_e=\sum_{e\ni v}p_e\le1.
\]
The first inequality says that the items on one pair form a packet.
The second uses the harmonic-charge bound for a light-only star.
The third uses the fact that the entire star is itself a feasible
pattern, of price at most one.
Since $h_v\ge0$, its local price satisfies
\[
 \sum_{e\ni v}z_{ve}r_e
 \le h_v+\sum_{e\ni v}(z_{ve}-h_vm_e)_+r_e
 \le B_v.
\]
For the final inequality, the sum is maximized by putting one unit
on the largest positive coefficient and at most $A-1$ on the second
largest, since $1\le A\le2$.  The parameter $h_v$ allows us to use
the third constraint as well: we first pay $h_v$ for
$\sum m_er_e\le1$, and then apply the first two constraints to the
remaining coefficients.  The choice $h=1/2$ will be useful only in
the final range of $a$.

For a heavy type on $uv$, couple
its two endpoint states by type.  Let $\pi_u,\pi_v$ be the prices of their
entire central packets and $O_u,O_v$ those of their off-central
parts. Keeping the entire central packet of smaller weight and both
off-central parts gives a feasible pattern.
At $u$ and $v$, replacing the central packet by the lighter one
cannot increase either load.  At any other vertex, each of the two
off-central contributions is below $1/2$, because its original star
contained a heavy item.  This explains the feasibility even when
the two central packets have different light items.
Its central price is at least $\min(\pi_u,\pi_v)$, hence
$\min(\pi_u,\pi_v)+O_u+O_v\le1$.  Since each central packet costs at
most $m_{uv}$, the sum of the two original star prices is at most
$1+m_{uv}$.  Replacing unit light coefficients by endpoint shares
adds at most $R(D_u+D_v)$: each light remainder has harmonic charge
at most $R$, and on any incident pair the coefficient increase is
at most $(z_{ve}-m_e)_+$.  This proves $q^G_{uv}$ even when the two
conditional central packets differ.  Alternatively each heavy item
has price at most $m_{uv}$, and each endpoint's entire light
remainder contributes at most $RJ_v$, proving $q^S_{uv}$.
Thus, for any choice of either bound on each pair,
\begin{equation}\label{eq:subtwo-source-functional}
 2C\le F(m,z;H):=
 \sum_v(1-d_H(v))B_v+\sum_{uv}H_{uv}q_{uv}.
\end{equation}
Here $1-d_H(v)$ is exactly the probability that the star at $v$
has no heavy item.  The remaining states are paired by their heavy
type and charged through $H_{uv}q_{uv}$.  Thus every local state is
accounted for once in $F$.
A budget $P_v\ge B_v$ covers a light-only state at $v$.
When $q_{uv}\le P_u+P_v-\nu$, the two endpoint budgets cover the
paired heavy states with $\nu$ to spare. Since the total heavy mass
is at least $3/2$, these savings add up to at least $3\nu/2$.
Precisely, if $P_v\ge B_v$ and $\nu\ge0$ satisfy
$P_u+P_v-\nu\ge q_{uv}$, then
\begin{equation}\label{eq:subtwo-pot-criterion}
 F(m,z;H)\le\sum_vP_v-\frac32\nu.
\end{equation}
This is weak linear programming duality in a particularly simple
form.  The local state masses satisfy
$(1-d_H(v))+\sum_{e\ni v}H_e=1$ at every vertex.  The number $P_v$
is a budget for this unit of state mass: $P_v\ge B_v$ covers a
light-only state, and $P_u+P_v-\nu\ge q_{uv}$ covers a paired heavy
state.  The multiplier $\nu$ uses the global constraint
$\sum_eH_e\ge3/2$.  No optimization algorithm is needed: any
displayed budgets satisfying these inequalities prove an upper
bound.  Directly substituting the inequalities makes the cancellation
visible:
\[
\begin{split}
 F(m,z;H)
 &\le\sum_v(1-d_H(v))P_v
       +\sum_{uv}H_{uv}(P_u+P_v-\nu)\\
 &=\sum_vP_v-\nu\sum_eH_e
 \le\sum_vP_v-\frac32\nu.
\end{split}
\]
The rest of the proof supplies these budgets for the three model
arrays.  Each check has only two forms: a vertex inequality
$P_v\ge B_v$, and a pair inequality $P_u+P_v-\nu\ge q_{uv}$.

For the star model $S_0$, split each high spoke into root share
$2a-1/2$ and leaf share $1/2$; split a low pair equally at $s$.
The shares for $S_1$ are symmetric.  In $Q$ split every pair equally.
Average these shares with \eqref{eq:subtwo-model-weights}.  With fixed
$h$, both $B_v$ and either $q^G,q^S$ are convex in $(m,z)$:
the top-two expression is a nonnegative combination of the largest
entry and the sum of the two largest entries, since $1\le A\le2$.
All model bounds below use the same actual $H$.  We only have to
account explicitly for one pair on which different models choose
different heavy-state bounds.

To see why convexity is enough, the coefficients $1-d_H(v)$ and
$H_{uv}$ in $F$ are nonnegative.  If the same choice of $q^G$ or
$q^S$ is made on a pair in every model, the value at the averaged
array is at most the average of the three model values.  We can
therefore average their numerical upper bounds without solving a
new problem for the original ten pair prices.  In the first range,
pair $34$ uses different choices in the star and triangle models;
we will bound that one additional error explicitly.

The choices of heavy-state bounds are summarized below.
Every unlisted pair uses $q^G$. Here $R=7/10$; the budget
formulas and their checks are given in the corresponding cases.
\begin{center}
\small
\renewcommand{\arraystretch}{1.2}
\begin{tabular}{@{}ccccc@{}}
\toprule
Range of $a$ & $h$ &
$q^S$ pairs in $S_0,S_1$ &
$q^S$ pairs in $Q$ &
$\nu$ in $S_0,S_1$\\
\midrule
$[6/11,R]$ & $0$ & none & $34$ & $0$\\
$[R,3/4)$ & $0$ & $23,24,34$ & $23,24,34$ & $0$\\
$[3/4,1]$ & $1/2$ & $23,24,34$ & $23,24,34$ & positive\\
\bottomrule
\end{tabular}
\end{center}
The triangle model uses the same budgets in all three ranges:
\[
 P_0=P_1=P_2=Aa,\qquad
 P_3=P_4=Ms+\nu/2,\qquad
 \nu=(2A-1)a-1.
\]
The star budgets for $S_1$ are obtained from those for $S_0$
by exchanging vertices $0$ and $1$. At $a=R$, either of the
first two constructions applies. Only the first range needs
the separate averaging correction on pair $34$.

\paragraph*{The range $6/11\le a\le R$}
First let $6/11\le a\le R$.  In both stars use $q^G$ everywhere.
Put
\[
 B_0=A(2a-1/2),\quad B=1/2+(A-1)s,
 \quad q_H=1+a+R(a-1/2),\quad q_L=1+s.
\]
Potential $q_L/2$ at every leaf, $\max(B_0,q_H-q_L/2)$ at
the root, and $\nu=0$ satisfy all rows.  Here $B_0$ is the root's
light-state bound and $B$ is the common leaf bound.  Every pair
through the root has heavy-state bound $q_H$, and every pair of
leaves has bound $q_L$.  The leaf budgets cover their light states
because
\[
 q_L/2-B=(3/2-A)s=3s/40\ge0.
\]
Two leaf budgets sum to $q_L$, while the root budget was chosen to
cover both $B_0$ and the remaining amount $q_H-q_L/2$ on a high
pair.  Thus these few inequalities check all five vertices and all
ten pairs.  Their objective is
\begin{equation}\label{eq:subtwo-low-star}
 S(a)=\max\left\{\frac{17a}{20}+\frac{263}{80},
                  \frac a5+\frac{73}{20}\right\}
 \le\frac{1553}{400}<\frac{39}{10}.
\end{equation}
For the triangle model $Q$, use $q^S$ only on pair $34$ and $q^G$
elsewhere.  Give the three triangle vertices the same budget $P_H$
and the two outside vertices the same budget $P_O$:
\[
 P_H=Aa,\qquad P_O=Ms+\nu/2,\qquad \nu=(2A-1)a-1.
\]
The choice of $\nu$ makes the high-pair constraint an equality;
the choice of $P_O$ does the same for the outside pair.  There are
only three different pair constraints:
\[
\begin{array}{c|c|c}
 \text{pairs}&q_{uv}&P_u+P_v-\nu-q_{uv}\\ \hline
 01,02,12&1+a&0\\
 03,04,13,14,23,24&1+s&s/5\\
 34&2Ms&0
\end{array}
\]
The light bounds are $Aa$ inside the triangle and $As$ outside it.
Thus $P_H$ covers the first, and
$P_O-As=(26a-9)/40>0$ covers the second.  Also $\nu\ge0$ for
$a\ge6/11$.  Substituting these budgets into
\eqref{eq:subtwo-pot-criterion} gives the objective
\begin{equation}\label{eq:subtwo-low-triangle}
 Q(a)=\frac{39}{10}-\frac a{20}<\frac{39}{10}.
\end{equation}
We can now combine the star and triangle bounds.  Every pair uses
the same heavy-state bound in all three models except $34$:
the stars use $q^G$, while the triangle uses $q^S$.  The minimum of
these two convex functions need not be convex, so this pair requires
a separate calculation.  Write $q=\lambda_Q$,
$\mu=\max(\lambda_0,\lambda_1)$, and
$\sigma=\min(\lambda_0,\lambda_1)$.  At the averaged shares,
$D_3=D_4=0$ and $J_3=J_4=s+(a-1/2)\mu$.  If
$E=(1-q)(1+s)+q(2Ms)$, $b_0=(12a-7)/5$, and $d=R-a\ge0$,
direct subtraction gives
\[
q^G_{34}-E=qb_0,
\qquad q^S_{34}-E=(1-q)d-R(2a-1)\sigma.
\]
The quantity $E$ is exactly the average bound already paid for on
pair $34$.  Either of the two displayed differences is therefore
an upper bound on the additional amount needed at the averaged
shares; we may use whichever is smaller.
If $b_0\le0$, there is no loss and the bounds $S(a),Q(a)\le39/10$ suffice.
Otherwise $a>7/12$,
and the loss is at most $(1-q)d$, since $H_{34}\le1$.
Moreover,
\[
 S(a)+d=
 \max\left\{\frac{319}{80}-\frac{3a}{20},
              \frac{87}{20}-\frac{4a}{5}\right\}
 \le\frac{39}{10}.
\]
Indeed both expressions decrease with $a$; at $a=7/12$
they equal $39/10$ and $233/60$, respectively.
Since $Q(a)\le39/10$, the combined bound is
\[
 2C\le(1-q)S(a)+qQ(a)+(1-q)d\le\frac{39}{10}.
\]

\paragraph*{The range $R\le a<3/4$}
In this interval we can avoid the preceding exceptional-pair
calculation.  Use $q^S$ on $23,24,34$ in all three models
and $q^G$ elsewhere.  Since $2s+R\le1+s$, the preceding triangle
certificate remains valid.  In $S_0$ set
\[
 q_M=2s+R,\quad P_B=\max(B,q_M/2),\quad
 P_A=q_L-P_B,\quad P_0=\max(B_0,q_H-P_B).
\]
Take $P_1=P_A$, $P_2=P_3=P_4=P_B$, and $\nu=0$.
The root is vertex $0$; leaf $1$ has budget $P_A$, and the other
three leaves have budget $P_B$.  The latter need to cover both
their light bound $B$ and half of the bound $q_M$ on pairs among
$2,3,4$.  This explains the maximum defining $P_B$.
Since $P_B\le q_L/2$, we have $P_A\ge P_B\ge B$.
The pairs $1j$, for $j\in\{2,3,4\}$, are covered exactly by
$P_A+P_B=q_L$.  Finally the root's maximum covers its light state
and every root pair, since the smallest leaf budget is $P_B$.
The inequalities used here are $2B\le q_L$ and $q_M\le q_L$;
the latter is just $a\ge R$.  The objective
\[
\max(B_0+q_L+2P_B,\ q_H+q_L+P_B)
\]
is convex in $a$, with endpoint values $1553/400$ at $R$ and
$311/80$ at $3/4$, both below $39/10$.  The $S_1$ certificate is
symmetric. Evaluating these formulas at $a=3/4$ bounds the preceding
interval; the next case supplies the certificate at $a=3/4$ itself.

\paragraph*{The range $3/4\le a\le1$}
Finally let $3/4\le a\le1$ and $s=1-a$.  Use $q^S$ on
$23,24,34$ in all three models.  In $S_0$ the light bounds are
\[
B_0=\frac12+A-\frac{3A}{2}s,\quad
B_A=\frac12+\frac A2s,\quad
B_B=\frac12+(A-1)s.
\]
Set
\[
\nu=\min\{B_0+B_B-q_H,\ 2B_B-(2s+R)\},\quad
(P_0,P_1,P_2,P_3,P_4)=(B_0,B_A+\nu,B_B,B_B,B_B).
\]
The change from the preceding interval is that $h_0=h_1=1/2$.
It uses the fact that a feasible light star has total price at most
one to improve the bound at the root.  Using the same $h_v$ in
$S_0$ and $S_1$ lets us average the two models; the other three
leaf bounds are unchanged.  We can now take a positive
$\nu$, obtaining a saving from each heavy pair.

The two terms defining $\nu$ are $3/40-s/80$ and
$3/10-23s/20$, both nonnegative for $0\le s\le1/4$.
Their definitions check two of the pair classes directly:
\[
\begin{array}{c|c|l}
 \text{pairs}&q_{uv}&\text{required inequality}\\ \hline
 02,03,04&q_H&B_0+B_B-\nu\ge q_H\\
 23,24,34&2s+R&2B_B-\nu\ge2s+R\\
 12,13,14&q_L&B_A+B_B\ge q_L
\end{array}
\]
The last row follows from $B_A+B_B-q_L=11s/80\ge0$.
For the remaining pair $01$, the left-hand side is $B_0+B_A$,
which is at least $B_0+B_B$ because
$B_A-B_B=(1-A/2)s\ge0$.  All vertex budgets cover their light
bounds as well, since $\nu\ge0$.  The objective is
\[
\sum_vP_v-\frac32\nu
=\max\left\{\frac{311}{80}-\frac{23s}{160},
              \frac{151}{40}+\frac{17s}{40}\right\}
\le\frac{311}{80}<\frac{39}{10}.
\]
The $S_1$ argument swaps $0$ and $1$.  For $Q$, the preceding
triangle budgets still work.  At vertices $0,1$, the new light
bound is $1/2+Aa/2\le Aa$, since $Aa\ge1$ in this interval;
the other vertex bounds are unchanged.  Thus choose
$P_H=Aa$ on $0,1,2$, $P_O=Ms+\nu/2$ on $3,4$, and
$\nu=(2A-1)a-1$. The cross-pair slack is at least $s/5$;
the triangle and outer rows hold with equality.  Its objective is
\[
 2Aa+2Ms+\frac{1+a}{2}
 =\frac{39}{10}-\frac a{20}\le\frac{39}{10}.
\]
The common-model convexity and \eqref{eq:subtwo-pot-criterion} now give
$2C\le39/10$, hence $C\le39/20$, in every range, as required.
\end{proof}

\subsubsection{Contracting Dense Sets}
Only triangles and five-vertex sets can violate the target odd-set bound: larger odd sets already satisfy it by degree counting.
The following lemma shows that all necessary contractions can be chosen at once.

\begin{lemma}\label{lem:subtwo-blocks}
Let $t_{uv}\in[0,1]$ be pair costs with incident sum at most $M$ at
each vertex and cost at most two in every triangle.  Call a five-set
dense if its internal cost exceeds $2\kappa$, and a triangle dense
if its internal cost exceeds $\kappa$.  Select every dense five-set,
then every dense triangle disjoint from them, and contract these sets
simultaneously, discarding internal pairs.  The selected sets are
disjoint.  In the quotient every vertex has degree cost at most $M$,
and every odd set $S$ has internal cost at most
$\kappa(|S|-1)/2$.
\end{lemma}
\begin{proof}
For a vertex set $U$, write $t(U)$ for the sum of costs on pairs
inside $U$, and $\partial t(U)$ for the sum on pairs with exactly
one endpoint in $U$. Summing the degree bounds over $U$ gives
\[
 2t(U)+\partial t(U)\le |U|M.
\]
Thus a set with large internal cost has little cost leaving it.
We first show that the selected sets are disjoint, and then use
their small boundary costs to check the quotient.

\paragraph*{The dense five-sets are disjoint}
Suppose two distinct dense five-sets $U,V$ overlap. Their combined
internal cost is greater than $4\kappa$, with pairs in $U\cap V$
counted twice. If they meet in one vertex, no pair is counted twice,
and degree counting on their nine-vertex union bounds this cost by
$9M/2<4\kappa$, a contradiction.

If they meet in $j\in\{2,3,4\}$ vertices, the degrees of the common
vertices account for all pairs touching the intersection, including
both copies of every pair inside it. This contributes at most $jM$.
The only uncounted pairs lie inside $U\setminus V$ or $V\setminus U$.
The resulting bounds are
\[
\begin{array}{c|c|c|c}
 j&|U\setminus V|&\text{cost inside each difference}&
                   t(U)+t(V)\text{ is at most}\\ \hline
 2&3&2&2M+4\\
 3&2&1&3M+2\\
 4&1&0&4M
\end{array}
\]
The first row uses the triangle hypothesis, and the second uses
the unit bound on each pair. Since $M<2$, all three bounds are at
most $2M+4=37/5<4\kappa$. Hence no two dense five-sets overlap.

\paragraph*{The remaining dense triangles are disjoint}
A dense five-set has boundary cost
\[
 \partial t(U)<5M-4\kappa=\frac{63}{125}<\frac M3.
\]
A triangle crossing its boundary contains two boundary pairs and
one other pair. The other pair costs at most one, so a dense crossing
triangle would need boundary cost greater than
$\kappa-1=999/1000>63/125$. This is impossible. A dense triangle
contained in a selected five-set needs no separate contraction.

It remains to compare dense triangles outside the five-sets.
Two sharing just one vertex have combined cost at most $M+2$:
the common vertex accounts for at most $M$, and the two opposite
pairs account for at most one each. Two sharing an edge have
combined cost at most $2M$, by summing degrees at the two common
vertices. Both bounds are at most $2\kappa$, contradicting density.
The selected triangles are therefore disjoint, and each has
\[
 \partial t(U)<3M-2\kappa\le\frac{2M}{3}.
\]
After contraction, the degree cost of a block is exactly its former
boundary cost. Ordinary vertices keep degree cost at most $M$;
three-blocks have degree cost below $2M/3$, and five-blocks below $M/3$.

\paragraph*{Triangles in the quotient}
We check the four possibilities according to the blocks present.
\begin{enumerate}[label=(\roman*)]
\item \emph{No block.} A dense triangle on three ordinary vertices
would already have been selected, so its cost is at most $\kappa$.
\item \emph{One five-block and two ordinary vertices.} The pairs
touching the block cost less than $M/3$ altogether, and the pair
between the ordinary vertices costs at most one. The total is
less than $M/3+1\le\kappa$.
\item \emph{One three-block and two ordinary vertices.} If this
quotient triangle had cost greater than $\kappa$, expand the block.
Its internal triangle also costs greater than $\kappa$, so the
resulting original five-set would have cost greater than $2\kappa$.
It would have been selected before the three-block, a contradiction.
\item \emph{At least two blocks.} Sum the three quotient degrees
and divide by two. Using $2M/3$ for each of two blocks and $M$ for
the third vertex gives cost at most $7M/6<\kappa$.
\end{enumerate}

\paragraph*{Larger odd sets in the quotient}
A five-set with no block was already tested in the original graph.
One containing a block has internal cost less than
\[
 \frac{2M/3+4M}{2}=\frac{7M}{3}<2\kappa.
\]
Here we used the weaker boundary bound $2M/3$, valid for either
kind of block. Finally, for odd $|S|\ge7$, every quotient degree
is at most $M$, so
\[
 t(S)\le\frac{|S|M}{2}
 =\frac{|S|}{|S|-1}M\frac{|S|-1}{2}
 \le\frac{7M}{6}\frac{|S|-1}{2}
 <\kappa\frac{|S|-1}{2}.
\]
The last step is \eqref{eq:subtwo-margins}. Singletons have no
internal pairs, which completes the odd-set check.
\end{proof}

\begin{figure}[htbp]
\centering
\begin{tikzpicture}[x=1cm,y=1cm,font=\small]
\draw[fill=softblue,draw=inkblue,dashed,rounded corners]
 (-.45,-.45) rectangle (3.2,2.2);
\node[vertex,label=below:$v_1$] (a) at (0,0) {};
\node[vertex,label=below:$v_2$] (b) at (2.7,0) {};
\node[vertex] (c) at (2.7,1.5) {};
\node[vertex] (d) at (1.35,1.85) {};
\node[vertex] (e) at (0,1.5) {};
\foreach \x/\y in {a/b,b/c,c/d,d/e,e/a,a/d,b/d,c/e}{
 \draw[black!35] (\x)--(\y);
}
\draw[very thick,warmred] (a)--(-1.4,-1.0);
\draw[very thick,inkblue] (b)--(3.8,-1.0);
\node[vertex] at (-1.4,-1.0) {};
\node[vertex] at (3.8,-1.0) {};
\node at (1.35,2.55) {A dense five-vertex set};
\draw[-{Latex},thick] (4.05,.85)--node[above] {contract}(5.45,.85);
\node[circle,draw=inkblue,fill=softblue,minimum size=10mm] (B) at (7,.85) {$B$};
\draw[very thick,warmred] (B)--node[above,sloped] {label $v_1$}(5.9,-1.0);
\draw[very thick,inkblue] (B)--node[above,sloped] {label $v_2$}(8.3,-1.0);
\node[vertex] at (5.9,-1.0) {};
\node[vertex] at (8.3,-1.0) {};
\node at (7,2.55) {The quotient vertex};
\end{tikzpicture}
\caption{Each external color at a contracted set records the original
vertex it uses. When restoring a color labeled $v_i$, the internal
pattern avoids $v_i$. The two external colors shown are distinct;
the grey internal edges have not yet been colored.}
\label{fig:reserved-endpoints}
\end{figure}
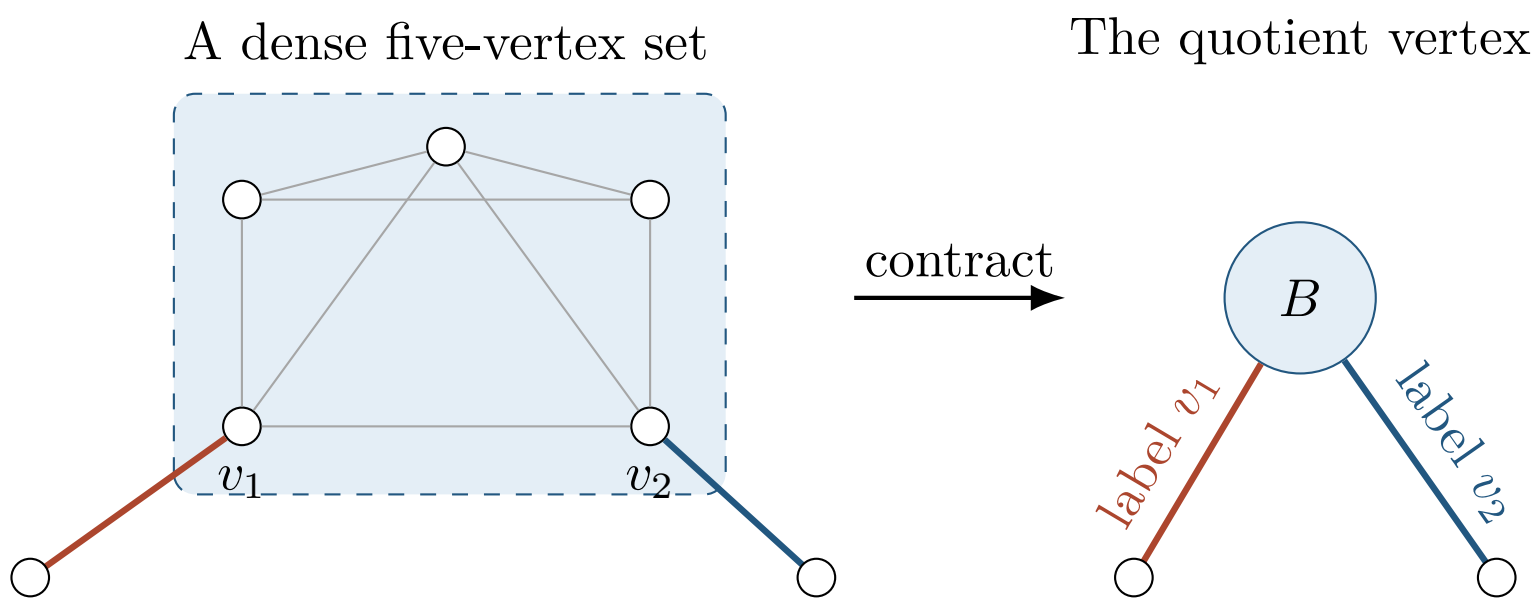

\subsubsection{Restoring a Dense Block}\label{sec:subtwo-restoration}
Work in the normalized local model, in which each star is an average of feasible integer star patterns and palettes are measured by their mass.
For a triangle, let $0\le t_i\le1$ be the cost of the bundle opposite vertex $i$, and suppose that $\kappa<T=\sum_i t_i\le2$.
The external colors labeled $i$ form a reserved class of mass $0\le q_i\le M-T+t_i$ in which only bundle $i$ may be used.
The restriction has a simple meaning: an external edge already uses
vertex $i$, but an internal packet on the opposite pair does not
touch $i$. We first place as much of that bundle as the reserved
class permits, then combine portions of the three bundles in colors
that have no external restriction.
Put $s=T-M>0$; the budget inequality implies $t_i\ge s>0$.
Since $q_i\le t_i-s\le(1-s)t_i$, a fraction $q_i/t_i$ of bundle $i$ fits in its reserved class while leaving at least an $s$ fraction of the bundle.
Here a fraction of a bundle means the same fraction of each of its
type demands and of its fractional packing. In particular, it has
the corresponding fraction of the bundle's packing cost.

For the common $s$ fraction of the three bundles, take an $s/2$ multiple of each of the three feasible star packings.
Each bundle then appears twice, so this uses palette mass $3s/2$ and covers exactly the required demand.
These star packings may use all three vertices and are assigned to
unreserved colors. The other $1-s$ fraction of each bundle is packed
separately, with the designated fraction already in its reserved
class. Thus the reserved classes are included in, rather than added
to, the separate-packing cost $(1-s)T$.
Packing the other fractions separately, including those placed in reserved classes, uses total palette mass at most
\[
T-sT+\frac32s\le1+\frac M2<\kappa.
\]
The quadratic expression increases on $T\le2$, since its
derivative is $M+5/2-2T\ge M-3/2>0$; its value at $T=2$ is
$1+M/2$.  Empty patterns fill unused colors.

For a five-block, let $q_v$ be its external endpoint budget and
$d_v=\sum_{e\ni v}t_e$.  A color labeled $v$ is available for
internal patterns avoiding $v$.
The next lemma measures how much this restriction can cost.
After normalizing the best unrestricted pattern to have price one,
the best pattern avoiding $v$ has price $r_v\le1$. Reserving $q_v$
colors for such patterns can therefore lose at most $q_v(1-r_v)$
in this price calculation. We must bound the internal demand price
and all these losses together.

\Needspace{13\baselineskip}
\begin{lemma}\label{lem:subtwo-five-palette}
Consider a demand on five vertices whose star at each vertex is the mean of feasible integer star patterns, as in \Cref{lem:subtwo-five-local}.
Suppose $T=\sum_e t_e\ge2\kappa$, every internal triangle has cost at most two, and the nonnegative external budgets satisfy $d_v+q_v\le M$.  For every nonnegative price
vector, let $C$ be the demand price when the maximum whole-pattern
price is normalized to one, and let $r_v$ be the maximum price of a
pattern avoiding $v$.  Then
\begin{equation}\label{eq:subtwo-five-palette}
 C+\sum_vq_v(1-r_v)\le\kappa-\frac1{1000}.
\end{equation}
Consequently, a palette of mass $\kappa$ fractionally packs the
internal demand when a disjoint class of mass $q_v$ is reserved for
patterns avoiding each $v$.
\end{lemma}
\begin{proof}
We introduce a sixth vertex $\infty$ to record the endpoint
restrictions using ordinary matchings. If a matching uses the edge
$\infty v$, all of its internal edges avoid $v$; its internal
packets can therefore share an external color labeled $v$.
The artificial edge will be charged the lost price $1-r_v$.
We give these edges larger amounts than their actual budgets, so that
the resulting matching inequality bounds the loss strongly enough.

\paragraph*{Choosing the artificial edge amounts}
Put $q=5M-2T\ge\sum_vq_v=:Q$ and
$\lambda=\max(2,T/2)$.  If $q=0$, all $q_v$ vanish and
\Cref{lem:subtwo-five-local} gives $C\le c<\kappa-1/1000$.
The inequality $Q\le q$ comes from summing $d_v+q_v\le M$ over
the five vertices, since $\sum_vd_v=2T$. Thus $q$ is an upper bound
on the total external budget, whether or not every vertex uses its
full allowance.

Otherwise set $k=1+(\lambda-M)/q>1$.  On the five internal vertices
and a new vertex $\infty$, assign amount $t_e$ to each internal pair
and $kq_v$ to $\infty v$. The choice of $k$ gives
$(k-1)q=\lambda-M$: it uses the gap between $M$ and $\lambda$
to increase all the external budgets by the same factor.

\paragraph*{Checking the matching inequalities}
We claim that these edge amounts can be covered by matchings of
total mass $\lambda$. After division by $\lambda$, Edmonds'
matching-polytope theorem requires degree at most one and internal
amount at most $(|S|-1)/2$ for every odd vertex set $S$.
On six vertices the only nontrivial odd-set sizes are three and
five, so the following list is exhaustive.

Each original vertex has degree
$d_v+kq_v\le M+(k-1)q_v\le M+(k-1)Q\le\lambda$.
At the new vertex the degree is
\[
kQ\le q+\lambda-M=\lambda+4M-2T\le\lambda-1.
\]
For the last inequality, $2T\ge4\kappa>4M+1$.
\begin{enumerate}[label=(\roman*)]
\item A triangle avoiding $\infty$ has cost at most $2\le\lambda$
by hypothesis.
\item A triangle containing $\infty$ has one internal pair, of
cost at most one: restrict a feasible local star distribution to
that pair to obtain a packing of mass one. Its two artificial edges have combined amount at
most $kQ$. Its total is at most $\lambda$.
\item The five-set avoiding $\infty$ is the whole original block.
Its cost is $T\le2\lambda$ by the choice of $\lambda$.
\item A five-set containing $\infty$ omits one original vertex,
say $v$. Its internal amount is $T-d_v+k(Q-q_v)$.
Summing the other four budget inequalities gives
\[
 Q-q_v\le4M-2T+d_v,
\]
because the sum of their internal degrees is $2T-d_v$.
Consequently its amount is at most
\begin{align*}
 T-d_v+k(Q-q_v)
 &\le T+k(4M-2T)+(k-1)d_v\\
 &\le T+k(5M-2T)-M\\
 &=\lambda+3M-T\le2\lambda.
\end{align*}
The second line uses $d_v\le M$ and $k>1$. For the equality,
$5M-2T=q$ and $kq=q+\lambda-M$; for the final inequality,
$T+\lambda\ge3T/2\ge3\kappa>3M$.
\end{enumerate}
Edmonds' theorem therefore supplies the claimed matching cover.

\paragraph*{Charging the reserved colors}
Let $\beta_e$ be the maximum internal packet price and put
$Z=\sum_vq_v(1-r_v)$.  Price edge $e$ by $\beta_e$ and edge
$\infty v$ by $1-r_v$. A matching without an artificial edge
has price at most one, since its internal packets are on disjoint
pairs and together form a feasible whole pattern. A matching using
$\infty v$ has internal packet price at most $r_v$, since all
those packets avoid $v$. Adding its artificial price $1-r_v$
again gives at most one.

The cover of mass $\lambda$ thus gives
\[
 \sum_e t_e\beta_e+kZ\le\lambda.
\]
A fractional packing of pair $e$ reproduces its type demand using
mass $t_e$, and each of its packets has price at most $\beta_e$.
Summing over pairs yields $C\le\sum_e t_e\beta_e$, so
$C+kZ\le\lambda$.
Combining this with $C\le c$ gives
\[
C+Z\le\left(1-\frac1k\right)C+\frac\lambda k
\le F:=c+\frac{\lambda-c}{k}.
\]
Both inputs are needed here: the five-vertex price lemma bounds
the unrestricted internal demand by $c$, while the matching cover
bounds the additional cost of respecting endpoint labels.

\paragraph*{The numerical margin}
For $2\kappa\le T\le4$, we have
$q\le63/125<3/5$ and $2-M=3/10$. Hence $k>3/2$ and
\[
 F<c+\frac23(2-c)=\frac{119}{60}<\kappa-\frac1{1000}.
\]
For $T\ge4$, write $\lambda=(5M-q)/4$ and
$k=(M+3q)/(4q)$. Then
\[
(\kappa-F)(M+3q)
=q^2-\frac{553}{1000}q+\frac{833}{10000}
=\left(q-\frac7{25}\right)^2
 +\frac{7q}{1000}+\frac{49}{10000}.
\]
Here $q\le1/2$ and $M+3q\le16/5$, so
$\kappa-F\ge49/32000>1/1000$.
Finally, for every nonnegative price vector the
reserved palette has maximum price
$\kappa-\sum_vq_v+\sum_vq_vr_v=\kappa-Z$.
Indeed, its unrestricted class has mass $\kappa-Q$ and can use
patterns of price at most one; the class labeled $v$ has mass $q_v$
and can use patterns of price at most $r_v$. These are exactly the
separate class constraints in the fractional packing program.
The inequality just proved says that the demand price $C$ is less
than this available price for every normalized nonnegative choice
of type prices.
The strict dual inequality proves fractional feasibility by finite
linear programming duality.  Nonnegative prices suffice because
deleting items preserves feasibility.
\end{proof}

\subsubsection{Proof of the Bound}\label{sec:subtwo-proof}
\begin{proof}[Proof of \Cref{thm:general-bp-subtwo}]
The empty graph is immediate. Put $n=\bp(G)$ and choose a sufficiently large integer $T_0\ge200$, with $b=1/T_0$ small enough for the desired $\varepsilon$.
Fix also $0<\sigma<1$, to be chosen sufficiently small after $b$.
Apply the grouping, rank rounding, and packet preparation of \Cref{lem:subtwo-preprocess}.  The retained
objects fit at every vertex in
\[
N=\left\lceil\frac{n}{1-2b}\right\rceil+1
\]
bins; each original pair has a bounded number of non-tiny types
and at most one tiny remainder object.  The set-asides have
$O(\sigma n/b)$ weighted degree.  On each original pair put
$t_{uv}=F(P_{uv})/N$.  Projection of local bins gives $t_{uv}\le1$;
the two-collection packing lemma gives triangle cost at most two.
The charge estimate of \Cref{lem:subtwo-charge} gives incident
pair cost below $M$.

For a block $B$ and $v\in B$, let $q_v$ be the charge of retained
external objects incident with $v$, divided by $N$.  Charge a tiny
external object by $w/(1-b)$.  In each local bin, keep each internal
neighbor projection together and external objects separately.  The
margin after \Cref{lem:subtwo-charge} gives
\begin{equation}\label{eq:subtwo-local-charge}
 \sum_{u\in B\setminus\{v\}}t_{uv}+q_v\le M.
\end{equation}
Select the blocks by \Cref{lem:subtwo-blocks}.  Form integer
packets separately on each original pair, then contract the blocks
once and delete internal packets from the quotient.  The quotient
packet multigraph has fractional ordinary chromatic index at most
$\kappa N$, by the matching-polytope formula and
\Cref{lem:subtwo-blocks}.  The ordinary multigraph coloring
theorem gives a proper coloring with $(\kappa+o(1))N$ colors,
uniformly in graph order. Pad the palette to at least $\lceil\kappa N\rceil$ colors, then add at most one further unused color to make its size even.  Residual multiplicity on a quotient pair is at most
$25(\lceil(D+1)/b\rceil+1)$, since two blocks contain at most
25 original pairs.  Packet count plus residual charge at an ordinary
quotient vertex is at most $MN$.  At a three-block it is at most
$(3M-2\kappa)N$, and at a five-block at most
$(5M-4\kappa)N$, by \eqref{eq:subtwo-local-charge} and their density.
Each of these constants is strictly below $\kappa$, supplying the
fixed missing-color margin needed in the next step.

Write $H$ for the properly colored quotient packet graph and $R$ for
the remaining external objects. A quotient vertex $B$ is an original
vertex or a contracted block of at most five vertices. Every packet
incidence at $B$ has its original endpoint $v\in B$ as a fixed label.
The packet coloring is proper, and a two-color path switch only
recolors entire packet edges, so each packet color at $B$ has exactly
one original endpoint label. Let $P_v$ be the number of packet colors
with label $v$, and let
$Q_{R,v}=\sum_{e\in R:\ell_B(e)=v}q_e$. The packet charge inequality,
applied to retained external objects at their original endpoint,
gives
\begin{equation}
 P_v+Q_{R,v}\le q_vN. \label{eq:subtwo-label-charge}
\end{equation}

Apply Lemma~\ref{lem:weighted-moment-mix} to $R$, using a fixed
rescaling of its size parameter so that the bounded-multiplicity
charges satisfy its hypotheses. Its fractions are supported only on
colors missing from $H$ at both endpoints, satisfy
$x_{e,c}=O(1/N)$ and $\sum_cx_{e,c}\le1$, and have the fixed charge
margin needed below. At $(B,c)$, let $I_{v,k}$ contain the residual
edges with endpoint label $v\in B$ and harmonic size class $k$,
where $k=0$ denotes the tiny class. There are $J=T_0+1$ classes.
We construct a random local set by choosing the pair $(v,k)$ with
joint probability $p_{v,k}$ and then sampling a subset conditional
on that pair. For $k\ge1$, set
\[
 p_{v,k}=\max\left\{k^{-1}\sum_{e\in I_{v,k}}x_{e,c},
                       \max_{e\in I_{v,k}}x_{e,c}\right\},
\]
and use the analogous tiny-class probability
\[
 p_{v,0}=\max\left\{\sum_{e\in I_{v,0}}q_ex_{e,c},
                       \max_{e\in I_{v,0}}x_{e,c}\right\}.
\]
Use probability zero for an empty class and the empty set, with the
null label, for the remaining probability. The displayed probabilities
sum to at most one, as verified below.
Conditional subset decompositions are those of
Lemmas~\ref{lem:bp-classchoice} and
\ref{lem:algorithmic-local-load}. Each chosen set is feasible,
contains only edges with label $v$, and the resulting distribution has inclusion probability
$x_{e,c}$ for every edge. The probability of label $v$ is bounded by
\begin{equation}
 \pi_{B,c,v}=\sum_{k=0}^{T_0}p_{v,k}
 \le \sum_{e:\ell_B(e)=v}q_ex_{e,c}
       +J\max_{e\in E_R(B)}x_{e,c}. \label{eq:subtwo-label-probability}
\end{equation}
The total class probability at $(B,c)$ is at most
$\sum_{e\in E_R(B)}q_ex_{e,c}+5J\max_e x_{e,c}<1$ for all large
$N$, because the missing-color margin is fixed. Empty probability
receives the null label. Summing \eqref{eq:subtwo-label-probability} over the $O(N)$ palette colors
and using $\sum_cx_{e,c}\le1$ yields
$\Pi_{B,v}\le Q_{R,v}+O(1)$, uniformly in graph order.

Apply Lemma~\ref{lem:bp-labelled-round} to $R$ with these distributions, using size parameter $5N$, at most five labels, and accuracy parameters divided by five. Indeed, each contracted block contains at most five original vertices, so $W_R(B)\le5N$; the multiplicity is bounded and the other constants change only by fixed factors. It gives a partial residual coloring in the old
palette, using only missing colors, such that every color at $B$ has
at most one original endpoint label and
\begin{equation}
 \#\{\text{old-palette external colors at $B$ labeled $v$}\}
 \le P_v+Q_{R,v}+\eta N+O(N^{3/4})
 \le q_vN+\eta N+O(N^{3/4}). \label{eq:subtwo-label-count}
\end{equation}
The packet colors and residual colors in this count are disjoint,
because residual support consists of packet-missing colors.
Uncolored residual objects, internal remainders, and rank set-asides
are colored later with fresh colors disjoint from the old palette;
those fresh colors are not reserved for internal patterns. Thus \eqref{eq:subtwo-label-count}
is exactly the count needed to restore each block. Choosing $\eta$
small relative to the fixed palette and price margins gives the
claimed $q_vN+\delta N+o_\delta(N)$ bound after renaming the accuracy parameter.

Restore the disjoint blocks using the original full costs and the reserved label classes.
After scaling by $N$, the fractional construction for a triangle uses palette mass at most $(1+M/2)N<\kappa N$; a five-block uses palette mass $\kappa N$ by \Cref{lem:subtwo-five-palette}.
For either kind of block, keep the external coloring fixed. If its actual class labeled $v$ has $m_v$ colors, reserve exactly
$r_v^*=\min\{m_v,\lfloor q_vN\rfloor\}$ of them and apply the palette argument with the smaller budget $q_v'=r_v^*/N\le q_v$.
Reducing a budget preserves its hypotheses and only increases the available internal patterns. When constructing the internal packing, temporarily treat the omitted colors as unreserved.
If an internal pattern is assigned to one of these colors, move that pattern to a fresh color before combining the internal and external colorings.
The label excess therefore costs at most $5\delta N+o_\delta(N)+O(1)$ fresh colors per block.
These fresh colors are used only on internal edges, so the same palette can be reused by all disjoint blocks.

It remains to turn the fractional internal packings into assignments of the original objects.
First omit each block's at most ten tiny remainder objects for cleanup.  Extend non-tiny prices by zero when
applying the full-demand inequality; deleting tiny items changes
neither the maximum pattern price nor the maximum price avoiding a
specified vertex. The density condition is used for the original full
demand, before the tiny objects are omitted.
There are only accuracy-dependent boundedly many remaining types and
patterns per block.  Floor the masses in a basic feasible internal
pattern solution and fill its slots with distinct actual objects.
Only boundedly many internal objects per original vertex remain.
Internal patterns in a color labeled $v$ avoid $v$, so they are
compatible with all external objects of that color.

We color the remaining external objects, the bounded internal remainders, and the objects set aside during rank grouping with fresh palettes.  The greedy cleanup
lemma bounds their total cost by an arbitrarily small multiple of
$n$ plus an accuracy-dependent constant, uniformly in graph order.
Choose $b$, then the rank parameter $\sigma$, then $\delta$ and the remaining fixed accuracies sufficiently small in terms of $\varepsilon$.  Since
$N=n/(1-2b)+O(1)$, the total is
$(\kappa+\varepsilon)n+C_\varepsilon$.  Bounded $n$ is covered by
the same greedy bound.  Restoring the original weights and expanding packets and grouped objects can only decrease color loads.
\end{proof}

\Needspace{24\baselineskip}
\begin{algorithm}[The general $1.999$ construction]\label{alg:general-bp-subtwo}
Input: a weighted multigraph without loops and a desired accuracy.
Choose the parameters in the order specified in the proof.
\begin{enumerate}[leftmargin=*,nosep]
\item Find approximate local packings as in
\Cref{sec:local-input-algorithms}. Group and round each original
parallel bundle as in \Cref{lem:subtwo-preprocess}.
\item Compute the fractional bundle costs and use them to select the
dense five-sets and triangles. Form whole packets on each original
pair, then contract the selected blocks and remove their internal packets.
\item Color the quotient packets. Insert residual external objects
while recording one original endpoint per block and color, using the
retained-list sampler of \Cref{sec:label-algorithm}.
\item Set aside each block's tiny internal remainder objects for cleanup.
Solve the bounded internal pattern programs in the reserved palettes
and round their pattern masses down. Color the internal remainders,
including the tiny objects, and excess-label patterns with fresh colors reused across
disjoint blocks; greedily color the remaining external objects and
preprocessing set-asides with another fresh palette.
\item Expand every packet and grouped object. Check each vertex-color
load and return the resulting coloring.
\end{enumerate}
Every step has
expected polynomial running time for each fixed accuracy.
\end{algorithm}

\subsection{Keeping Endpoint Labels in Polynomial Time}
\label{sec:label-algorithm}

The extra requirement in \Cref{lem:bp-labelled-round} is a count test
on a vertex's own list. We can therefore use the sampler in \Cref{sec:computable-progress}
with a few more rejection tests; we never enumerate complete lists.
The rational class-choice and bounded-load decompositions give local
distributions of polynomial support. Conditioning on an endpoint
label, compressing, and smoothing them gives a bounded number of
decorated profiles with rational probabilities of polynomial encoding
length, for fixed parameters.

For each vertex $z$ and label $a$, compute the rational mean
\[
 \mu_{z,a}=\sum_{c\in\mathcal C}
       \sum_{P:\,\operatorname{label}(P)=a}\lambda'_{z,c,P}.
\]
Draw the profiles at different colors independently. Accept the list
if its usual progress test passes and its label counts satisfy
$C_{z,a}\le\mu_{z,a}+t$ for every $a$, where
$t=\lceil n^{3/4}\rceil$. Computing the means, counting labels,
and evaluating progress all take polynomial time with rational
arithmetic. By \eqref{eq:label-rejection-bound}, the rejection
probability is at most the computable rational number $p'_k=o(1)$.
Thus, for sufficiently large $n$, an accepted list takes at most two
draws in expectation.

No two vertices in the current group are adjacent. Acceptance at one
vertex depends only on its own list and the earlier groups, so these
samplers give independent lists with exactly the required conditioned
distributions. The Harris--Srinivasan procedure in \Cref{sec:computable-progress} needs
only these independent samples and the coefficients of the current
assignment. It does not need the acceptance probabilities or the
normalizing constants of the conditioned distributions. In its mean
bounds, replace $1-p_k$ by $1-p'_k$: a load row uses
$m_k/(1-p'_k)$, and an error row with unchanged contribution $C$ uses
$(Q_v^{(k-1)}-C+1/n)/(1-p'_k)$, with the same scaling as before.
The extra $O(n^{-1/2})$ in $p'_k$ is smaller than the existing
$O(n^{-1/4}\sqrt{\log n})$ bound on $p_k$, so the load and progress
estimates keep their stated orders. The row tests and subset sampler
are unchanged. This proves the expected polynomial running time in
\Cref{lem:bp-labelled-round}. The later option matching and small-edge
rounding preserve labels because they only remove options.

\paragraph*{The other steps of the $1.999$ construction}
Use \Cref{sec:local-input-algorithms} to find local packings and let
$B$ be their maximum bin count, padding with empty bins. For each
fixed $\xi>0$,
\[
 B\le(1+\xi)\bp(G)+O_\xi(1).
\]
Run the construction with
$N=\lceil B/(1-2b)\rceil+1$. Its estimates require actual local
packings with this bin bound, so this substitution preserves every
charge and density inequality. Choosing $\xi$ and $b$ small enough
absorbs the factor $(1+\xi)/(1-2b)$ into the requested accuracy.
In particular, $N$ is at most a constant times the number of input
edges plus an accuracy-dependent constant. Enumerating all triangles
and five-sets takes $O(|V(G)|^5)$ tests. The bundle programs, packet
coloring, and missing-color mixing have the implementations already
given. The preceding labeled sampler handles the residual external
objects.

For internal restoration, first remove the at most ten tiny objects
in a block, as in the proof. At most $10D$ types remain, all of
weight at least $b$. A feasible internal pattern has at most
$\lfloor5/(2b)\rfloor$ objects: summing its vertex-load bounds counts
each object twice. Thus all its integer type vectors can be listed
in a number depending only on accuracy. There are at most six palette classes:
one unrestricted class and at most five endpoint classes. List the
allowed patterns in each class and solve the resulting constant-size
rational linear program. The palette argument proves it feasible.
Flooring the class-pattern masses uses no more than the available
integer number of colors in any class and leaves only boundedly many
objects per block. Assign actual objects to the remaining slots and
perform the stated cleanup and expansion. All numbers have polynomial
encoding length, and every operation is polynomial for fixed accuracy.
Together these steps give the running-time claim in
\Cref{thm:general-bp-subtwo}.

\end{document}